\documentclass[12pt,a4paper]{article}
\usepackage{hyperref}
\usepackage{amsthm,amsmath,amssymb}
\usepackage{verbatim}
\usepackage{graphicx}
\usepackage{enumitem}
\usepackage{abstract}
\usepackage{url}
\usepackage{xcolor}
\usepackage{float}
\usepackage[T1]{fontenc}
\usepackage{cite}

\theoremstyle{plain}
\newtheorem{theorem}{Theorem}[section]
\newtheorem{lemma}[theorem]{Lemma}
\newtheorem{proposition}[theorem]{Proposition}
\newtheorem{corollary}[theorem]{Corollary}

\theoremstyle{definition}
\newtheorem{definition}[theorem]{Definition}

\theoremstyle{remark}

\numberwithin{equation}{section}

\def\IN{{\mathbb N}}
\def\IZ{{\mathbb Z}}
\def\ID{{\mathbb D}}

\def\IR{{\mathbb R}}

\def\dom{\mathrm{dom}}
\def\range{\mathrm{range}}
\def\prefixeq{\subseteq} 
\def\truepath{\mathcal{T}}
\newcommand{\dotminus}{\buildrel\textstyle.\over{\hbox{\vrule height3pt depth0pt width0pt}{\smash-}}}
\def\sign{\mathrm{sign}}
\def\restrictedto{\upharpoonright}
\newcommand{\leqS}{\leq_{\mathrm{S}}}

\newcommand{\equivS}{\equiv_{\mathrm{S}}}
\newcommand{\modd}{\,\mathrm{mod}\,}
\def\EC{\IR_{\mathrm{c}}}
\def\LC{\IR_{\mathrm{lc}}}
\def\NC{\IR_{\mathrm{nc}}}
\def\WC{\IR_{\mathrm{wc}}}
\def\WNC{\IR_{\mathrm{wnc}}}
\def\LNC{\IR_{\mathrm{lnc}}}
\def\DLNC{\IR_{\mathrm{dlnc}}}

\begin{document}

\title{Nearly Computable Numbers with Additional Properties}

\author{Peter Hertling and Philip Janicki\\[2mm]
Fakult\"at f\"ur Informatik,
Universit\"at der Bundeswehr M\"unchen, \\
85579 Neubiberg, Germany \\[2mm]
{\tt peter.hertling@unibw.de, philip.janicki@unibw.de}}

\date{\today}

\maketitle


\begin{abstract}
A real number is called \emph{nearly computable} if it is the limit of a computable sequence of rational (or dyadic) numbers
that has the additional property that
for any computable subsequence the sequence of differences of consecutive elements converges to $0$ computably.
The main result of this article provides methods for constructing nearly computable numbers with additional properties. 
The main result is proved by an infinite injury priority argument.
Using these methods, we construct nearly computable numbers with several different additional properties:
\begin{itemize}
\item
a nearly computable number that is not divergence bounded computable,
\item
for any left-computable number $\alpha$ that is neither in the least nor in the highest Solovay degree,
a left-computable, nearly computable number that is Solovay incomparable with $\alpha$,
\item
a nearly computable number $\gamma$ that can be written as the difference of two left-computable numbers
but has the property that for any two left-computable numbers $\alpha$ and $\beta$ with
$\gamma=\alpha-\beta$ both $\alpha$ and $\beta$ are Martin-L{\"o}f random.
\end{itemize}
We also look at two classes of nearly computable numbers:
the class $\DLNC$ of differences of real numbers that are nearly computable and left-computable,
and the class $\WNC$ of nearly computable numbers that can be written as the difference of two left-computable numbers.
Both classes form real closed fields,
in fact, $\DLNC$ is closed under computable functions with open domain that are Lipschitz continuous or differentiable.
Finally, we show that $\DLNC$ is a proper subset of $\WNC$.

\end{abstract}

\section{Introduction}
\label{section:intro}

A real number is called \emph{nearly computable} if it is the limit of a computable sequence $(x_n)_n$ of dyadic numbers 
that has the additional property that for any computable, strictly increasing function $s:\IN\to\IN$ the sequence
$(x_{s(n+1)}-x_{s(n)})_n$ converges to $0$ computably~\cite[Definition 6.1]{HJ2024}.
This number class ---let us call it $\NC$--- has many natural properties; compare~\cite{HJ2024}.
The main result of this article provides methods for constructing nearly computable numbers with additional properties. 
The main result is proved by an infinite injury priority argument.
It generalizes an old construction by Downey and LaForte~\cite[Theorem 8]{DL2002}
and a recent construction by H{\"o}lzl and Janicki~\cite[Theorem 37]{HJ2026}.
Using the main result, we construct nearly computable numbers with several different additional properties:
\begin{itemize}
\item
a nearly computable number that is not divergence bounded computable;
this strengthens a result by Zheng, Rettinger, Gengler~\cite[Proposition 5.2(3)]{ZRG2005}.
\item
for any left-computable number $\alpha$ 
(a real number is \emph{left-computable} if it is the limit of a computable, increasing sequence of dyadic numbers)
that is neither in the least nor in the highest Solovay degree, 
that means, it is neither computable nor Martin-L{\"o}f random,
a left-computable, nearly computable number that is Solovay incomparable with $\alpha$,
\item
a nearly computable number $\gamma$ that is weakly computable
(a real number is \emph{weakly computable} if it is the difference of two left-computable real numbers)
but has the property that for any two left-computable numbers $\alpha$ and $\beta$ with
$\gamma=\alpha-\beta$ both $\alpha$ and $\beta$ are Martin-L{\"o}f random;
this strengthens a result by Miller~\cite[Theorem 3.5]{Mil2017}.
\end{itemize}
In ~\cite{HJ2024} the authors have shown that the nearly computable numbers form a real closed field.
Here we investigate the real numbers that are nearly computable and weakly computable.
These numbers form a real closed field as well, called $\WNC$.
Furthermore, we also consider the set $\DLNC$ of differences of real numbers that are nearly computable and left-computable.
Obviously, it is a subset of $\WNC$. We show that it is a real closed field as well.
In fact, we show that $\DLNC$ is closed under computable functions with open domain that are Lipschitz continuous or differentiable.
Furthermore, we show that it is a proper subset of $\WNC$.
Thus, there exists a nearly computable real number that can be written as the difference of two left-computable numbers
but not as the difference of two nearly-computable, left-computable numbers.
In fact, this follows from the third of the three results listed above.

In the following section we introduce basic notions and formulate basic results that will be needed later on.
Section~\ref{section:master} contains the main theorem, and in Section~\ref{section:proof-master} we prove it.
As an example of how to apply it, in Section~\ref{section:not-right-computable} we use the theorem in order to show the known result that there exists a left-computable and nearly computable number that is not computable.
This follows from the main result by Downey and LaForte~\cite[Theorem 8]{DL2002}
and an observation by Stephan and Wu~\cite[Theorem 5]{SW2005}; see~\cite[Theorem 7.3]{HJ2024}.
In Section~\ref{section:divergence} we deduce the result saying that there exists a nearly computable number that is 
not divergence bounded computable. 
As the remaining two new results that we wish to deduce from the main result
make use of Solovay reducibility, in Section~\ref{section:Solovay-reducibility}
we formulate several characterizations of this reducibility notion and present some fundamental, well-known facts about it.
In Section~\ref{section:Solovay-incomparable} we show that for any left-computable number $\alpha$ that is neither in the smallest
nor in the largest Solovay degree there exists a nearly computable number that is Solovay incomparable with $\alpha$.
In Section~\ref{section:not-regainingly-approximable} we demonstrate that a recent result by H{\"o}lzl and Janicki~\cite[Theorem 37]{HJ2026},
namely the existence of a left-computable number which is nearly computable
but not regainingly approximable, can be deduced from our main theorem.
Section~\ref{section:decomp-random} contains our strengthening of Miller's result stated above.
Finally, in Section~\ref{section:wnc-dlnc} we discuss the sets of numbers $\DLNC$ and $\WNC$.
We show that they are real closed fields, and from the result in the previous section we deduce that the first is a proper subset of the second.

\section{Basic Notions and Observations}
\label{section:basic}

Let $X$ and $Y$ be sets.
By $f:\subseteq X\to Y$ we denote a function defined on a subset $\dom(f):=\{x\in X~:~(\exists y \in Y) (x,y) \in f\}$ of $X$, the \emph{domain} of $f$,
with values in $Y$.
For $x\in X$ we write $f(x)\downarrow$ if $x\in\dom(f)$ and $f(x)\uparrow$ if $x\not\in\dom(f)$.
If $\dom(f)=X$ then we may write the function $f$ as $f:X\to Y$ and call it a {\em total} function.
Let $Y^X$ be the set of all total functions from $X$ to $Y$.
Note that in the case $X=\emptyset$ there is exactly one function in $Y^X$, the empty function.
By $\IN=\{0,1,2,\ldots\}$ we denote the set of natural numbers.
For $k\in\IN$ let $X^k :=X^{\{0,\ldots,k-1\}}$ be the set of all \emph{sequences of length $k$} of elements of $X$.
Sometimes we call them \emph{strings} or \emph{vectors}.
For $k=0$ this set contains exactly one element, the empty sequence.
Let $X^*:=\bigcup_{k=0}^\infty X^k$ be the set of all finite sequences of elements of $X$.
The length of a finite sequence $\sigma\in X^*$ will be written as $|\sigma|$.
The elements of $X^\IN$ are called {\em sequences} or \emph{infinite sequences} (if we wish to stress that we are not considering a finite sequence)
with elements in $X$, and a sequence $p:\IN\to X$
is often written in the form $(p(n))_n$ or $(p_n)_n$.
For any $k\in\IN$ and $p\in X^\IN$ let $p\upharpoonright k$ be the prefix of length $k$
of the sequence $p:\IN\to X$,
that is, the uniquely determined element $g\in X^k$ with $g\subseteq p$, in other words,
$p\restrictedto k=p(0)\ldots p(k-1) = p_0\ldots p_{k-1}$.
Similarly, for any subset $A\subseteq\IN$ and any $k\in\IN$ we define $A\restrictedto k:=A \cap\{0,\ldots,k-1\} = \{n\in A\colon n <k\}$. 

Let $\IR$ be the set of all real numbers.
We call a finite or infinite sequence $a:\dom(a)\subseteq\IN\to\IR$ of real numbers
\emph{increasing} if, for all $n\in\IN$ with $n+1\in\dom(a)$, $a_n < a_{n+1}$, and 
\emph{nondecreasing} if, for all $n\in\IN$ with $n+1\in\dom(a)$, $a_n \leq a_{n+1}$, and we define the terms
\emph{decreasing} and \emph{nonincreasing} analoguosly.
A sequence $(x_n)_n$ of real numbers \emph{converges to a} real number $\xi$ if 
there exists a total function $g:\IN\to\IN$ such that, 
for all $m,n\in\IN$
\[ \text{if } m\geq g(n) \text{ then } d\left( x_m, \xi \right) \leq 2^{-n} . \]
Then $\xi$ is called a \emph{limit of $(x_n)_n$}, and we write $\lim_{i\to\infty} x_i = \xi$,
and $g$ is called a \emph{modulus of convergence (to $\xi$) of the sequence $(x_n)_n$}.
A sequence is said to \emph{converge computably (to $\xi$)} if it converges (to $\xi$)
and possesses a computable modulus of convergence.
Sometimes it is more convenient to work with the Cauchy condition.
We say that a sequence $(x_n)_n$ of real numbers is \emph{computably Cauchy}
if there exists a computable function $g:\IN\to\IN$ such that for all $n,k_1,k_2$ with $g(n)\leq k_1 \leq k_2$ we have
$\left| x_{k_1} - x_{k_2} \right|\leq 2^{-n}$
(one can assume without loss of generality that $g$ is increasing).
Then we will say that such a function \emph{witnesses} that the sequence $(x_n)_n$ is computably Cauchy.
Of course, a sequence $(x_n)_n$ is computably Cauchy if, and only if, it converges computably.
A function $s:\IN\to\IN$ {\em tends to $\infty$} if for every $M\in\IN$ there exists some $N\in\IN$ such that
for all $n\in\IN$,
\[ \text{if } n \geq N \text{ then } s(n) \geq M . \]


\begin{definition}[{\cite[Definition 5.1]{HJ2024}}]
We say that a sequence $(x_n)_n$ of real numbers is \emph{nearly computably Cauchy}
if for every total, computable, and increasing function $s:\IN\to\IN$
the sequence of real numbers $(x_{s(n+1)}-x_{s(n)})_n$ converges computably to $0$.
\end{definition}

The function $\langle \cdot,\cdot\rangle:\IN^2\to\IN$ defined by
$\langle i,j\rangle:= \frac{(i+j)\cdot (i+j+1)}{2} + i$, for all $i,j\in\IN$,
is a computable bijection between $\IN^2$ and $\IN$.
Also its projections $\pi_1:\IN\to\IN$ and $\pi_2:\IN\to\IN$ defined by
$\pi_1(\langle i,j\rangle):=i$ and $\pi_2(\langle i,j\rangle):=j$, for all $i,j\in\IN$, are computable.
%
Let 
\[ \ID:=\left\{x\in\IR\colon (\exists i,j,k\in\IN) \, x=\frac{i-j}{2^k} \right\} \]
be the set of all \emph{dyadic} numbers.
In order to transfer computability notions from $\IN$ to the set $\ID$ we use
the following surjective function $\nu_\ID:\IN\to\ID$ defined by
\[ \nu_\ID( \langle a,b \rangle) :=\begin{cases} 
	\frac{-a/2}{2^b} & \text{if } a \text{ is even}, \\
	\frac{(a+1)/2}{2^b} & \text{if } a \text{ is odd},
\end{cases} \]
for $a,b\in\IN$.
Instead of $\nu_\ID$ one may use any equivalent numbering of $\ID$; see~\cite[Section 4]{HJ2024}.
For introducing useful computability notions on $\IR$ even any equivalent numbering
of a dense subset of $\IR$ would work as well; 
see again~\cite[Section 4]{HJ2024}.
A sequence $(x_n)_n$ of dyadic numbers is called {\em computable}
if there exists a computable function $f\colon\IN\to\IN$ such that $x_n=\nu_\ID(f(n))$, for all $n$.

\begin{definition}
A real number $\xi$ is called \emph{computable} (resp. \emph{left-computable})
if there exists a computable sequence $(x_n)_n$ of dyadic numbers that converges to $\xi$ and
converges computably (resp. is nondecreasing).
\end{definition}

It is easy to see and well known that a real number $\xi$ is computable if, and only if, both $\xi$ and $-\xi$ are left-computable.

\begin{definition}[Ambos-Spies, Weihrauch and Zheng \cite{AWZ2000}]
\label{definition:wc}
A real number $\gamma$ is called \emph{weakly computable} if 
there exists a computable sequence $(c_n)_n$ of dyadic numbers that converges to $\gamma$ and has the property
that the series $\sum_n |c_{n+1} - c_n|$ converges.
\end{definition}

\begin{proposition}[Ambos-Spies, Weihrauch and Zheng \cite{AWZ2000}]
\label{proposition:wc}
A real number is weakly computable if, and only if, it can be written as the difference of two left-computable numbers.
\end{proposition}

\begin{lemma}[{\cite[Prop.~5.5]{HJ2024}}]
\label{lem:computable-convergent}
Let $(x_n)_n$ be a computable sequence of dyadic numbers that is nearly computably Cauchy.
Then $(x_n)_n$ converges.
\end{lemma} 

\begin{definition}[{\cite[Definition 6.1]{HJ2024}}]
A real number $\xi$ is called \emph{nearly computable}
if there exists a computable sequence $(x_n)_n$ of dyadic numbers 
that converges to $\xi$ and is nearly computably Cauchy.
\end{definition}

It is not difficult to show that every computably convergent sequence is nearly computably Cauchy~\cite[Cor.~3.5(1) and Lemma 5.2]{HJ2024}.
So we have the implication that every computable number is nearly computable.
But the converse is not true; compare Corollary~\ref{corollary:nearly-left-not-right}.

\begin{lemma}[{\cite[Prop.~7.2]{HJ2024}}]
A real number $\xi$ is left-computable and nearly computable
if, and only if, there exists a computable sequence $(x_n)_n$ of dyadic numbers 
that is nondecreasing, nearly computably Cauchy, and converges to $\xi$.
\end{lemma}

Let $\varphi_0,\varphi_1,\varphi_2,\ldots$ be a standard enumeration
of the set $P^{(1)}$ of all possibly partial computable functions with domain and range in $\IN$.
As usual, we write $\varphi_i(n)[t]{\downarrow}$ to express that the $i$-th Turing machine
(which computes $\varphi_i$)
stops after at most~$t$~steps on input $n$.
For $i,t\in\IN$ let
\[ \tilde{\ell}_i(t) := \max\{l\in\IN \colon l\leq t \text{ and } (\forall k < l) \, \varphi_i(k)[t]\downarrow \} .\]
Note that, for all $i$, the function $\tilde{\ell}_i:\IN\to\IN$ is nondecreasing and satisfies, for all $t$,
$\tilde{\ell}_i(t)\leq t$.
The function $\tilde{\ell}_i$ tends to infinity if, and only if, $\varphi_i$ is a total function.
The function $(i,t)\mapsto \tilde{\ell}_i(t)$ is computable.
For $i,t\in\IN$ let
\[ \ell_i(t) := \max\{l\in\IN \colon \begin{array}[t]{l}
    (\forall k < l) \, (\varphi_i(k)[t]\downarrow \text{ and } \varphi_i(k)< t) \text{ and } \\
    (\forall k) \, (0 < k < l \Rightarrow \varphi_i(k-1)< \varphi_i(k)) \} .
    \end{array}  \]
Note that, for all $i$, the function $\ell_i:\IN\to\IN$ is nondecreasing and satisfies, for all $t$,
$\ell_i(t)\leq t$ and,
if $\ell_i(t)>0$, then $\varphi_i(\ell_i(t)-1)<t$.
The function $\ell_i$ tends to infinity if, and only if, $\varphi_i$ is a total and increasing function.
The function $(i,t)\mapsto \ell_i(t)$ is computable.

\section{Main result}
\label{section:master}

Let $(\mathcal{R}_{2i}(\xi))_i$ be an infinite sequence of assertions about a real number $\xi$;
we shall call them \emph{requirements}.
In the applications of the main result these requirements will be chosen
such that any real number $\xi$ satisfying all of them has the property that we wish it to have.
In addition we want $\xi$ to be nearly computable.
This will be ensured by an infinite sequence of requirements $(\mathcal{R}_{2i+1}(\xi))_i$ with odd indices, which will be formulated in the
beginning of the following section.
In order to obtain such a number $\xi$, we will construct
a computable and converging sequence $(x_t)_t$ of dyadic numbers
such that its limit $\xi$ satisfies all of these requirements and is nearly computable.
In the construction we shall use a computable function $A:\IN^3\times\ID\to\ID$, which will depend upon the application.
In order to make sure that the construction works,
we require the function $A$ to have certain properties that will be described below.
In order to give the reader some intuition concerning $A$,
let us first explain the purpose of the input parameters of $A(d,t,\ell,q)$ and of its output.
\begin{itemize}
\item
The natural number $d$ is a precision parameter.
Note that any precision parameter $d$ will be used by at most one requirement $\mathcal{R}_{2i}(\xi)$
and the index $2i$ of this requirement is determined by the precision parameter $d$ via a computable function.
In fact, for later purposes, let $b:\{0,1\}^{\mathrm{even}}\times\IN\to\IN$
(where $\{0,1\}^\mathrm{even}$ is the set of all finite binary strings of even length)
be a bijective computable function such that, for all $m,n\in\IN$ and $\sigma,\tau\in\{0,1\}^{\mathrm{even}}$,
if $m\leq n$ and $\sigma \prefixeq\tau$ then $b(\sigma,m)\leq b(\tau,n)$.
The function $\iota\colon\IN\to\IN$ is defined by $\iota(b(\sigma,n)):=|\sigma|/2$, 
for all $\sigma\in\Sigma^{\mathrm{even}}$ and $n\in\IN$.
Note that it is computable.
Then for any possible precision parameter $d$ 
the requirement $\mathcal{R}_{2\iota(d)}(\xi)$ is the unique requirement
that can use this precision parameter.
\item
The natural number $t$ will be the current stage in the construction of a dyadic sequence $(x_s)_s$. 
\item
The natural number $\ell$ will be the number of stages in the construction of a dyadic sequence $(x_s)_s$
in which the function $A$ ``has been applied'' already with the parameter $d$ and in which some action 
has been required (this will be made precise below).
\item
The dyadic number $q$ will be the current value $x_{t-1}$ of the dyadic sequence $(x_s)_s$.
\item
The output $A(d,t,\ell,q)$ serves two purposes: to indicate whether a certain condition is satisfied
(this is a Boolean value; we encode this by making sure that $A(d,t,\ell,q)\neq 0$ iff this is true)
and to say which amount should be added to $q$ (this is the numerical value $A(d,t,\ell,q)$)
if this condition is satisfied. 
There may be cases where the respective condition is satisfied
but it is not really necessary to add something to $q$.
In these cases we will nevertheless demand that $A(d,t,\ell,q)$ is a small nonzero value,
that means, we will demand that, instead of not changing $q$ at all,
some small nonzero value should be added to $q$.
In the examples that we will consider, this is possible.
This allows us to encode the Boolean value (whether the condition is satisfied)
and the numerical value (if the condition is satisfied, then which amount should be added to $q$)
by just one value, the value $A(d,t,\ell,q)$.
\end{itemize}

We shall demand that the function $A$ satisfies several properties.
First of all,
\begin{equation}
\label{equation:range-powers-of-2}
\mathrm{range}(A) \subseteq \{0\} \cup\{-2^{-n}\mid n\in\IN\} \cup \{2^{-n} \mid n\in\IN\} . 
\end{equation}
In order to describe the other properties that we demand of $A$, we need to introduce some more notation.
Let us fix a sequence $(x_t)_t$ of dyadic numbers,
and let us fix a nonempty, finite or infinite, increasing sequence $(T_j)_j$ of natural numbers with $T_0>0$. 
For $d\in\IN$ we define a function $h_d\colon\dom(T)\to\IN$ recursively by
\begin{align*}
  h_d(0) &:=0, \\ 
	h_d(j) &:= \#\{k<j \colon A(d,T_k,h_d(k),x_{T_k-1})\neq 0 \}, \text{ for } j\in\dom(T) \text{ with } j>0 .
\end{align*}
The value $h_d(j)$ says in how many stages $T_k<T_j$ some action was required with the parameter $d$.
Obviously, the function $h_d$ is nondecreasing.

We shall consider two types of requirements $\mathcal{R}_{2i}(\xi)$:
on the one hand those that need to be taken care of only finitely often until they are satisfied forever,
and on the other hand those that need to be taken care of infinitely often.
We shall always assume that the requirements $\mathcal{R}_{2i}(\xi)$ with even $i$ are the ones
that need to be taken care of only finitely often until they are satisfied forever
while the requirements $\mathcal{R}_{2i}(\xi)$ with odd $i$ are the ones that need to be taken care of infinitely often.

Given a sequence $(x_t)_t$ of dyadic numbers and a number $d\in\IN$,
we say that a nonempty, finite or infinite, increasing sequence $(T_j)_j$ with $T_0>0$
is \emph{$(d,j,(x_t)_t)$-good} for some $j\in\dom(T)$
with $A(d,T_j,h_d(j),x_{T_j-1})\neq 0$ if the following is true:
for all $j'\in\dom(T)$ with $h_d(j')=h_d(j)+1$
\begin{equation}
\label{equation:i-d-good}
   \left| x_{T_{j'}-1} - x_{T_j-1} - A(d,T_j,h_d(j),x_{T_j-1}) \right|
	   \leq \frac{1}{2} \cdot \left|A(d,T_j,h_d(j),x_{T_j-1})\right| . 
\end{equation}
This condition means that the movement of $x_{T_j-1}$ by the amount of $A(d,T_j,h_d(j),x_{T_j-1})\neq 0$, 
that has been demanded in stage $T_j$,
has been carried out approximately until the next stage $T_{j+1}$ of this sequence,
and, as long as with the parameter $d$ no new movement is demanded, the value $x_{T_{j'}-1}$ will approximately stay there.
Given a sequence $(x_t)_t$ of dyadic numbers and a number $d\in\IN$,
we say that a nonempty, finite or infinite, increasing sequence $(T_j)_j$ with $T_0>0$
is \emph{$(d,(x_t)_t)$-good} if 
the sequence $(T_j)_j$ is $(d,j,(x_t)_t)$-good for all $j$ in the set
\[ \{j\in\dom(T) \colon A(d,T_j,h_d(j),x_{T_j-1})\neq 0 \} \setminus \mathrm{First}( d,(T_k)_k ) \]
where the set $\mathrm{First}( d,(T_j)_j )$ is defined by
\begin{eqnarray*}
  \lefteqn{\mathrm{First}( d,(T_j)_j )} && \\
	&:= & 
   \begin{cases}
	  \ \emptyset
		& \text{ if } \iota(d) \text{ is even}, \\
    \begin{array}[b]{l}
		\{j \in\dom(T) \colon A(d,T_j,h_d(j),x_{T_j-1})\neq 0 \text{ and } \\
		(\forall j'\in\dom(T)) \, (\text{if } j'<j \text{ and } A(d,T_{j'},h_d(j'),x_{T_{j'}-1}) \neq 0 \\
		\text{then } |A(d,T_{j'},h_d(j'),x_{T_{j'}-1})| \neq |A(d,T_j,h_d(j),x_{T_j-1})| ) \}
		\end{array}
		& \text{ if } \iota(d) \text{ is odd}.
	\end{cases}
\end{eqnarray*}
So, the sequence $(T_j)_j$ must be $(d,j,(x_t)_t)$-good for all those $j\in\dom(T)$
such that the value $A(d,T_j,h_d(j),x_{T_j-1})$ is nonzero and not the first with its absolute value.

In the case of odd $\iota(d)$ (that means, in the case when the requirement $\mathcal{R}_{2\iota(d)}$ needs to be taken care of infinitely often),
why do we not demand that it is $(d,j,(x_t)_t)$-good for all $j\in\dom(T)$ with $A(d,T_j,h_d(j),x_{T_j-1})\neq 0$?
This is because in this case we must make sure that the responsible strategy allows lower priority strategies to act
infinitely often, and this will be in stages $T_j$ such that the value $A(d,T_j,h_d(j),x_{T_j-1})$ is not zero 
and the first with its absolute value; compare step (II$_\mathrm{even}$) in the formal construction
in Section~\ref{subsection:3}. In these stages we will not even try to ensure that \eqref{equation:i-d-good} becomes true.
In the case of even $\iota(d)$ (that means, in the case when the requirement $\mathcal{R}_{2\iota(d)}$ needs to be taken care of only finitely often),
this problem does not appear as the responsible strategy will be active only finitely often.

Concerning the function $A$, we distinguish two cases
and we consider only functions $A$ for which one of these two cases applies:
on the one hand the \emph{nonnegative} case in which
\begin{equation}
\label{eq:nonnegative-case}
\mathrm{range}(A)\subseteq\{x\in\ID\colon x\geq 0\},
\end{equation}
on the other hand the \emph{sign switching} case in which 
\begin{equation}
\label{eq:signswitching-case}
\mathrm{range}(A)\cap\{x\in\ID\colon x>0\}\neq\emptyset \text{ and } 
\mathrm{range}(A)\cap\{x\in\ID\colon x<0\}\neq\emptyset .
\end{equation}

In the nonnegative case the sequences $(x_t)_t$ to be considered resp. that are obtained are nondecreasing.
In this case, sometimes it is useful to consider a weaker version of \eqref{equation:i-d-good}, namely the following condition:
\begin{equation}
\label{equation:i-d-leftgood}
   x_{T_{j+1}-1} - x_{T_j-1} \geq A(d,T_j,h_d(j),x_{T_j-1}) . 
\end{equation}
This condition means that the movement of $x_{T_j-1}$ by the amount of $A(d,T_j,h_d(j),x_{T_j-1})>0$, 
that has been demanded in stage $T_j$,
has been carried out until the next stage $T_{j+1}$ of this sequence,
and it is alright if the value $x_{T_{j+1}}$ has been moved even further to the right.
Given a nondecreasing sequence $(x_t)_t$ of dyadic numbers and a number $d\in\IN$,
a nonempty, finite or infinite, increasing sequence $(T_j)_j$ with $T_0>0$
is called \emph{$(d,(x_t)_t)$-left-good} if 
for all $j$ in the set
\[ \{j\in\dom(T) \colon A(d,T_j,h_d(j),x_{T_j-1})\neq 0 \text{ and } j+1\in\dom(T) \} \setminus \mathrm{First}( d,(T_k)_k ) \]
the condition \eqref{equation:i-d-leftgood} is satisfied.

We shall demand further properties
for a computable (in the sign switching case)
respectively for a computable and nondecreasing (in the nonnegative case)
sequence $(x_t)_t$ of dyadic numbers,
and for a computable, nonempty, finite or infinite sequence $(T_j)_j$, which may be $(d,(x_t)_t)$-good or $(d,(x_t)_t)$-left-good :
\begin{equation}
\label{equation:bound-on-sum}
\text{for any } J\in\dom(T), \ 
   \left| \sum_{j\in\{0,\ldots,J\} \setminus \mathrm{First}( d,(T_j)_j )} A(d,T_j,h_d(j),x_{T_j-1}) \right| \leq 2^{-d}.
\end{equation}
In the case when we wish to construct a weakly computable number, then we shall demand the following property,
which is stronger than \eqref{equation:bound-on-sum}:
\begin{equation}
\label{equation:bound-on-absolute-sum}
   \sum_{j\in\dom(T)} \left|A(d,T_j,h_d(j),x_{T_j-1}) \right| \leq 2^{-d}.
\end{equation}
And also:
\begin{itemize}
\item
in the case when $\iota(d)$ is even (that means that the requirement $\mathcal{R}_{2\iota(d)}$ needs to be taken care of only finitely often),
the following condition is important:
\begin{equation}
\label{equation:bounded}
\text{the function } h_d \text{ is bounded,}
\end{equation}
\item
in the case when $\iota(d)$ is odd (that means that the requirement $\mathcal{R}_{2\iota(d)}$ needs to be taken care of infinitely often)
and \eqref{eq:nonnegative-case} is true (that means we are in the nonnegative case;
actually this condition will be assumed to be satisfied only in the fourth case in Theorem~\ref{theorem:master}, the ``nonnegative and left-good case''),
the sequence $(x_t)_t$ is nondecreasing, and the sequence $(T_j)$ is infinite,
the following condition is important:
\begin{equation}
\label{equation:computable-convergence}
\text{the series } \sum_j A(d,T_j,h_d(j),x_{T_j-1}) \text{ is computably Cauchy}, 
\end{equation}
Of course, the statement that a series $\sum_j y_j$ is \emph{computably Cauchy}
simply means that the sequence $\left( \sum_{j=0}^{J-1} y_j\right)_J$ is computably Cauchy,
\item
in the case when $\iota(d)$ is odd (that means that the requirement $\mathcal{R}_{2\iota(d)}$ needs to be taken care of infinitely often)
and \eqref{eq:signswitching-case} is true (that means we are in the sign switching case)
the following condition is important: using the notation
\[ J_m := \{j\in\dom(T) \mid \left|A(d,T_j,h_d(j),x_{T_j-1})\right|=2^{-m}\}, \]
\begin{equation}
\label{equation:cancellation}
\begin{array}{l}
\bullet \text{ for any } m\in\IN, \text{ for any two consecutive numbers } j,j'\in J_m \text{ we have } \\
A(d,T_{j'},h_d(j'),x_{T_{j'}-1}) = - A(d,T_{j},h_d(j),x_{T_{j}-1}), \\
\bullet \text{ for any } m,m'\in\IN \text{ with } m<m'  \text{ if } J_{m'}\neq\emptyset \text{ then either } J_m \text{ is empty} \\
\text{or } J_m \text{ has an odd number of elements and } \max J_m < \min J_{m'}.
\end{array}
\end{equation}
\item
in the case when $\iota(d)$ is odd (that means that the requirement $\mathcal{R}_{2\iota(d)}$ needs to be taken care of infinitely often)
and \eqref{eq:signswitching-case} is true (that means we are in the sign switching case)
also the following condition is important: using the notation as in \eqref{equation:cancellation}
\begin{equation}
\label{equation:Jm-finite}
\text{ for every }  m \text{ the set } J_m \text{ is finite.}
\end{equation}
\end{itemize}

\begin{lemma}
\label{lemma:cancel-comp-conv}
Let $d\in\IN$ be such that $\iota(d)$ is odd.
Let $A:\IN^3\times\ID\to\ID$ be a computable function that satisfies \eqref{equation:range-powers-of-2} and \eqref{eq:signswitching-case}.
Let $(x_t)_t$ be a computable sequence of dyadic numbers,
and let $(T_j)_j$ be an infinite, increasing, and computable sequence of natural numbers with $T_0>0$
such that \eqref{equation:cancellation} and \eqref{equation:Jm-finite} are true.
Then \eqref{equation:computable-convergence} is true as well.
\end{lemma}

\begin{proof}
The assertion is clear if there are only finitely many $j$ with $A(d,T_j,h_d(j),x_{T_j-1})\neq 0$.
Let us assume that there are infinitely many $j$ with $A(d,T_j,h_d(j),x_{T_j-1})\neq 0$.
It follows from \eqref{equation:range-powers-of-2}, \eqref{equation:cancellation}, and \eqref{equation:Jm-finite}
that $\lim_{j\to\infty} A(d,T_j,h_d(j),x_{T_j-1}) = 0$. 
Given a number $n$, one can compute a number $J$ such that $|A(d,T_j,h_d(j),x_{T_j-1})| < 2^{-n}$ for all $j\geq J$.
It follows from \eqref{equation:cancellation} that for all for all $j_1,j_2$ with $J\leq j_1 \leq j_2$ and all $m$
\[ \left|\sum_{j\in \{j_1,\ldots,j_2-1\} \cap J_m} A(d,T_j,h_d(j),x_{T_j-1})\right| 
 \begin{cases} 
 =0 & \text{if } m\leq n, \\
 \leq 2^{-m} & \text{if } m > n .
\end{cases} \]
We obtain for all $j_1,j_2$ with $J\leq j_1 \leq j_2$
\[ \left| \sum_{j\in\{j_1,\ldots,j_2-1\}} A(d,T_j,h_d(j),x_{T_j-1}) \right| \leq \sum_{m>n} 2^{-m} \leq 2^{-n} . \]
Thus, the series $\sum_{j \in \IN\setminus\mathrm{First}(d, (T_j)_j)} A(d,T_j,h_d(j),x_{T_j-1})$ is computably Cauchy.
\end{proof}

As before, in the following theorem we treat the case where the requirements $\mathcal{R}_{2i}(\xi)$
with even $i$ need to be taken care of only finitely often
while the requirements $\mathcal{R}_{2i}(\xi)$
with odd $i$ need to be taken care of infinitely often.

\begin{theorem}
\label{theorem:master}
Let $(\mathcal{R}_{2i}(\xi))_i$ be an infinite sequence of requirements of a real number $\xi$.
\begin{enumerate}
\item
Sign switching case: 
Let us assume that there exists a computable function $A:\IN^3\times\ID\to\ID$ 
satisfying 
\eqref{equation:range-powers-of-2} and \eqref{eq:signswitching-case} and,
for any $d\in\IN$, for any computable sequence $(x_t)_t$ of dyadic numbers,
and for any computable, nonempty, finite or infinite, increasing sequence $(T_j)_j$ with $T_0>0$:
\begin{itemize}
\item
if $\iota(d)$ is even and $(T_j)_j$ is $(d,(x_t)_t)$-good, then also \eqref{equation:bound-on-sum} and \eqref{equation:bounded},
\item
if $\iota(d)$ is odd, then also \eqref{equation:cancellation}, and if $(T_j)_j$ is additionally $(d,(x_t)_t)$-good
then additionally \eqref{equation:bound-on-sum} and \eqref{equation:Jm-finite}.
\end{itemize}
Furthermore, let us assume that $A$ satisfies the following condition:
for every $i\in\IN$, if $(x_t)_t$ is a computable and converging sequence of dyadic numbers
such that there exists some $d\in\IN$ such that $\iota(d)=i$ and 
such that there exists a computable, infinite, increasing, $(d,(x_t)_t)$-good sequence $(T_j)_j$ with $T_0>0$,
then the requirement $\mathcal{R}_{2i}(\lim_{t\to\infty} x_t)$ is true.

Then there exists a nearly computable real number $\xi$ such that $\mathcal{R}_{2i}(\xi)$ is true for all $i$.
\item
Sign switching and weakly computable case:
Let us assume that there exists a computable function $A:\IN^3\times\ID\to\ID$ 
as in the sign switching case such that additionally for any $d\in\IN$,
for any computable sequence $(x_t)_t$ of dyadic numbers, 
and for any computable, nonempty, finite or infinite, increasing, $(d,(x_t)_t)$-good 
sequence $(T_j)_j$ with $T_0>0$
the condition \eqref{equation:bound-on-absolute-sum} is satisfied.

Then there exists a nearly computable and weakly computable real number $\xi$ such that $\mathcal{R}_{2i}(\xi)$ is true for all $i$.
\item
Nonnegative case:
Let us assume that all requirements $(\mathcal{R}_{2i}(\xi))_i$ for odd $i$ are true for all real numbers $\xi$
(that means that the only requirements that one has to take care of need to be taken care of only finitely often).
Let us assume that there exists a computable function $A:\IN^3\times\ID\to\ID$ 
satisfying 
\eqref{equation:range-powers-of-2} and \eqref{eq:nonnegative-case} and,
for any $d\in\IN$, for any computable, nondecreasing sequence $(x_t)_t$ of dyadic numbers,
and for any computable, nonempty, finite or infinite, increasing sequence $(T_j)_j$ with $T_0>0$:
\begin{itemize}
\item
if $\iota(d)$ is even and $(T_j)_j$ is $(d,(x_t)_t)$-good, then also \eqref{equation:bound-on-sum} and \eqref{equation:bounded},
\item
if $\iota(d)$ is odd, then $A(d,t,\ell,q)=0$ for all $t,\ell\in\IN$ and $q\in\ID$.
\end{itemize}
Furthermore, let us assume that $A$ satisfies the following condition:
for every even $i\in\IN$, if $(x_t)_t$ is a nondecreasing, computable, and converging sequence of dyadic numbers
such that there exists some $d\in\IN$ such that $\iota(d)=i$ and 
such that there exists a computable, infinite, increasing, $(d,(x_t)_t)$-good sequence $(T_j)_j$ with $T_0>0$,
then the requirement $\mathcal{R}_{2i}(\lim_{t\to\infty} x_t)$ is true.

Then there exists a nearly computable and left-computable real number $\xi$ such that $\mathcal{R}_{2i}(\xi)$ is true for all $i$.
\item
Nonnegative and left-good case:
Let us assume that there exists a computable function $A:\IN^3\times\ID\to\ID$ 
satisfying 
\eqref{equation:range-powers-of-2} and \eqref{eq:nonnegative-case} and,
for any $d\in\IN$, for any computable, nondecreasing sequence $(x_t)_t$ of dyadic numbers,
and for any computable, nonempty, finite or infinite, increasing sequence $(T_j)_j$ with $T_0>0$:
\begin{itemize}
\item
if $\iota(d)$ is even and $(T_j)_j$ is $(d,(x_t)_t)$-left-good, then also \eqref{equation:bound-on-sum} and \eqref{equation:bounded},
\item
if $\iota(d)$ is odd and $(T_j)_j$ is $(d,(x_t)_t)$-left-good then also \eqref{equation:bound-on-sum},
and if $(T_j)_j$ is additionally infinite then \eqref{equation:computable-convergence} is true as well.
\end{itemize}
Furthermore, let us assume that $A$ satisfies the following condition:
for every $i\in\IN$, if $(x_t)_t$ is a computable, nondecreasing and converging sequence of dyadic numbers
such that there exists some $d\in\IN$ such that $\iota(d)=i$ and 
such that there exists a computable, infinite, increasing, $(d,(x_t)_t)$-left-good sequence $(T_j)_j$ with $T_0>0$,
then the requirement $\mathcal{R}_{2i}(\lim_{t\to\infty} x_t)$ is true.

Then there exists a nearly computable and left-computable real number $\xi$ such that $\mathcal{R}_{2i}(\xi)$ is true for all $i$.
\end{enumerate}
\end{theorem}

\section{Proof of Theorem~\ref{theorem:master}}
\label{section:proof-master}

As the proof is fairly long, we present it in the following subsections:
\begin{enumerate}[itemsep=1pt,
labelindent=\parindent,
leftmargin=*,
align=left]
\item[\ref{subsection:1}]
The odd requirements and a rough sketch of the proof
\item[\ref{subsection:2}]
The infinite injury priority argument:
notation, first observations, and a sketch of the interplay of the various strategies
\item[\ref{subsection:3}]
The formal construction
\item[\ref{subsection:4}]
The correctness proof, starting with two general lemmas
\item[\ref{subsection:5}]
Proofs by induction concerning the triggered movements: either executed immediately or split up and stored for later execution
\item[\ref{subsection:6}]
About strategies of even length
\item[\ref{subsection:7}]
General lemmas concerning, e.g., only finitely many initializations of the strategies on the true path and convergence of the constructed sequence
\item[\ref{subsection:8}]
Satisfaction of the even requirements
\item[\ref{subsection:9}]
The sign switching and weakly computable case
\item[\ref{subsection:10}]
Satisfaction of the odd requirements
\end{enumerate}

\subsection{The odd requirements and a rough sketch of the proof}
\label{subsection:1}

We wish to construct a converging sequence $(x_t)_t$ of dyadic numbers 
such that its limit $\xi$ satisfies all of the requirements $(\mathcal{R}_{2i}(\xi))_i$, for $i\in\IN$.
In addition, we want the number $\xi$ to be nearly computable.
For the second point, we claim that it is sufficient to make sure that the following \emph{odd requirements}
\begin{equation}
\label{eq:left-req}
\tag{$\mathcal{R}_{2i+1}$}
\begin{array}{l}
\text{if $\varphi_i\colon\IN\to\IN$ is a total and increasing function then}\\
\text{the sequence }
 (x_{\varphi_i(n+1)} - x_{\varphi_i(n)})_n  \text{ converges computably to } 0,
\end{array}
\end{equation}
for $i\in\IN$, are satisfied. 
Indeed, let us consider an arbitrary total, computable, and increasing function $s:\IN\to\IN$. 
The sequence $(x_{s(n+1)} - x_{s(n)})_n$ converges computably to $0$ because there exists some
$i$ with $\varphi_i=s$ and because $(x_t)_t$ satisfies requirement $\mathcal{R}_{2i+1}$.

The computable function $A$, that is supposed to satisfy the assumptions
of Theorem~\ref{theorem:master}, 
has the duty to make sure that the even requirements $(\mathcal{R}_{2i})_i$ are satisfied.
The idea is always the following.
When at the beginning of a stage $t>0$ the constructed number $x_{t-1}$ is close to a certain
computed dyadic number (let us call is $q_{t-1}$) and some additional conditions are satisfied
then the value $A(d,t,h_{d}(\ldots),x_{t-1})$ that is computed in stage $t$,
will be different from $0$, and the algorithm will try to move $x_{t-1}$ by this value.
When this succeeds, that is, when there is a $t'>t$ such that $x_{t'}$ is close to
$x_{t-1}+A(d,t,h_{d}(\ldots),x_{t-1})$, then there are two possibilities: either the other computed sequence
of dyadic numbers $(q_t)_t$ does not follow this move and never comes close to $x_{t''}$ for any $t''\geq t'$ again.
Then the even requirement $\mathcal{R}_{2i}$ should be satisfied trivially.
Or it comes close to $x_{t''}$ for some $t''\geq t'$ again. Then the same procedure is followed again.
In this way the other computed sequence $(q_t)_t$ is forced to make certain jumps to the left or right.
Now the even requirement $\mathcal{R}_{2i}$ should be of such a kind that,
if $(q_t)_t$ makes a certain number (this may be finitely many or infinitely many)
of jumps of certain sizes, then $\mathcal{R}_{2i}$ should be satisfied.

How is this procedure combined with a strategy for satisfying the odd requirements $(\mathcal{R}_{2i+1})_i$?
The strategy for satisfying an odd requirement $\mathcal{R}_{2i+1}$ is as follows.
We try to make sure that from a certain stage $t_0$ on (that means: for all $t\geq t_0$)
the value $x_{t-1}$ is always quite close to $x_{\varphi_i(\ell_i(t)-1)}$.
In this way we make sure that the sequence $(x_{\varphi_i(n)})_n$
can make only small steps to the right or to the left.
By decreasing the possible difference from time to time we will ensure that the
sequence $(x_{\varphi_i(n+1)} - x_{\varphi_i(n)})_n$ converges computably.
But there seems to be a conflict.
In order to satisfy an even requirement, $x$ should make a possibly big jump to the right or to the left,
but an odd requirement $\mathcal{R}_{2i+1}$ of higher priority may demand that $x$
can move to the right or the left only a bit
so that it still is only slightly different from $x_{\varphi_i(\ell_i(t)-1)}$.
The idea for resolving this conflict is to split the big move that $x$ is supposed
to make into many small jumps. Each of these small jumps is done only when 
$x_{\varphi_i(\ell_i(t)-1)}$ is close enough to $x_{t-1}$
(remember that always $\varphi_i(\ell_i(t)-1)\leq t-1$; when we do not change $x_t$ much from a certain stage $t_0$
on then, once $\varphi_i(\ell_i(t)-1)\geq t_0$, the value $x_{\varphi_i(\ell_i(t)-1)}$ will be close enough to $x_{t-1}$;
so, we just have to wait for $\varphi_i(\ell_i(t)-1)$ to be large enough).
The strategies are combined by an infinite injury priority argument.
As usual, whenever a higher priority strategy has achieved a certain partial goal,
then it leaves open some window for lower priority strategies to act.

\subsection{The infinite injury priority argument:
notation, first observations, and a sketch of the interplay of the various strategies}
\label{subsection:2}

In order to define a priority between the strategies, we introduce some more notation.
Let $\Sigma := \{0, 1\}$ be the binary alphabet.
For any $n\in\IN$, any element of $\Sigma^n$ is called a {\em binary string of length $n$}.
There is exactly one binary string of length $0$; we write is as $\lambda$.
Remember that $\Sigma^*=\bigcup_{n\in\IN}\Sigma^n$ is the set of all finite binary strings.
Let $\Sigma^\mathrm{odd}$ be the set of all binary strings of odd length, 
let $\Sigma^\mathrm{even}$ be the set of all binary strings of even length,
let $\Sigma^{\equiv 0 \mathrm{mod} 4}$ be the set of all binary strings whose length is divisible by $4$,
and let $\Sigma^{\equiv 2 \mathrm{mod} 4}$ be the set of all binary strings  whose length is divisible by $2$ but not by $4$.
The elements of $\Sigma^\IN$ will be called {\em infinity binary strings}.
We say that a finite or infinite binary string $\sigma$ is a {\em prefix} of a finite or infinite binary string $\tau$ if $\sigma\subseteq\tau$.
It is a {\em proper prefix} if $\sigma\subsetneq\tau$.
If $\sigma$ is a finite binary string and $a\in\Sigma$ then $\sigma a$ is the finite binary string
of length $|\sigma|+1$ with $\sigma \subseteq\sigma a$ and $(\sigma a)(|\sigma|)=a$.
For arbitrary $\sigma, \tau \in \Sigma^*$ we say that $\sigma$ is \emph{lexicographically smaller} than $\tau$ 
($\sigma <_L \tau$ for short) if there exists a $\rho \in \Sigma^*$ with $\rho 0 \subseteq \sigma$ and $\rho 1 \subseteq \tau$.

Any two finite binary strings $\sigma,\tau$ are either identical or one of them is a proper prefix of the other one or one of them
is lexicographically smaller than the other one.
For two arbitrary finite binary strings $\sigma,\tau$ we can apply this case distinction to the two strings $\sigma0$ and $\tau0$
and conclude that exactly one of the following five conditions is satisfied:
\begin{enumerate}
\item
$\sigma0 <_L \tau$,
\item
$\sigma0\subseteq\tau$,
\item
$\sigma=\tau$,
\item
$\tau0\subseteq\sigma$,
\item
$\tau0 <_L \sigma$.
\end{enumerate}

The construction of the sequence $(x_t)_t$ is done in infinitely many stages $0,1,2,\ldots$.
Every stage $t>0$ consists of at most $t$ substages.
In each substage a certain strategy determined by a binary string $\sigma$ of length smaller than $t$ is applied.
In each stage we start with the empty string $\lambda$. 
The strategy corresponding to some binary string $\sigma$ (short: the strategy $\sigma$)
is supposed to check whether something needs to be done with respect to the requirement
$\mathcal{R}_{|\sigma|}$ and to take some action if this is necessary.
Depending on a certain result of the strategy $\sigma$ we either set $\truepath[t]:=\sigma$,
end stage $t$ and continue with stage $t+1$, or we extend $\sigma$ by either $0$ or $1$ to a string $\tau=\sigma0$ or $\tau=\sigma1$.
In the second case, if the length of $\tau$ is smaller than $t$ then the next strategy after $\sigma$ applied in stage $t$
is the one corresponding to $\tau$.
If the length of $\tau$ is equal to $t$ then we set $\truepath[t]:=\tau$, end stage $t$ and continue with stage $t+1$.
We say that a strategy $\sigma$ \emph{is applied} in stage $t$ if $\sigma \prefixeq \truepath[t]$.
Based on the sequence of binary strings $(\truepath[t])_t$ we define the \emph{true path} $\truepath \in \Sigma^{\IN}$
recursively as follows:
\begin{align*}
\truepath(e) &:= \begin{cases}
0 &\text{if $(\truepath \restrictedto e) 0 \subseteq \truepath[t]$ for infinitely many $t \in \IN$}, \\
1 &\text{otherwise}.
\end{cases}
\end{align*}

\begin{lemma}
\label{lemma:infinite-injury-applied}
For every binary string $\tau\subseteq\truepath$ there exists a number $t_\tau$ such that,
for all $t\geq t_\tau$ and all $\sigma \in \Sigma^*$, if $\sigma <_L \tau$ then $\sigma \not\subseteq \truepath[t]$. 
\end{lemma}

\begin{proof}
Let $\tau$ be a prefix of $\truepath$. If $\tau$ does not contain a $1$ then there is no $\sigma$ with 
$\sigma <_L \tau$. In that case the number $t_\tau:=0$ does the job.
Otherwise, for every string $\rho$ such that $\rho 1 \subseteq \tau$ let
$s_\rho:= \min\{s\in\IN ~:~ (\forall t\geq s) \ \rho0\not\subseteq \truepath[t]\}$.
Then the number $t_\tau:= \max\{s_\rho ~:~ \rho \in\Sigma^*,\ \rho 1\subseteq \tau\}$
has the desired property.
\end{proof}

Before we come to the details of the construction, we try to give an intuition concerning the infinite injury priority argument
and the interplay between strategies of even or odd length and of shorter or larger length.
Every strategy $\sigma$ takes care of the requirement $\mathcal{R}_{|\sigma|}$.
At the end of any stage $t$, all strategies $\tau$ with $\mathcal{T}[t] <_L \tau$ will be initialized,
and in some cases, even all strategies $\tau$ with $\mathcal{T}[t]0 <_L \tau$ will be initialized.
Let us look at the strategies $\sigma$ on the true path, that is, with $\sigma \prefixeq \mathcal{T}$.
In the end, these strategies will be the important ones.
We will see that they will be applied infinitely often. Let us look at the interplay between these strategies.
\begin{itemize}
\item
If $\sigma$ has even length, is on the true path, and is active only finitely often (that will be the case
when the even requirement $\mathcal{R}_{|\sigma|}$ needs to be taken care of only finitely often until it is satisfied forever)
then the string $\sigma1$ will be a prefix of the true path, and at the end of a stage $t$
in which $\sigma$ is applied, that is, with $\sigma\prefixeq\mathcal{T}[t]$,
all longer strategies on the true path will be initialized, the even ones and the odd ones.
But that is all right as $\sigma$ will be active only finitely often.
\item
If $\sigma$ has even length, is on the true path, and is active infinitely often (that will be the case
when the even requirement $\mathcal{R}_{|\sigma|}$ needs to be taken care of infinitely often)
then the string $\sigma0$ will be a prefix of the true path, and at the end of the corresponding stage
the longer strategies on the true path will not be initialized.
In this case $\sigma$ will from time to time allow the longer strategies on the true path to act.
But, most of its own actions will not be disturbed by the actions of longer strategies on the true path, and this is all right, as far
as $\sigma$ is concerned.
The question arises whether the actions of the longer strategies on the true path will be disturbed too much by the actions of $\sigma$.
We claim that this is not the case.

Let us first look at longer prefixes of even length of the true path.
First, in the sign switching case and the sign switching and weakly computable case 
(the first and the second case in Theorem~\ref{theorem:master}) the jumps of the constructed sequence $(x_t)_t$ 
triggered by $\sigma$ go to the right direction and to the left direction in turns. Therefore, as long as they have the same absolute value
they cancel each other out. This has the effect that the total amount added to the sequence $(x_t)_t$ by $\sigma$ between two stages in which $\sigma$
allows lower priority strategies to act, is equal to zero. 
This is sufficient for longer prefixes of even length of the true path.
Secondly, in the third case of Theorem~\ref{theorem:master} we excluded the possibility of strategies
$\sigma$ of even length on the true path that are active infinitely often.
Thirdly, in the nonnegative and left-good case, the fourth case of Theorem~\ref{theorem:master},
the strategies of even length do not disturb each other at all, because all of them just want to move $x_t$ to the right,
and if $x_t$ is moved to the right even more, then that is all right.
So, in this case there is no problem for longer prefixes of even length of the true path.

For longer prefixes of odd length of the true path (such strategies $\tau$ have to make sure that the odd requirement $\mathcal{R}_{|\tau|}$ is satisfied,
that is that for $|\tau|=2i+1$ the sequence $(x_{\varphi_i(n+1)} - x_{\varphi_i(n)})_n$ converges computably to $0$)
we observe additionally that the jumps triggered by $\sigma$ will become smaller and smaller and tend to $0$.
This has the effect that for larger stages the jumps triggered by $\sigma$ will only lead to small oscillations of $x_t$ around some point,
which is enough for the computable convergence.
\item
If $\sigma$ has odd length, is on the true path, and is active infinitely often (that will be the case
when the odd requirement $\mathcal{R}_{|\sigma|}$ needs to be taken care of infinitely often,
that is, when that for $|\sigma|=2i+1$ the function $\varphi_i$ is total and increasing and $\sigma$ has to ensure that
the sequence $(x_{\varphi_i(n+1)} - x_{\varphi_i(n)})_n$ converges computably to $0$)
then the longer strategies on the true path will not be initialized.
In this case $\sigma$ will from time to time allow the longer strategies on the true path to act.
Its own actions will not be disturbed by actions of longer strategies on the true path because
the windows during which $\sigma$ will allow longer strategies on the true path to act
will be open only when $\sigma$ has finished the tasks it had been given and has not yet received new tasks.
We already explained how such a $\sigma$ interacts with a shorter prefix of even length of the true path.
How does it interact with a longer prefix of even length of the true path?
Well, when a longer prefix $\tau$ of even length of the true path tries to trigger a jump of the constructed sequence $(x_t)_t$
then $\sigma$ will be asked to cut this jump into many small hops and to take care of these itself or to delegate them further to 
even shorted prefixes of odd length of the true path. 
Finally, the same happens when a longer prefix $\tau$ of odd length of the true path tries to trigger a jump of the constructed sequence $(x_t)_t$.
\end{itemize}

\subsection{The formal construction}
\label{subsection:3}

Each strategy $\sigma\in\Sigma^*$ is equipped with certain \emph{$\sigma$-parameters}.
We have for every $\sigma \in \Sigma^\mathrm{even}$ and for every $t \in \IN$:
\begin{itemize}
\item
a number $d(\sigma)[t]\in\IN$ (this is a precision parameter),
\end{itemize}
and for every $\sigma \in \Sigma^\mathrm{odd}$ and for every $t \in \IN$:
\begin{itemize}
\item
a dyadic number $q(\sigma)[t] \in\ID$ ($\sigma$ wants to add $q(\sigma)[t-1]$ to the current value $x_{t-1}$),
\item
a \emph{restraint} $r(\sigma)[t] \in\IN$ ($\sigma$ forbids jumps of absolute size bigger then $2^{-r(\sigma)[t]}$).
\end{itemize}
The algorithm computes a sequence of dyadic numbers $(x_t)_t$.
We will show that it converges and that its limit $\xi:=\lim_{t\to \infty} x_t$ has the desired properties.
For each $t\in\IN$, the number $x_t$ is defined in stage $t$.
Additionally, in each stage $t$ a dyadic number $v_t\neq 0$ is defined.
If some action is required in stage $t$, then this number $v_t$ is just the new number that the algorithm would like to add to $x_{t-1}$, 
and either this number is directly added to $x_{t-1}$ or the task to add this number (at some later stage) is delegated
to another strategy.
Otherwise, that is, when no action is required in stage $t$,
then $v_t$ is simply set to $1$ (so that in the formula below for the initialization of $d(\sigma)$ 
it does not disturb).

We say that a strategy $\sigma \in \Sigma^*$ \emph{is initialized} in a stage $t > 0$, if the parameters of the strategy $\sigma$ (with suffix $[t]$) are set to the following values:
for $\sigma$ of even length:
\[
d(\sigma)[t]	:= \min\{b(\sigma,m) \colon \begin{array}[t]{l}
    m\in\IN, b(\sigma,m)>\max (\{d(\sigma)[t-1],2-\log_2(|v_t|)\}\cup \\
    \{1+r(\gamma)[t] \mid |\gamma| \text{ is odd and either } \gamma0\prefixeq\mathcal{T}[t] \text{ or } 
		 \gamma=\mathcal{T}[t]\})\} ,
		\end{array} 
\]
and for $\sigma$ of odd length:
\[ q(\sigma)[t] 	:= 0. \]
Remember that we fixed the function $b$ in the beginning of Section~\ref{section:master} when we explained the purpose of the input parameters of the function $A$. The function $b$ is injective. So, every strategy $\sigma$ of even length has its own sequence of possible values for the precision parameter $d(\sigma)[t]$.
Of course, the new value $d(\sigma)[t]$ must be chosen large enough.
The condition that the new value should be larger than $2-\log_2(|v_t|)$ 
makes sure that the contributions of initialized strategies are so small that they do not disturb the movement triggered by the current strategy (which has higher priority) too much;
see quite at the end of the first part of the proof of Lemma~\ref{lemma:even-movement}.
The condition that the new value should be larger than 
$\max \{1+r(\gamma)[t] \mid |\gamma| \text{ is odd and either } \gamma0\prefixeq\mathcal{T}[t] \text{ or } 
		 \gamma=\mathcal{T}[t]\}$
makes sure that the contributions of initialized strategies are so small that they do not get into conflict with strategies of higher priority and of odd length;
see the end of the proof of Lemma~\ref{lemma:main-odd}.
Note that the restraint parameter $r(\sigma)$ for $\sigma$ of odd length is never initialized.
Any parameter not explicitly set during some stage gets the same value assigned to it
as at the end of the previous stage.

Every substage consists of several steps.
It depends on a strategy $\sigma$ which of these steps are carried out for $\sigma$.
In the following we say that some strategy $\sigma$ \emph{reaches a step (X) in stage $t$},
for $X\in\{\mathrm{I},\mathrm{II}_\mathrm{even},\mathrm{III}_\mathrm{even},\mathrm{II}_\mathrm{odd},\mathrm{III}_\mathrm{odd},\mathrm{IV}\}$ and $t>0$,
if $\sigma \subseteq \mathcal{T}[t]$
and during the substage of stage $t$ corresponding to $\sigma$ step (X) is reached.

\subsection*{Formal construction of the sequence $(x_t)_t$}

{\bf Stage $0$:}
We set $x_0:=0$ and $v_t:=1$ and $\mathcal{T}[0]:=\lambda$.
We set $d(\sigma)[0]:=b(\sigma,0)$, 
for all $\sigma\in\Sigma^\mathrm{even}$,
and $q(\sigma)[0]:=0$ and $r(\sigma)[0]:=0$, for all $\sigma\in\Sigma^\mathrm{odd}$.
Then we continue with stage $1$.

\bigskip

\noindent
{\bf Stage $t$ for $t>0$:}

We start with $\sigma:=\lambda$ and then proceed as follows.
\begin{enumerate}[labelindent=\parindent,
align=left]
\item[(I)]
If $|\sigma|=t$ 
then we set $v_t:=1$,
we initialize all $\tau$ with $\sigma<_L \tau$,
we set $\truepath[t]:=\sigma$, and we continue with stage $t+1$.
Otherwise, we continue with (II$_\mathrm{even}$) if $\sigma$ is of even length (then we set $i:=|\sigma|/2$)
and with (II$_\mathrm{odd}$) if $\sigma$ is of odd length (then we set $i:=(|\sigma|-1)/2$).
\end{enumerate}
The following steps (II$_\mathrm{even}$) and (III$_\mathrm{even}$) are for strategies $\sigma$ of even length $2i$ and, in combination with step (IV), 
take care of the even requirement $\mathcal{R}_{2i}$.
\begin{enumerate}[labelindent=\parindent,
align=left]
\item[(II$_\mathrm{even}$)]
Let $d_0:=d(\sigma)[t-1]$. 
Let $(T_j)_j$ be the nonempty, finite, increasing sequence of stages $s$ with $0< s \leq t$ and $\sigma\subseteq \mathcal{T}[s]$
and $d(\sigma)[s-1]=d_0$, 
and in which $\sigma$ reaches step (II$_\mathrm{even}$).
This sequence is nonempty because $t$ is such a stage,
so for $j_t:=|\dom(T)|-1$ we have $\dom(T)=\{0,\ldots,j_t\}$ and $T_{j_t}=t$.
We compute
\[ v:=A(d_0,t,h_{d_0}(j_t),x_{t-1}) . \]
If $j_t \in \mathrm{First}( d_0,(T_j)_j )$
(note that this implies that $i=\iota(d_0)$ is odd (in other words: $\sigma\in\Sigma^{\equiv 2 \mathrm{mod } 4}$) and $v\neq 0$)
then we set $\sigma:=\sigma0$ and continue with (I). 
Otherwise we continue with (III$_\mathrm{even}$).

{\em Explanation: 
This step is important if $\sigma$ reaches step (IV) in infinitely many stages $s$ (this can be the case only if $i$ is odd).
In this case we wish to make sure that there are infinitely many stages in which extensions $\tau\in\Sigma^\mathrm{even}$ of $\sigma0$ are applied.
The current stage $t$ is such a stage
if $i$ is odd, if the current value $A(d_0,\ldots)$ is nonzero,
and if the current value $A(d_0,\ldots)$ is the first one with its absolute value.
}
\item[(III$_\mathrm{even}$)]
If $v\neq 0$ then we set $v_t:=v$ and continue with (IV).
Otherwise we set $\sigma:=\sigma1$ and continue with (I).

{\em Explanation: 
If $v\neq 0$ then we would like to move $x$ by adding $v$ to $x_{t-1}$. 
This will be taken care of in step (IV).
If $v=0$ then at present we do not need to take care of the requirement $\mathcal{R}_{2i}$.
}
\end{enumerate}
The following steps (II$_\mathrm{odd}$) and (III$_\mathrm{odd}$) are for strategies $\sigma$ of odd length $2i+1$ and,
in combination with step (IV), take care of the odd requirement $\mathcal{R}_{2i+1}$.
\begin{enumerate}[labelindent=\parindent,
align=left]
\item[(II$_\mathrm{odd}$)]
If $\ell_i(t)>0$ and
\begin{equation}
\label{equation:III-1}
   \varphi_i(\ell_i(t)-1) \geq \max \left( \{0\} \cup \{s<t \mid \sigma \text{ reaches step (III$_\mathrm{odd}$) in stage } s \} \right)
\end{equation}
then we continue with (III$_\mathrm{odd}$),
otherwise we set $\sigma:=\sigma1$ and continue with (I).

{\em Explanation: In the ``otherwise'' case either $\ell_i(t)=0$ or $\ell_i(t)>0$ but $\varphi_i(\ell_i(t)-1)$ is still too small;
it has not (yet) reached the previous stage in which $\sigma$ reached step (III$_\mathrm{odd}$).}
\item[(III$_\mathrm{odd}$)]
If $q(\sigma)[t-1]=0$ then we set 
\[ r(\sigma)[t]:= r(\sigma)[t-1]+1,
\]
we set $\sigma:=\sigma0$, and we continue with (I).
Otherwise, we set 
\[ v_t := \sign(q(\sigma)[t-1]) \cdot \min\{ |q(\sigma)[t-1])|, 2^{-r(\sigma)[t-1]} \} \]
and
\begin{eqnarray*}
 q(\sigma)[t] &:=& q(\sigma)[t-1]-v_t \\ 
 &=& \begin{cases} 
    q(\sigma)[t-1] - 2^{-r(\sigma)[t-1]} & \text{if } q(\sigma)[t-1] > 2^{-r(\sigma)[t-1]}, \\
    q(\sigma)[t-1] + 2^{-r(\sigma)[t-1]} & \text{if } q(\sigma)[t-1] < - 2^{-r(\sigma)[t-1]}, \\
    0 & \text{if } |q(\sigma)[t-1]| \leq 2^{-r(\sigma)[t-1]}
    \end{cases}
\end{eqnarray*}
and continue with (IV).

{\em Explanation: If the number $q(\sigma)[t-1]$ is equal to zero, then one may say
that the strategy $\sigma$ has reached its current goal;
now it allows lower priority strategies to act; we just increase its restraint $r(\sigma)$.
If the number $q(\sigma)[t-1]$ is not equal to zero, then
$\sigma$ has not yet reached its current goal.
On the one hand, we will see that the condition checked in step (II$_\mathrm{odd}$) implies that
$x_{t-1}$ is so close to $x_{\varphi_i(\ell_i(t)-1)}$ that $\sigma$ allows us to move $x$ somewhat,
namely by $\pm 2^{-r(\sigma)[t-1]}$.
On the other hand, we would in fact like to add the amount $q(\sigma)[t-1]$ to $x_{t-1}$.
But we do not dare to add more than $\sign(q(\sigma)[t-1]) \cdot 2^{-r(\sigma)[t-1]}$.
Therefore, in step (IV) we will try to add the one of these two numbers 
that has the smallest absolute value.
}
\end{enumerate}
The following step (IV) is for strategies $\sigma$ both of even or of odd length.
\begin{enumerate}[labelindent=\parindent,
align=left]
\item[(IV)]
{\em Explanation:
We would like to add $v_t$ to $x_{t-1}$.
If there is some higher priority strategy that might restrain us from doing so,
then we formulate it as a ``task'' for the lowest higher priority strategy of this kind
that $x$ is moved by this amount.}

We check whether there exists a string $\gamma$ of odd length
such that $\gamma0 \subseteq \sigma$ 
(note that this implies that $q(\gamma)[t-1]=0$).
If there does not exist such a string $\gamma$ then we set
\[ x_t := x_{t-1} + v_t . \]
If there exists such a string $\gamma$
then we fix the longest string $\gamma$ with this property
and set $q(\gamma)[t]:=v_t$.
Finally, we initialize all $\tau$ with $\sigma0 <_L \tau$,
we set $\truepath[t]:=\sigma$ and continue with stage $t+1$.
\end{enumerate}
This completes the formal construction of the sequence $(x_t)_t$.

\subsection{The correctness proof, starting with two general lemmas}
\label{subsection:4}

We come to the correctness proof.
We wish to show that the limit $\xi:=\lim_{t\to\infty} x_t$ exists and has all the desired properties.
It is clear that the sequence $(x_t)_t$ is a computable sequence of dyadic numbers.
Furthermore, it is clear that in the nonnegative case the sequence $(x_t)_t$ is nondecreasing.
Thus, once we have shown that the sequence converges, in the nonnegative case its limit will be a left-computable number.
We assume that the function $A$ satisfies \eqref{equation:range-powers-of-2},
and in the respective cases also \eqref{equation:bound-on-sum}, \eqref{equation:bounded}, \eqref{equation:computable-convergence}, \eqref{equation:cancellation}, \eqref{equation:Jm-finite}.
When it becomes important whether \eqref{eq:signswitching-case} or \eqref{eq:nonnegative-case} is assumed, we shall point out which one is assumed.
We shall also demonstrate that, if in the sign switching state aditionally the condition \eqref{equation:bound-on-absolute-sum} is satisfied
then the sequence $(x_t)_t$ converges weakly computably.

We observe that in any stage at most one strategy can reach step (IV)
and if in a stage $t$ a strategy $\sigma$ reaches step (IV) then $\sigma=\mathcal{T}[t]$.
And a strategy can only be initialized in a stage in which it is not applied.

The following two lemmas contain assertions on strategies that can never be applied or must be applied
infinitely often, under certain circumstances.

\begin{lemma}
\label{lemma:div4-no-0}
Let $\sigma$ have length divisible by $4$.
For any $t\in\IN$ we have $\sigma0\not\prefixeq\mathcal{T}[t]$.
\end{lemma}

\begin{proof}
For the sake of a contradiction, let us assume that there are a $t\in\IN$ and a strategy $\sigma$ of length divisible by $4$ with $\sigma0\prefixeq\mathcal{T}[t]$.
Then we have $\sigma\prefixeq\mathcal{T}[t]$.
In stage $t$ the strategy $\sigma$ must first reach (I).
Then it must reach (II$_\mathrm{even}$) because otherwise we would obtain $\sigma=\mathcal{T}[t]$ and $\sigma0\not\prefixeq\mathcal{T}[t]$.
From step (II$_\mathrm{even}$) it must reach step (III$_\mathrm{even}$) because the length of $\sigma$ is divisible by $4$.
If then the algorithm would continue with (I) then we would get $\sigma1\prefixeq\mathcal{T}[t]$ and $\sigma0\not\prefixeq\mathcal{T}[t]$.
Thus, the algorithm must continue with (IV). But this leads to $\sigma=\mathcal{T}[t]$ and $\sigma0\not\prefixeq\mathcal{T}[t]$, contradiction.
\end{proof}

\begin{lemma}
\label{lemma:odd-infinitely-oftenIV}
Let $\sigma$ have odd length.
If there are infinitely many $t$ such that 
$\sigma$ reaches step (III$_\mathrm{odd}$) in stage $t$ then
there are infinitely many $t$ with $\sigma 0 \prefixeq \mathcal{T}[t]$.
\end{lemma}

\begin{proof}
Let us assume that there are infinitely many $t$ such that 
$\sigma$ reaches step (III$_\mathrm{odd}$) in stage $t$.
We wish to show that then there are infinitely many $t$ with $\sigma 0 \prefixeq \mathcal{T}[t]$.
For the sake of a contradiction,
let us assume that there are only finitely many $t$ with $\sigma 0 \prefixeq \mathcal{T}[t]$.
Let us consider a number $s>0$ so large that there are no $t\geq s$ with $\sigma 0 \prefixeq \mathcal{T}[t]$.
This implies that for all $t\geq s$ such that $\sigma$ reaches step (III$_\mathrm{odd}$) in stage $t$ 
we have $q(\sigma)[t-1]\neq 0$. In particular, $r(\sigma)[t-1]=r(\sigma)[s-1]$ for all $t\geq s$
(remember that initializations do not change $r(\sigma)$).
On the one hand, after stage $s-1$ the absolute value of the number $q(\sigma)$
can never increase because this could happen only if at some stage
$t'\geq s$ a strategy $\tau$ with $\sigma 0\prefixeq\tau$ would reach step (IV).
But this would imply $\sigma 0\prefixeq\tau\prefixeq\mathcal{T}[t']$.
On the other hand, in each stage $t\geq s$ in which $\sigma$ reaches step (III$_\mathrm{odd}$),
the absolute value of the number $q(\sigma)$ will be decreased by $2^{-r(\sigma)[t-1]}$ or will be set to $0$.
But then, after finitely many steps, it must have reached the value zero. 
This is a contradiction.
\end{proof}

Later we shall formulate a similar lemma (Lemma~\ref{lemma:4mod2-0}) for $\sigma$ whose length is divisible by $2$ but not by $4$.

\subsection{Proofs by induction concerning the triggered movements: either executed immediately or split up and stored for later execution}
\label{subsection:5}

If a strategy $\sigma$ reaches step (IV) in a stage $t$ with a value $v_t$ then
$v_t\neq 0$ and it is the goal to add $v_t$ to $x$.
If $\sigma$ has even length then the value $v_t$ has been defined in step (III$_\mathrm{even}$) of stage $t$.
If $\sigma$ has odd length then the value $v_t$ has been defined in step (III$_\mathrm{odd}$) of stage $t$
and we must have $q(\sigma)[t-1]\neq 0$. But then the nonzero value $q(\sigma)[t-1]$ must have been triggered
in a stage earlier than stage $t$. 
In the following we analyze when and under which circumstances this happens.

We define a partial function $\mathrm{trigd}:\subseteq \Sigma^\mathrm{odd}\times\IN\to\IN$ as follows:
\begin{align*}
    \dom(\mathrm{trigd}) &:= \{(\gamma,t) \in \Sigma^\mathrm{odd}\times\IN \colon 
		q(\gamma)[t]\neq 0\}, \\
		\mathrm{trigd}(\gamma,t) &:= 1+ \max \{s < t \colon q(\gamma)[s] = 0 )\}, 
    \text{ for } (\gamma,t)\in \dom(\mathrm{trigd}) .
\end{align*}
Note that $q(\gamma)[0]=0$ for all $\gamma\in\Sigma^\mathrm{odd}$.
Hence, the function $\mathrm{trigd}$ is well defined.
The assertions in the following two lemmas are immediate from the definitions and the construction.

\begin{lemma}
\label{lemma:trig-1}
Let us consider a pair $(\sigma,t_3)\in\dom(\mathrm{trigd})$, and let us set
$t_2:=\mathrm{trigd}(\sigma,t_3)$.
\begin{enumerate}
\item
$1 \leq t_2 \leq t_3$.
\item
If $t_1\leq t_3$ is a number such that $q(\sigma)[t_1]=0$ then $t_1 < t_2$.
\item
For all $s$ with $t_2 \leq s \leq t_3$
we have $q(\sigma)[s]\neq 0$ and $\sigma$ is not initialized in stage $s$.
\item
The stage $t_2$ is the largest stage $t\leq t_3$ such that the strategy $\mathcal{T}[t]$ has
the following two properties: it reaches step (IV) in stage $t$
and $\sigma$ is the longest strategy of odd length with $\sigma0\prefixeq\mathcal{T}[t]$.
\end{enumerate}
\end{lemma}

\begin{lemma}
\label{lemma:trig-2}
Let us consider a pair $(\sigma,t_2)\in\Sigma^\mathrm{odd} \times \IN$.
If $\sigma$ reaches step (IV) in stage $t_2$ then $(\sigma,t_2-1)\in\dom(\mathrm{trigd})$.
\end{lemma}

Let 
\[ \mathrm{StepIV} := \{t\in\IN\colon \mathcal{T}[t] \text{ reaches step (IV) in stage } t\} . \]
We define a function $\mathrm{trigIV}:\subseteq\mathrm{StepIV}\to\mathrm{StepIV}$ by
\begin{align*}
   \dom(\mathrm{trigIV}) &:= \{t \in \mathrm{StepIV} \colon |\mathcal{T}[t]| \text{ is odd}\}, \\
	 \mathrm{trigIV}(t) &:= \mathrm{trigd}(\mathcal{T}[t],t-1), \text{ for } t\in  \dom(\mathrm{trigIV}).
\end{align*}
By Lemma~\ref{lemma:trig-2}, for any $t\in\mathrm{StepIV}$ such that $\mathcal{T}[t]$ has odd length, the number
$ \mathrm{trigIV}(t)$ is defined.
It is clear that $1 \leq \mathrm{trigIV}(t)<t$ for every $t\in\dom(\mathrm{trigIV})$.
By Lemma~\ref{lemma:trig-1} for every $t\in\dom(\mathrm{trigIV})$, $\mathrm{trigIV}(t)\in\mathrm{StepIV}$ 
and $\mathcal{T}[t]$ is the longest string $\gamma$ of odd length with $\gamma0\prefixeq\mathcal{T}[\mathrm{trigIV}(t)]$.
The $j$-th iteration $\mathrm{trigIV}^{\circ j}\subseteq\mathrm{StepIV}\to\mathrm{StepIV}$ of $\mathrm{trigIV}$, for $j\in\IN$, is defined recursively:
$\mathrm{trigIV}^{\circ 0} := \mathrm{id}_\mathrm{StepIV}$ or, in more detail:
\begin{align*}
   \dom(\mathrm{trigIV}^{\circ 0}) &:= \mathrm{StepIV}, \\
	 \mathrm{trigIV}^{\circ 0}(t) &:=  t, \text{ for } t\in \mathrm{StepIV},
\end{align*}
and, for $j\in\IN$,
\begin{align*}
   \dom(\mathrm{trigIV}^{\circ j+1}) &:= \{t \in  \dom(\mathrm{trigIV}^{\circ j}) \colon |\mathcal{T}[ \mathrm{trigIV}^{\circ j}(t)]| \text{ is odd}\}, \\
	 \mathrm{trigIV}^{\circ j+1}(t) &:=  \mathrm{trigIV}( \mathrm{trigIV}^{\circ j}(t)) \text{ for } t\in \dom(\mathrm{trigIV}^{\circ j+1}).
\end{align*}

For any strategy $\sigma\in\Sigma^*$ we define
\[ \ell(\sigma):= \text{ the number of strings } \gamma \text{ of odd length with } \gamma0\prefixeq\sigma \]
Let us consider a number $t\in\mathrm{StepIV}$. 
We define $\gamma_0:=\mathcal{T}[t]$.
If $\ell(\mathcal{T}[t])>0$ then let $\gamma_1,\ldots,\gamma_{\ell(\mathcal{T}[t])}$ be
the strings $\delta$ of odd length with $\delta0\prefixeq\mathcal{T}[s]$ in decreasing order of length.
Thus, for $i=0,\ldots,\ell(\mathcal{T}[t])-1$ the string $\gamma_{i+1}$ is the longest string $\delta$ of odd length with $\delta0\prefixeq\gamma_{i}$.
For $j=0,\ldots,\ell(\mathcal{T}[t])$ we see:
if $s$ is a stage with $\mathrm{trigIV}^{\circ j}(s)=t$ then $\mathcal{T}[s]=\gamma_j$.
And there can be no stage $s$ with $\mathrm{trigIV}^{\circ j}(s)=t$ for any $j>\ell(\mathcal{T}[t])$.
If $s_{\ell(\mathcal{T}[t])}$ is a stage such that 
$\mathrm{trigIV}^{\circ \ell(\mathcal{T}[t])}(s_{\ell(\mathcal{T}[t])})=t$ then the stages
$s_j:=\mathrm{trigIV}^{\circ \ell(\mathcal{T}[t])-j}(s_{\ell(\mathcal{T}[t])})$ for $j=0,\ldots,\ell(\mathcal{T}[s])$ satisfy
$\mathrm{trigIV}^{\circ j}(s_j)=t$ and $\mathcal{T}[s_j]=\gamma_j$, thus, we have the following situation:
\[
\begin{array}{rllllllll}
     t = s_0 & < & s_1 & < & \ldots & < & s_{\ell(\mathcal{T}[t])-1} & < &  s_{\ell(\mathcal{T}[t])}\\
     \mathcal{T}[t] = \gamma_0 & \supseteq & \gamma_10 & \supseteq & \ldots & \supseteq &  \gamma_{\ell(\mathcal{T}[t])-1} 0 & \supseteq & \gamma_{\ell(\mathcal{T}[t])} 0
\end{array}
\]
In particular, we obtain
\begin{equation}
\label{equation:TrigVcirc-1}
 \begin{array}{l}
  \{s\in\mathrm{StepIV} \colon (\exists j \geq 0) \, \mathrm{trigIV}^{\circ j}(s)=t\} \\
  = \{s\in\mathrm{StepIV} \colon (\exists j) \, (0\leq j \leq \ell(\mathcal{T}[t]) \text{ and } \mathrm{trigIV}^{\circ j}(s)=t)\}
	\end{array}
\end{equation}
The set in \eqref{equation:TrigVcirc-1} can be described as the set of all stages $s$ such that $\mathcal{T}[s]$
reaches step (IV) in stage $s$ and such that this fact was triggered by $\mathcal{T}[t]$ in step (IV) of stage $t$.

\begin{lemma}
\label{lemma:triggered-StepV-finite}
For every $t_1\in\mathrm{StepIV}$ the set 
$\{s\in\mathrm{StepIV} \colon (\exists j \geq 0) \, \mathrm{trigIV}^{\circ j}(s)=t_1\}$
is finite.
\end{lemma}

\begin{proof}
First we prove the following

\emph{Claim: The set $\{s\in\mathrm{StepIV} \colon \mathrm{trigIV}(s)=t_1\}$ is finite.}

Let us prove this. 
If there is no string $\gamma$ of odd length with $\gamma0\prefixeq\mathcal{T}[t_1]$ then the set $\{s\in\mathrm{StepIV} \colon \mathrm{trigIV}(s)=t_1\}$ is empty.
Otherwise, let $\gamma$ be the longest string of odd length with $\gamma0\prefixeq\mathcal{T}[t_1]$.
For any element $s$ of the set in the claim we see $\mathcal{T}[s]=\gamma$ and $t_1<s$ and
$q(\gamma)[s']\neq 0$, for every $s'$ with $t_1\leq s' < s$.
Thus, if the set in the claim were infinite then we would have $q(\gamma)[s]\neq 0$ for all $s\geq t_1$.
In particular, we would get $r(\gamma)[s]=r(\gamma)[t_1]$ for all $s\geq t_1$
(remember that initializations do not change $r(\gamma)$).
On the one hand, in any stage $s>t_1$ the absolute value of the number $q(\gamma)$
cannot increase because this could happen only if in this stage a strategy $\tau$ with $\gamma 0\prefixeq\tau$ would reach step (IV).
But this would imply $\gamma 0\prefixeq\tau\prefixeq\mathcal{T}[s]$ and $q(\gamma)[s-1]=0$.
On the other hand, in step (III$_\mathrm{odd}$) of any stage $s$ in the set in the claim
the absolute value of the number $q(\gamma)$ will be decreased by $2^{-r(\gamma)[t_1]}$ or will be set to $0$.
But then, after finitely many steps, it must have reached the value zero. 
This is a contradiction.
We have proved the claim above.

Next, by induction one concludes that for each $j$ with $0\leq j \leq \ell(\mathcal{T}[t_1])$ the set
$\{t\in\mathrm{StepIV} \colon \mathrm{trigIV}^{\circ j}(s)=t_1\}$ is finite.
Now the assertion follows from \eqref{equation:TrigVcirc-1}.
\end{proof}

We define a function $\mathrm{trigult}:\mathrm{StepIV}\to\mathrm{StepIV}$ by
\[ \mathrm{trigult}(t):=\mathrm{trigIV}^{\circ m}(t) \text{ where } m:=\max\{j\in\IN \colon \mathrm{trigIV}^{\circ j}(t) \text{ is defined}\} . \]
This function is well defined because for every $t\in\mathrm{StepIV}$ the value $\mathrm{trigIV}^{\circ 0}(t)$ is defined
and equal to $t$ and, whenever $\mathrm{trigIV}(t)$ is defined, then $\mathrm{trigIV}(t)<t$.
Furthermore, for any $t\in\mathrm{StepIV}$ the strategy $\mathcal{T}[\mathrm{trigult}(t)]$ has even length and reaches step (IV) in stage $\mathrm{trigult}(t)$.
One may say that it is this strategy of even length that in stage $\mathrm{trigult}(t)$ ``ultimately'' triggered that $\mathcal{T}[t]$ reaches step (IV) in stage $t$.
Later we shall formulate two lemmas about the function $\mathrm{trigult}$.

We wish to give some intuition for the following two lemmas before we state them formally.
If a strategy $\sigma$ (it may have even length or odd length) reaches step (IV) in some stage $t_1$ then this means that the value $v_{t_1}\neq 0$ should be added to $x$.
If there is no strategy $\gamma$ of odd length with $\gamma0\prefixeq\sigma$ then this can be done immediately. 
Otherwise, the task to add $v_{t_1}$ is delegated to the longest strategy $\gamma$ of odd length with $\gamma0\prefixeq\sigma$,
and this value is stored in $q(\gamma)[t_1]$. This may be repeated. 
Now, as long as after stage $t_1$ the strategy $\sigma$ is not initialized and does not reach stage (IV) again,
the sum of movements of $x$ triggered in stage $t_1$ and of the stored values $q(\gamma)$ triggered in stage $t_1$ is equal to $v_{t_1}$.
And if $\sigma$ and also some $\gamma$'s of odd length with $\gamma0\prefixeq\sigma$ are initialized at some stage then,
with respect to this sum, this has only the effect that the corresponding values $q(\gamma)$ are set to $0$.

In order to state this more formally, we introduce some more notation.
We define for $t_1,t_2\in\IN$
\[ \mathrm{Trig}(t_1)[t_2] := \{ \gamma \in\Sigma^\mathrm{odd} \colon q(\gamma)[t_2]\neq 0 \text{ and } (\exists j\geq 0) \ \mathrm{trigIV}^{\circ j}(\mathrm{trigd}(\gamma,t_2))=t_1 \} . \]
This set contains all strategies $\gamma$ with odd length and nonzero value $q(\gamma)[t_2]$
such that the stage $t_1$ is one of the stages in which this nonzero value $q(\gamma)[t_2]$ was ``triggered''.
Clearly, if $t_2 < t_1$ or if $\mathcal{T}[t_1]$ does not reach step (IV) in stage $t_1$ then $\mathrm{Trig}(t_1)[t_2]=\emptyset$.
For every $t_2\geq t_1$ and every $\gamma\in\mathrm{Trig}(t_1)[t_2]$ we have $\gamma\in\Sigma^\mathrm{odd}$ and $\gamma0\prefixeq\mathcal{T}[t_1]$.
Hence, if $\ell(\mathcal{T}[t_1])=0$ then $\mathrm{Trig}(t_1)[t_2]=\emptyset$.
We observe that for $t_1,t_2\in\IN$ with $\ell(\mathcal{T}[t_1])>0$
\begin{equation}
\label{equation:ell}
\mathrm{Trig}(t_1)[t_2] = \{ \gamma \in\Sigma^\mathrm{odd} \colon \begin{array}[t]{l}
       q(\gamma)[t_2]\neq 0 \text{ and }  \\
       (\exists j) \, ( 0 \leq j \leq \ell(\mathcal{T}[t_1])-1 \text{ and }  \mathrm{trigIV}^{\circ j}(\mathrm{trigd}(\gamma,t_2))=t_1) \} . 
			\end{array}
\end{equation}
For $t_1\in\IN$ we define
\[ \mathrm{Jumps}(t_1) := \{s \in\IN \colon x_s \neq x_{s-1} \text{ and } (\exists j\geq 0)\, \mathrm{trigIV}^{\circ j}(s)=t_1\} . \]
Note that, if $x_s\neq x_{s-1}$ then $s\in \mathrm{StepIV}$ and there is no $\gamma\in\Sigma^\mathrm{odd}$ with $\gamma0\prefixeq\mathcal{T}[s]$,
hence $\ell(\mathcal{T}[s])=0$.

\begin{lemma}
\label{lemma:Jumps}
\begin{enumerate}
\item
If $t_1\not\in\mathrm{StepIV}$ then $\mathrm{Jumps}(t_1) = \emptyset$.
\item
Every $s\in \mathrm{Jumps}(t_1)$ satisfies $s\geq t_1$.
\item
If $t_1<t_2$ and $\mathcal{T}[t_1] \prefixeq \mathcal{T}[t_2]$ then $\mathrm{Jumps}(t_1) \subseteq \{t_1,\ldots,t_2-1\}$.
\end{enumerate}
\end{lemma}

\begin{proof}
The first assertion is clear.
As for every $s$ and every $j$ we have $\mathrm{trigIV}^{\circ j}(s) \leq s$, we conclude that every $s\in \mathrm{Jumps}(t_1)$
satisfies $s \geq t_1$.
Finally, let us consider $t_1<t_2$ with $\mathcal{T}[t_1] \prefixeq \mathcal{T}[t_2]$.
We wish to show that every $s\in \mathrm{Jumps}(t_1)$ satisfies $s<t_2$.
For the sake of a contradiction, let us assume that there is an $s\in \mathrm{Jumps}(t_1)$ with $s\geq t_2$.
Then there exists a number $k\geq 1$ with $\mathrm{trigIV}^{\circ k} (s)=t_1$.
But we prove by induction over $j$ that $\mathrm{trigIV}^{\circ j}(s)\geq t_2$ for all $j$ such that $\mathrm{trigIV}^{\circ j}(s)$ is defined.
For $j=0$ this is clear because $\mathrm{trigIV}^{\circ 0}(s)=s\geq t_2$.
Let us assume that $s_j:=\mathrm{trigIV}^{\circ j}(s)\geq t_2$ for some $j<k$.
This implies that $s_{j+1}:=\mathrm{trigIV}^{\circ j+1}(s)$ exists.
The strategy $\gamma:=\mathcal{T}[s_j]$ must have odd length and must satisfy $\gamma0\prefixeq\mathcal{T}[t_1]$
because $\mathrm{trigIV}^{\circ k-j}(s_j) = t_1$.
As $\mathcal{T}[t_1]\prefixeq \mathcal{T}[t_2]$,
we must have $q(\gamma)[t_2-1]=0$ which, by Lemma~\ref{lemma:trig-1}, implies that
$s_{j+1}=\mathrm{trigIV}(s_j)=\mathrm{trigd}(\gamma,s_j-1)\geq t_2$.
We have shown that a number $s\geq t_2$ cannot be an element of $\mathrm{Jumps}(t_1)$.
\end{proof}

\begin{lemma}
\label{lemma:triggered-movement}
Let us assume that some strategy $\sigma$ reaches step (IV) in some stage $t_1$.
\begin{enumerate}
\item
Let $t_3\geq t_1$ be a stage such that there is no stage $t_2$ with $t_1 < t_2 \leq t_3$ in which $\sigma$ is initialized or reaches step (IV).
Then
\begin{equation}
\label{eq:movement-not-initialized}
  v_{t_1}
  = \sum_{s\in \mathrm{Jumps}(t_1)\restrictedto (t_3+1)} (x_s - x_{s-1})
  + \sum_{\gamma \in\mathrm{Trig}(t_1)[t_3]} q(\gamma)[t_3] .
\end{equation}
\item
Let $t_3\geq t_1$ be a stage such that there is no stage $t_2$ with $t_1 < t_2 \leq t_3$ in which $\sigma$ reaches step (IV).
Then
\begin{equation}
\label{eq:movement-possibly-initialized}
  \left| v_{t_1} \right|
  \geq \left| \sum_{s\in \mathrm{Jumps}(t_1)\restrictedto (t_3+1)} (x_s - x_{s-1})
  + \sum_{\gamma \in\mathrm{Trig}(t_1)[t_3]} q(\gamma)[t_3] \right| . 
\end{equation}
\item
Furthermore, if $v_{t_1}>0$ then all of the numbers $(x_s - x_{s-1})$ and $q(\gamma)[t_3]$ appearing in the sums in the right hand side of
\eqref{eq:movement-not-initialized} and of \eqref{eq:movement-possibly-initialized} are positive, if $v_{t_1}<0$ then all of these numbers are negative.
\end{enumerate}
\end{lemma}

\begin{proof}
We assume $v_{t_1}>0$. The other case, $v_{t_1}<0$, is treated analogously.
We prove the assertions by induction over $t_3$. We start with $t_3=t_1$.
In stage $t_1$ the strategy $\sigma=\mathcal{T}[t_1]$ reaches step (IV).
If there is no $\gamma$ of odd length with $\gamma0\prefixeq\sigma$ then $t_1\in\mathrm{Jumps}(t_1)\restrictedto (t_1+1)$ and $x_{t_1} = x_{t_1-1}+v_{t_1}$
and $\mathrm{Trig}(t_1)[t_1]=\emptyset$,
hence, the assertions are true for $t_3=t_1$.
If there is a $\gamma$ of odd length with $\gamma0\prefixeq\sigma$, then for the longest $\gamma$ of this kind in step (IV) we obtain
$q(\gamma)[t_1]=v_{t_1}>0$ while for all shorter $\gamma$ of this kind we have $q(\gamma)[t_1]=q(\gamma)[t_1-1]=0$.
Furthermore, we obtain $x_{t_1}=x_{t_1-1}$. Hence, in this case we have $\mathrm{Jumps}(t_1)[t_1]=\emptyset$ and
$\mathrm{Trig}(t_1)[t_1]=\{\gamma\}$ where $\gamma$ is the longest strategy of odd length with $\gamma0\prefixeq\sigma$.
Thus, in this case for $t_3=t_1$ the claimed assertions are true as well.

We come to the induction step. Let us assume that $t_3>t_1$. We set $\tau:=\mathcal{T}[t_3]$. 
We distinguish several cases. 

\emph{Case I:}
$t_3\in\mathrm{Jumps}(t_1)\restrictedto (t_3+1)$.

As $t_1<t_3$ we have $\tau\in\mathrm{Trig}(t_1)[t_3-1]$.
By induction hypothesis (third assertion) $q(\tau)[t_3-1]>0$, hence, $v_{t_3}>0$.
In step (III$_\mathrm{odd}$) of stage $t_3$ we obtain $q(\tau)[t_3]=q(\tau)[t_3-1]-v_{t_3} \geq 0$.
This implies that $\tau\in\mathrm{Trig}(t_1)[t_3]$ if $q(\tau)[t_3]>0$ and $\tau\not\in\mathrm{Trig}(t_1)[t_3]$ if $q(\tau)[t_3]=0$.
We also note that due to $\tau0\prefixeq\mathcal{T}[t_1]$ and $\ell(\tau)=0$ the strategy $\tau$ is the shortest strategy $\gamma$ of odd length with $\gamma0\prefixeq\mathcal{T}[t_1]$.
For all other strategies $\gamma$ of odd length with $\gamma0\prefixeq\mathcal{T}[t_1]$ we obtain
$q(\gamma)[t_3]=q(\gamma)[t_3-1]$, and this implies
$\mathrm{Trig}(t_1)[t_3] \setminus\{\tau\} = \mathrm{Trig}(t_1)[t_3-1] \setminus\{\tau\}$.
Thus, the right one of the two sums in the right hand side of
\eqref{eq:movement-not-initialized} and of \eqref{eq:movement-possibly-initialized} decreases by $v_{t_3}>0$.
Furthermore, $x_{t_3} = x_{t_3-1} + v_{t_3}$, hence, the left one of the two sums in the right hand side of
\eqref{eq:movement-not-initialized} and of \eqref{eq:movement-possibly-initialized}
increases by $v_{t_3}>0$.
Hence, in this case the claimed assertions are true.

Let us assume that Case I does not apply, that is, that $t_3\not\in\mathrm{Jumps}(t_1)\restrictedto (t_3+1)$. 
Then the left one of the two sums in the right hand side of
\eqref{eq:movement-not-initialized} and of \eqref{eq:movement-possibly-initialized} does not change when going from $t_3-1$ to $t_3$.
Therefore, we only need to look at the right one of the two sums.

\emph{Case II:} 
$t_3\not\in\mathrm{Jumps}(t_1)\restrictedto (t_3+1)$ but $t_3\in\mathrm{StepIV}$ and there exists a $j\geq 0$ with $\mathrm{trigIV}^{\circ j}(t_3)=t_1$.

As $t_1<t_3$ we have $j>0$ and $\tau\in\mathrm{Trig}(t_1)[t_3-1]$. By induction hypothesis (third assertion) $q(\tau)[t_3-1]>0$, hence, $v_{t_3}>0$.
We also obtain $\tau0\prefixeq\sigma$.
But we also have $\ell(\tau)>0$. Let $\delta$ be the longest strategy of odd length with $\delta0\prefixeq\tau$.
Remember that the sets $\mathrm{Trig}(t_1)[t_3-1]$ and $\mathrm{Trig}(t_1)[t_3]$ contain only strings $\gamma\in\Sigma^\mathrm{odd}$
with $\gamma0\prefixeq\sigma$.
There are four possible types of such strings, and we observe the following facts about them.
\begin{enumerate}
\item
$q(\tau)[t_3]=q(\tau)[t_3-1]-v_{t_3}\geq 0$ in step (III$_\mathrm{odd}$) of stage $t_3$.
\item
On the one hand, $q(\delta)[t_3-1]=0$ because $\delta0\prefixeq\tau\prefixeq\mathcal{T}[t_3]$,
on the other hand, $q(\delta)[t_3]=v_{t_3}>0$ in step (IV) of stage $t_3$.
\item
For all $\gamma\in\Sigma^\mathrm{odd}$ with $\gamma0\prefixeq\sigma$ that are shorter than $\delta$
and, hence, satisfy $\gamma0\prefixeq\delta\prefixeq\tau=\mathcal{T}[t_3]$, we obtain $q(\gamma)[t_3-1]=0=q(\gamma)[t_3]$.
\item
For all $\gamma\in\Sigma^\mathrm{odd}$ with $\gamma0\prefixeq\sigma$ that are longer than $\tau$
we obtain $\tau0\prefixeq\gamma$. Hence, these strategies are not initialized in stage $t_3$, and we obtain
$q(\gamma)[t_3]=q(\gamma)[t_3-1]$.
\end{enumerate}
We conclude that the right one of the two sums in the right hand side of
\eqref{eq:movement-not-initialized} and of \eqref{eq:movement-possibly-initialized} does not change when going from $t_3-1$ to $t_3$.
In this case the claimed assertions are true as well.

\emph{Case III:}
$t_3\in\mathrm{StepIV}$ but there is no $j\geq 0$ with $\mathrm{trigIV}^{\circ j}(t_3)=t_1$.

We observe $\tau\not\in\mathrm{Trig}(t_1)[t_3-1]$.
Furthermore, then $t_3\not\in\mathrm{Jumps}(t_1)\restrictedto (t_3+1)$.
Remember that the sets $\mathrm{Trig}(t_1)[t_3-1]$ and $\mathrm{Trig}(t_1)[t_3]$ contain only strings $\gamma\in\Sigma^\mathrm{odd}$
with $\gamma0\prefixeq\sigma$.
There are three cases when for a strategy $\gamma\in\Sigma^\mathrm{odd}$ with $\gamma0\prefixeq\sigma$ the value $q(\gamma)[t_3]$ can be different from $q(\gamma)[t_3-1]$:
\begin{enumerate}
\item
When $\gamma=\tau$. Then $\tau$ has odd length and $q(\tau)[t_3-1]\neq 0$ (because $t_3\in\mathrm{StepIV}$).
But $\tau\not\in\mathrm{Trig}(t_1)[t_3-1]$.
We claim that also $\tau\not\in\mathrm{Trig}(t_1)[t_3]$.
If $q(\tau)[t_3]\neq 0$ then this is true because then $\mathrm{trigd}(\tau,t_3)=\mathrm{trigd}(\tau,t_3-1)$.
If $q(\tau)[t_3]=0$ then this is trivial.
\item
When there is a strategy $\delta$ of odd length with $\delta0\prefixeq\tau$ and $\gamma$ is the longest such $\delta$.
On the one hand, $\gamma0\prefixeq\tau=\mathcal{T}[t_3]$ implies $q(\gamma)[t_3-1]=0$, hence, $\gamma\not\in\mathrm{Trig}(t_1)[t_3-1]$.
On the other hand, $\mathrm{trigd}(\gamma,t_3)=t_3$, and as there is no $j\geq 0$ with $\mathrm{trigIV}^{\circ j}(t_3)=t_1$
we conclude that also $\gamma\not\in\mathrm{Trig}(t_1)[t_3]$.
\item
When $\gamma$ is initialized in stage $t_3$. But then $\sigma$ is initialized in stage $t_3$ as well.
Furthermore, then $q(\gamma)[t_3]=0$ and $\gamma\not\in\mathrm{Trig}(t_1)[t_3]$.
If $\gamma\in\mathrm{Trig}(t_1)[t_3-1]$ then by induction hypothesis $q(\gamma)[t_3-1]>0$.
\end{enumerate}
We conclude that, if $\sigma$ is not initialized in stage $t_3$, then the right one of the two sums in the right hand side of
\eqref{eq:movement-not-initialized} and of \eqref{eq:movement-possibly-initialized} does not change when going from $t_3-1$ to $t_3$.
Otherwise, some values in the right one of the two sums are decreased from positive to zero, and the others remain unchanged.
Hence, in Case III the claimed assertions are true as well.

\emph{Case IV:}
$t_3\not\in\mathrm{StepIV}$.
Then the only possibility for a strategy $\gamma\in\Sigma^\mathrm{odd}$ with $\gamma0\prefixeq\sigma$ that $q(\gamma)[t_3]$ can be different from $q(\gamma)[t_3-1]$
is when $\gamma$ is initialized in stage $t_3$. But then $\sigma$ is initialized in stage $t_3$ as well.
So, if $\sigma$ is not initialized in stage $t_3$ then the right one of the two sums in the right hand side of
\eqref{eq:movement-not-initialized} and of \eqref{eq:movement-possibly-initialized} does not change when going from $t_3-1$ to $t_3$.
If some strategy $\gamma\in\mathrm{Trig}(t_1)[t_3-1]$ is initialized in stage $t_3$
then we obtain $q(\gamma)[t_3-1]> 0 = q(\gamma)[t_3]$.
Thus, in this case some values in the right one of the two sums may be decreased from positive to zero, and the others remain unchanged.
Hence, in Case IV the claimed assertions are true as well.
\end{proof}

\subsection{About strategies of even length}
\label{subsection:6}

Let us now turn to strategies $\sigma$ of even length.
We define
\[ \mathit{D} := \{d\in\IN\colon 
   (\exists \sigma\in\Sigma^\mathrm{even}) \, (\exists t>0) \,
	 ( d=d(\sigma)[t-1] \text{ and } 
	 \sigma \text{ reaches step (II$_\mathrm{even}$) in stage } t)\}.
\]
Let us consider a number $d\in \mathit{D}$.
There is exactly one $\sigma\in\Sigma^\mathrm{even}$ such that there exists a $t>0$ with $d=d(\sigma)[t-1]$.
Let us call it $\sigma_d$. 
Let $(T_j^{(d)})_j$ be the finite or infinite increasing sequence of all stages $t>0$ with $d=d(\sigma_d)[t-1]$ and $\sigma_d\prefixeq\mathcal{T}[t]$ and in which $\sigma_d$ reaches step (II$_\mathrm{even}$).
This sequence is not empty and satisfies $T_0^{(d)}>0$.

For real numbers $a,b,c$, let us say that the number $b$ \emph{lies between $a$ and $c$}
if $\min\{a,c\}\leq b \leq \max\{a,c\}$.

\begin{corollary}
\label{corollary:even-triggered-movement}
Let $d\in \mathit{D}$. 
\begin{enumerate}
\item
For every $k$ such that $\sigma_d$ reaches step (IV) in stage $T_k^{(d)}$ and such that
$T_{k+1}^{(d)}$ exists we have $\mathrm{Jumps}(T_k^{(d)}) \restrictedto T_{k+1}^{(d)} = \mathrm{Jumps}(T_k^{(d)})$.
\item
For every $k$ such that $\sigma_d$ reaches step (IV) in stage $T_k^{(d)}$ and such that $T_{k+1}^{(d)}$ exists
\[ \sum_{s \in\mathrm{Jumps}(T_k^{(d)})} (x_s - x_{s-1}) 
  = \sum_{s \in\mathrm{Jumps}(T_k^{(d)})\restrictedto T_{k+1}^{(d)}} (x_s - x_{s-1}) 
  = A(d,T_k^{(d)},h_{d}(k),x_{T_k^{(d)}-1}) \]
and
\[ \sum_{s \in\mathrm{Jumps}(T_k^{(d)})} \left| x_s - x_{s-1} \right| 
  = \sum_{s \in\mathrm{Jumps}(T_k^{(d)})\restrictedto T_{k+1}^{(d)}} \left| x_s - x_{s-1} \right| 
  = \left| A(d,T_k^{(d)},h_{d}(k),x_{T_k^{(d)}-1}) \right| . \]
\item
Let $t\in\IN$. 
Let us assume that there is a number $k$ with 
$T_k^{(d)} < t$ and let $k$ be the largest such number.
Furthermore, let us assume that $\sigma_d$ reaches step (IV) in stage $T_k^{(d)}$.
Then the sum $\sum_{s \in\mathrm{Jumps}(T_k^{(d)})\restrictedto t} (x_s - x_{s-1})$
lies between the numbers $0$ and $A(d,T_k^{(d)},h_{d}(k),x_{T_k^{(d)}-1})$,
and the sum $\sum_{s \in\mathrm{Jumps}(T_k^{(d)})\restrictedto t} \left| x_s - x_{s-1} \right|$
lies between the numbers $0$ and $\left|A(d,T_k^{(d)},h_{d}(k),x_{T_k^{(d)}-1})\right|$.
\end{enumerate}
\end{corollary}

\begin{proof}
In the following for simplicity we write $T_k$ instead of $T_k^{(d)}$.
\begin{enumerate}
\item
If $\sigma_d$ reaches step (IV) in stage $T_k$ then $\mathcal{T}[T_k]=\sigma$.
If $T_{k+1}$ exists then $\sigma\prefixeq \mathcal{T}[T_{k+1}]$.
Thus, the assertion follows from Lemma~\ref{lemma:Jumps}.
\item
Let us assume that $\sigma_d$ reaches step (IV) in stage $T_k$ and that $T_{k+1}$ exists.
The equation 
\[ \sum_{s \in\mathrm{Jumps}(T_k)} (x_s - x_{s-1})
 = \sum_{s \in\mathrm{Jumps}(T_k)\restrictedto T_{k+1}} (x_s - x_{s-1}) \]
follows from the first assertion of the corollary; same for the corresponding equation with terms $|x_s - x_{s-1}|$.
We see $v_{T_k}=A(d,T_k,h_{d}(k),x_{T_k-1}) \neq 0$.
As $T_{k+1}$ is defined, $\sigma_d\prefixeq\mathcal{T}[T_{k+1}]$.
This implies that $q(\gamma)[T_{k+1}-1]=0$ for all strategies $\gamma$ of odd length with $\gamma0\prefixeq\sigma_d$.
By the first assertion of Lemma~\ref{lemma:triggered-movement} we obtain
\[  \sum_{s \in\mathrm{Jumps}(T_k)\restrictedto T_{k+1}} (x_s - x_{s-1}) 
  = v_{T_k} = A(d,T_k,h_{d}(k),x_{T_k-1}) , \]
and by the first and the third assertion of Lemma~\ref{lemma:triggered-movement} we obtain
\[  \sum_{s \in\mathrm{Jumps}(T_k)\restrictedto T_{k+1}} |x_s - x_{s-1}|
  = |v_{T_k}| = \left|A(d,T_k,h_{d}(k),x_{T_k-1})\right| . \]
\item
Again we have $v_{T_k}=A(d,T_k,h_{d}(k),x_{T_k-1}) \neq 0$.
The assertion follows from the second and the third assertion of Lemma~\ref{lemma:triggered-movement}.
\qedhere
\end{enumerate}
\end{proof}

\begin{lemma}
\label{lemma:triggered-2}
Let $d\in D$.
Let $j\in\dom(T^{(d)})$ be a number such that $\sigma_d$ reaches step (IV) in stage $t_1:=T^{(d)}_j$.
Let $j'\in\dom(T^{(d)})$ be a number with $j'>j$ such that $\sigma_d$ reaches step (II$_\mathrm{even}$) in stage $t_2:=T^{(d)}_{j'}$
and such that $h_{d}(j')=h_{d}(j)+1$.
Let $s$ be a number with $t_1<s<t_2$ and $s\in\mathrm{StepIV}$.
Then either $t_1=\mathrm{trigult}(s)$ or $t_1<\mathrm{trigult}(s)$ and
$\mathcal{T}[\mathrm{trigult}(s)]$ is either initialized in stage $t_1$ or satisfies $\mathcal{T}[\mathrm{trigult}(s)]0 \prefixeq \sigma_d$. 
\end{lemma}

\begin{proof}
First, we observe that the fact that $\sigma_d$ reaches step (IV) in stage $t_1=T^{(d)}_j$ implies
$A(d,T^{(d)}_j,h_d(j),x_{T^{(d)}_j-1})\neq 0$.

Let us consider the case when the strategy $\gamma:=\mathcal{T}[s]$ has odd length.
Then $s':=\mathrm{trigIV}(s)$ is defined.
We claim that then $t_1\leq s'$.
Let us prove this.
We cannot have $\gamma0<_L\sigma_d$ because otherwise $\sigma_d$ would be initialized in stage $s$
leading to $d(\sigma_d)[s]>d$ in contradiction to $s<t_2$ and $d=d(\sigma_d)[t_2]$.
If $\gamma0\prefixeq\sigma_d$ then $q(\gamma)[t_1-1]=0$, hence by Lemma~\ref{lemma:trig-1} we obtain $t_1-1<s'$, hence $t_1\leq s'$.
We cannot have $\gamma=\sigma_d$ because $\gamma$ has odd length and $\sigma_d$ has even length.
We cannot have $\sigma_d0\prefixeq\gamma$:
if $|\sigma_d|$ is divisible by $4$ then this is true by Lemma~\ref{lemma:div4-no-0}
because then we would have $\sigma_d0\prefixeq\gamma = \mathcal{T}[s]$;
if $|\sigma_d|$ is divisible by $2$ but not by $4$ then $\sigma_d0\prefixeq\gamma$ would imply that
$\sigma_d$ reaches step (II$_\mathrm{even}$) in stage $s$, there is a $j''$ with $s=T^{(d)}_{j''}$, and
we would have $A(d,s,h_d(j''),x_{s-1})\neq 0$ in contradiction to $T^{(d)}_j = t_1 < s=T^{(d)}_{j''} < t_2 = T^{(d)}_{j'}$ 
and $h_{d}(j')=h_{d}(j)+1$.
Finally, if $\sigma_d0 <_L \gamma$ then $\gamma$ is initialized in stage $t_1$ leading to $q(\gamma)[t_1]=0$.
By Lemma~\ref{lemma:trig-1} we obtain $t_1<s'$.
We have proved the claim $t_1\leq s'$.

By induction we obtain $t_1 \leq \mathrm{trigIV}^{\circ k}(s)$ for all $k$ such that $\mathrm{trigIV}^{\circ k}(s)$ is defined.
Ultimately, we obtain $t_1 \leq \mathrm{trigult}(s) \leq s$.
We set $\delta:=\mathcal{T}[\mathrm{trigult}(s)]$.
Then $\delta$ has even length and reaches step (IV) in stage $\mathrm{trigult}(s)$.
If $t_1=\mathrm{trigult}(s)$ then we are done. Let us assume $t_1<\mathrm{trigult}(s)$.
We claim that then $\delta$ is either initialized in stage $t_1$ or satisfies $\delta0 \prefixeq \sigma_d$. 
We cannot have $\delta0<_L\sigma$ because otherwise $\sigma$ would be initialized in stage $\mathrm{trigult}(s)$,
and we have seen above already that $\sigma$ cannot be initialized in any stage properly between $t_1$ and $t_2$.
We cannot have $\delta=\sigma_d$ because then $\sigma_d$ would reach step (IV) in stage $s$,
hence, $\sigma_d$ would reach step (II$_\mathrm{even}$) in stage $s$, there would be a $j''$ with $s=T^{(d)}_{j''}$, and
we would have $A(d,s,h_d(j''),x_{s-1})\neq 0$ in contradiction to $T^{(d)}_j = t_1 < s=T^{(d)}_{j''} < t_2 = T^{(d)}_{j'}$ 
and $h_{d}(j')=h_{d}(j)+1$.
We cannot have $\sigma_d0\prefixeq\delta$:
if $|\sigma_d|$ is divisible by $4$ then this is true by Lemma~\ref{lemma:div4-no-0}
because then we would have $\sigma_d0\prefixeq\delta = \mathcal{T}[s]$;
if $|\sigma_d|$ is divisible by $2$ but not by $4$ then $\sigma_d0\prefixeq\delta$ would imply that
$\sigma_d$ reaches step (II$_\mathrm{even}$) in stage $s$, there is a $j''$ with $s=T^{(d)}_{j''}$, and
we would have $A(d,s,h_d(j''),x_{s-1})\neq 0$ in contradiction to $T^{(d)}_j = t_1 < s=T^{(d)}_{j''} < t_2 = T^{(d)}_{j'}$ 
and $h_{d}(j')=h_{d}(j)+1$.
Finally, if $\sigma_d0 <_L \delta$ then $\delta$ is initialized in stage $t_1$.
Or we can have $\delta0\prefixeq\sigma_d$.
\end{proof}

For a finite or infinite increasing sequence $(S_j)_j$ of natural numbers and a natural number $t$
let $(S_j)_j|^{<t}$ be the finite initial segment of $(S_j)_j$ consisting of 
exactly the elements $S_j< t$ of the sequence $(S_j)_j$.
Furthermore, in the following by $(x_t)_t$ we shall always mean the sequence constructed above.
In particular, whenever in the following we speak about a ``$(d,j)$-good'' sequence resp. a ``$d$-good'' sequence,
we always mean a ``$(d,j,(x_t)_t)$-good'' sequence resp. a ``$(d,(x_t)_t)$-good'' sequence, where $(x_t)_t$ is the sequence constructed above.
For $d\in \mathit{D}$ we define
\[ J^{(d)} := \left\{j\in\dom(T^{(d)}) \colon T_j^{(d)}\in\mathrm{StepIV} \right\}
\]
and 
\[ \mathrm{Jumps}^{(d)} := \bigcup_{j \in J^{(d)}} \mathrm{Jumps}(T^{(d)}_j). \]

\begin{lemma}
\label{lemma:even-movement}
Let us assume that the assumptions of one of the first three cases in Theorem~\ref{theorem:master} are true
(not the fourth case).
\begin{enumerate}
\item
For every $t_3\in\IN$ and every $d\in \mathit{D}$ the finite sequence $(T_j^{(d)})_j|^{<t_3}$
is either empty or a $d$-good sequence.
\item
For every $t\in\IN$ and every $d\in \mathit{D}$
\[ \left| \sum_{s \in\mathrm{Jumps}^{(d)}\restrictedto t} (x_s - x_{s-1}) \right| \leq 2^{-d}. \]
\end{enumerate}
\end{lemma}

\begin{proof}
We prove this by induction over $t_3$.
Let us fix some $d\in \mathit{D}$ and some $t_3\in\IN$.
In both parts of the following proof we write $T_k$ instead of $T_k^{(d)}$.

\begin{enumerate}
\item
Let us fix a number $j\in\dom((T_k)_k|^{<t_3})$ with $A(d,T_j,h_{d}(j),x_{T_j-1})\neq 0$.
In the case when $\iota(d)$ is odd we additionally assume that $j \not\in \mathrm{First}( d,(T_k)_k|^{<t_3} )$.
It is sufficient to show that the finite sequence $(T_k)_k|^{<t_3}$ is $(d,j)$-good.
Let us fix a number $j'\in\dom((T_k)_k|^{<t_3})$ with $h_{d}(j')=h_{d}(j)+1$.
Let $t_1:=T_j$ and $t_2:=T_{j'}$. Note that $t_1<t_2<t_3$.
Then we have $d=d(\sigma_d)[t_1-1]=d(\sigma_d)[t_2-1]$.
In any stage $t$ with $t_1 \leq t < t_2$ the strategy $\sigma_d$ is not initialized
(otherwise the value $d(\sigma_d)$ would be increased in stage $t$ in contradiction to
$d(\sigma_d)[t_2-1]=d=d(\sigma_d)[t_1-1]$).
The strategy $\sigma_d$ reaches step (IV) in stage $t_1$
(it reaches step (II$_\mathrm{even}$) because $t_1=T_j$;
then it reaches step (III$_\mathrm{even}$) because either $\iota(d)$ is even or
$\iota(d)$ is odd and in this case our additional assumption $j \not\in \mathrm{First}( d,(T_k)_k|^{<t_3} )$ above implies this;
then it reaches step (IV) because $A(d,T_j,h_{d}(j),x_{T_j-1})\neq 0$), hence $\sigma_d=\mathcal{T}[t_1]$.
But $\sigma_d$ does not reach step (IV) in any stage $t$ with $t_1<t<t_2$ (because of $h_{d}(j')=h_{d}(j)+1$).
We see
\[    x_{T_{j'}-1} - x_{T_j-1} = x_{t_2-1} - x_{t_1-1} = \sum_{t=t_1}^{t_2-1} (x_t - x_{t-1}) . \]
Let us consider a number $t$ with $t_1\leq t < t_2$ and $x_t\neq x_{t-1}$.
Then $\mathcal{T}[t]$ reaches step (IV) in stage $t$.
If $t\not\in\mathrm{Jumps}(t_1)\restrictedto t_2$ then $\mathrm{trigult}(t)\neq t_1$,
hence, by Lemma~\ref{lemma:triggered-2}, $t_1 < \mathrm{trigult}(t)$
and $\mathcal{T}[\mathrm{trigult}(t)]$ is either initialized in stage $t_1$ or satisfies $\mathcal{T}[\mathrm{trigult}(t)]0\prefixeq\sigma_d$.
We obtain $t\in\mathrm{Jumps}^{(e)}\restrictedto t_2$ where 
\[ e := d(\mathcal{T}[\mathrm{trigult}(t)])[\mathrm{trigult}(t)-1] . \] 
Then $e\in D^\mathrm{init}(t_1,t_2) \cup D^\mathrm{cancel}(t_1,t_2)$ where we define
\begin{eqnarray*}
  D^\mathrm{init}(t_1,t_2) &:=& \begin{array}[t]{l} \{d' \in \mathit{D} \colon 
     (\exists \tau\in\Sigma^\mathrm{even}) \, (\tau \text{ is initialized in stage } t_1 \text{ and} \\
     (\exists s) \, (t_1 < s < t_2 \text{ and } d'=d(\tau)[s-1] ))\},
		\end{array} \\
  D^\mathrm{cancel}(t_1,t_2) &:=& \begin{array}[t]{l} \{d' \in \mathit{D} \colon 
     (\exists \tau\in\Sigma^\mathrm{even}) \, (\tau0\prefixeq\sigma_d \text{ and} \\
     (\exists s) \, (t_1 < s < t_2 \text{ and } d'=d(\tau)[s-1] ))\}.
		\end{array} 
\end{eqnarray*}
We obtain
\begin{align*}
   &x_{T_{j'}-1} - x_{T_j-1} 
	= x_{t_2-1} - x_{t_1-1} \\
	&= \sum_{t\in\mathrm{Jumps}(t_1)\restrictedto t_2} (x_t - x_{t-1})
	    + \sum_{e\in D^\mathrm{init}(t_1,t_2) \cup D^\mathrm{cancel}(t_1,t_2)}
			   \sum_{t\in \mathrm{Jumps}^{(e)} \land t_1 \leq t < t_2} (x_t - x_{t-1}) 
\end{align*}
By Corollary~\ref{corollary:even-triggered-movement} we have
\[ \sum_{t\in\mathrm{Jumps}(t_1)\restrictedto t_2} (x_t - x_{t-1})
	= A(d,T_j,h_{d}(j),x_{T_j-1}) . \]
For $e\in D^\mathrm{init}(t_1,t_2)$ we have
$\{t\in\IN\colon t \in \mathrm{Jumps}^{(e)} \land t_1 \leq t < t_2\} = \mathrm{Jumps}^{(e)}\restrictedto t_2$.
Furthermore, by induction hypothesis for every $e\in \mathit{D}$ we have
\[ \left| \sum_{t\in \mathrm{Jumps}^{(e)}\restrictedto t_2} (x_t - x_{t-1}) \right|
 \leq 2^{-e} . \]
We claim that for any $e \in D^\mathrm{cancel}(t_1,t_2)$
\begin{equation}
\label{claim:cancel}
  \sum_{t \in \mathrm{Jumps}^{(e)} \land t_1 \leq t < t_2} (x_t - x_{t-1}) \\ = 0 . 
\end{equation}
Let us consider an $e \in D^\mathrm{cancel}(t_1,t_2)$.
There exist a $\tau$ and an $s$ with $t_1 \leq s < t_2$ and $e=d(\tau)[s-1]$. This $\tau$ must be $\sigma_e$.
We observe $\sigma_e0\prefixeq\sigma_d=\mathcal{T}[t_1]$ and
$\sigma_e0\prefixeq\sigma_d\prefixeq \mathcal{T}[t_2]$.
We see that $\iota(e)$ must be odd by Lemma~\ref{lemma:div4-no-0}.
For $m\in\IN$ we define
\[ J_m^{(e)} := \left\{k \in \dom(T^{(e)}) \colon \left| A(e,T^{(e)}_k,h_e(k),x_{T^{(e)}_k-1})\right| = 2^{-m} \right\} . \]
The step (II$_\mathrm{even}$) in the construction and \eqref{equation:cancellation} imply that for every number $m$ such that
$J_m^{(e)} \cap \{t_1,\ldots, t_2-1\} \neq \emptyset$ we have $J_m^{(e)}\subseteq \{t_1,\ldots, t_2-1\}$.
Only for the first number $j\in J_m^{(e)}$ the strategy $\sigma_e$ does not reach step (IV) in step $T^{(e)}_j$.
In all other stages $T^{(e)}_j$ with $j\in J_m^{(e)}\setminus \{\min J_m^{(e)}\}$ it reaches step (IV) in step $T^{(e)}_j$.
With Corollary~\ref{corollary:even-triggered-movement} and \eqref{equation:cancellation} we obtain
\begin{eqnarray*}
  \sum_{t \in \mathrm{Jumps}^{(e)} \land t_1 \leq t < t_2} (x_t - x_{t-1}) 
	&=& \sum_{m \colon J_m^{(e)} \cap \{t_1,\ldots, t_2-1\} \neq \emptyset} 
	   \sum_{j\in J_m^{(e)}\setminus \{\min J_m^{(e)}\}}
		 \sum_{t\in \mathrm{Jumps}(T^{(e)}_j)} (x_t - x_{t-1}) \\
	&=& \sum_{m \colon J_m^{(e)} \cap \{t_1,\ldots, t_2-1\} \neq \emptyset} 
	   \sum_{j\in J_m^{(e)}\setminus \{\min J_m^{(e)}\}} A(e,T^{(e)}_j,h_e(j),x_{T^{(e)}_j-1}) \\
	&=& \sum_{m \colon J_m^{(e)} \cap \{t_1,\ldots, t_2-1\} \neq \emptyset}  0 \\
	&=& 0 .
\end{eqnarray*}
Thus, \eqref{claim:cancel} is true.
Putting all of this together we obtain 
\begin{eqnarray*}
 \lefteqn{\left| x_{T_{j'}-1} - x_{T_j-1}  - A(d,T_j,h_{d}(j),x_{T_j-1}) \right|} && \\
   &\leq & \sum_{e\in D^\mathrm{init}(t_1,t_2)} 2^{-e} \\
	 &\leq & 2^{-(1 - \log_2(|v_{t_1}|))} \\
	 &= & \frac{1}{2} \cdot |v_{t_1}|
	 = \frac{1}{2} \cdot A(d,T_j,h_{d}(j),x_{T_j-1}) ,
\end{eqnarray*}
where in the third line we used that for $e\in D^\mathrm{init}(t_1,t_2)$ we have $e>2 - \log_2(|v_{t_1}|)$.
Hence, $(T_j^{(d)})_j|^{<t_3}$ is a $d$-good sequence.
\item
If $\mathrm{Jumps}^{(d)}\restrictedto t=\emptyset$ then the assertion is trivial.
Let us assume that $\mathrm{Jumps}^{(d)}\restrictedto t\neq\emptyset$.
Let $k$ be the largest number in $J^{(d)}$ 
with $T_k \leq t$.
By Corollary~\ref{corollary:even-triggered-movement} the sum
\[ \sum_{s \in\mathrm{Jumps}^{(d)}\restrictedto t} (x_s - x_{s-1})
  = \sum_{j \in J^{(d)} \restrictedto k} \sum_{s \in\mathrm{Jumps}(T_j)\restrictedto T_{j+1}} (x_s - x_{s-1})
	   + \sum_{s \in\mathrm{Jumps}(T_k)\restrictedto t} (x_s - x_{s-1})
\]
lies between the following two sums:
\begin{itemize}
\item
$\sum_{j \in J^{(d)} \restrictedto k} A(d,T_j,h_{d}(j),x_{T_j-1})
  = \sum_{j\in\{0,\ldots,k-1\}\setminus \mathrm{First}( d,(T_i)_i|^{<t})} A(d,T_j,h_{d}(j),x_{T_j-1})$ and 
\item
$\sum_{j \in J^{(d)} \restrictedto k+1} A(d,T_j,h_{d}(j),x_{T_j-1})
  = \sum_{j\in\{0,\ldots,k\}\setminus \mathrm{First}( d,(T_i)_i|^{<t})} A(d,T_j,h_{d}(j),x_{T_j-1})$
\end{itemize}
By the first assertion of this lemma  the finite sequence $(T_j)_j|^{<t}$ is $d$-good.
By \eqref{equation:bound-on-sum} the absolute value of each of the previous two sums is bounded by $2^{-d}$.
Then this applies to the sum $\sum_{s \in\mathrm{Jumps}^{(d)} \restrictedto t} (x_s - x_{s-1})$ as well.
\qedhere
\end{enumerate}
\end{proof}

\begin{corollary}
\label{corollary:Jumps-interval-bound}
Let us assume that the assumptions of one of the first three cases in Theorem~\ref{theorem:master} are true
(not the fourth case).
\begin{enumerate}
\item
For every $d\in \mathit{D}$ the sequence $(T_j^{(d)})_j$ is a $d$-good sequence.
\item
For every $t_1,t_2\in\IN$ with $t_1<t_2$ and every $d\in \mathit{D}$
\[ \left| \sum_{s \in\mathrm{Jumps}^{(d)}\land t_1 \leq s < t_2} (x_s - x_{s-1}) \right| \leq 2^{-d+1}. \]
\end{enumerate}
\end{corollary}

\begin{proof}
\begin{enumerate}
\item
This follows directly from Lemma~\ref{lemma:even-movement}.1.
\item
The assertion follows from Lemma~\ref{lemma:even-movement}.2:
\begin{eqnarray*}
  \lefteqn{\left| \sum_{s \in\mathrm{Jumps}^{(d)}\land t_1 \leq s < t_2} (x_s - x_{s-1}) \right|} && \\
	&=& \left| \sum_{s \in\mathrm{Jumps}^{(d)}\restrictedto t_2} (x_s - x_{s-1}) 
	  - \sum_{s \in\mathrm{Jumps}^{(d)}\restrictedto t_1} (x_s - x_{s-1}) \right| \\
	&\leq & \left| \sum_{s \in\mathrm{Jumps}^{(d)}\restrictedto t_2} (x_s - x_{s-1})\right|
	  + \left| \sum_{s \in\mathrm{Jumps}^{(d)}\restrictedto t_1} (x_s - x_{s-1}) \right| \\
	 &\leq & 2 \cdot 2^{-d} .
\end{eqnarray*}
\end{enumerate}
\end{proof}

\begin{lemma}
\label{lemma:fourthcase}
Let us assume that the assumptions of the fourth case in Theorem~\ref{theorem:master} are true.
\begin{enumerate}
\item
For every $d\in \mathit{D}$ the sequence $(T_j^{(d)})_j$ is a $d$-left-good sequence.
\item
For every $t_1,t_2\in\IN$ with $t_1<t_2$ and every $d\in \mathit{D}$ 
\[ \sum_{s \in\mathrm{Jumps}^{(d)}\land t_1 \leq s < t_2} (x_s - x_{s-1}) \leq 2^{-d}. \]
\end{enumerate}
\end{lemma}

\begin{proof}
Let us fix a number $d\in\mathit{D}$. In this proof we write $T_j$ instead of $T^{(d)}_j$.
\begin{enumerate}
\item
Let us fix a number $j\in\dom(T)$ with $A(d,T_j,h_{d}(j),x_{T_j-1})\neq 0$,
with $j\not\in\mathrm{First}( d,(T_k)_k )$, and with $j+1\in\dom(T)$.
Then $\sigma_d$ reaches step (IV) in stage $T_j$.
We wish to show that \eqref{equation:i-d-leftgood} is true.
By Corollary~\ref{corollary:even-triggered-movement}
every number $s\in \mathrm{Jumps}(T_j)$ satisfies $T_j \leq s < T_{j+1}$.
We remember that the sequence $(x_t)_t$ is nondecreasing.
Using Corollary~\ref{corollary:even-triggered-movement} again we obtain
\[ x_{T_{j+1}-1} - x_{T_j-1} 
   = \sum_{s=T_j}^{T_{j+1}-1} (x_s - x_{s-1})
   \geq \sum_{s \in \mathrm{Jumps}(T_j) \restrictedto T_{j+1}} (x_s - x_{s-1})
	 = A(d,T_j,h_d(j),x_{T_j-1}) . \]
Thus, \eqref{equation:i-d-leftgood} is true.
\item
If there is a number $j\in\dom(T)$ with $t_2\leq T_j$ then let $J$ be the smallest such $j$.
Otherwise let $J:=|\dom(T)|-1$. In both cases we obtain
\begin{eqnarray*} \sum_{s \in\mathrm{Jumps}^{(d)}\land t_1 \leq s < t_2} (x_s - x_{s-1})
   &\leq & \sum_{s \in\mathrm{Jumps}^{(d)} \restrictedto t_2} (x_s - x_{s-1}) \\
	 &\leq & \sum_{j \in J^{(d)} \land j \leq J} A(d,T_j,h_d(j),x_{T_j-1})
	 \leq 2^{-d} ,
\end{eqnarray*}
where the second inequality is a consequence of Corollary~\ref{corollary:even-triggered-movement}
and the third inequality follows from \eqref{equation:bound-on-sum}.
\qedhere
\end{enumerate}
\end{proof}

\subsection{General lemmas concerning, e.g., only finitely many initializations of the strategies on the true path and convergence of the constructed sequence}
\label{subsection:7}

The following lemma is similar to Lemmas~\ref{lemma:div4-no-0} and \ref{lemma:odd-infinitely-oftenIV}.

\begin{lemma}
\label{lemma:4mod2-0}
Let $\sigma$ have length divisible by $2$ but not by $4$.
If $\sigma$ is initialized only finitely often but there are infinitely many $t$ such that 
$\sigma$ reaches step (IV) in stage $t$ then
there are infinitely many $t$ with $\sigma 0 \prefixeq \mathcal{T}[t]$.
\end{lemma}

\begin{proof}
Let us assume that $\sigma$ is initialized only finitely often
and that $\sigma$ reaches step (IV) in infinitely many stages.
This cannot happen in the third case of Theorem~\ref{theorem:master} 
because in that case we have $A(d_0,\ldots)=0$ for any $d_0$ such that there exists a number $t$ with $d(\sigma)[t]=d_0$.
There exists a number $t_0>0$ such that, for all $t\geq t_0$, $d(\sigma)[t-1]=d(\sigma)[t_0-1]$.
We set $d_0:=d(\sigma)[t_0-1]$.
The sequence $(T^{(d_0)}_j)_j$ is computable.
For simplicity of notation, in the rest of the proof we write $(T_j)_j$ instead of $(T^{(d_0)}_j)_j$.
By Corollary~\ref{corollary:Jumps-interval-bound} this sequence is $(d_0,(x_t)_t)$-good in the first two cases of Theorem~\ref{theorem:master},
and by Lemma~\ref{lemma:fourthcase} it is $(d_0,(x_t)_t)$-left-good in the fourth case of Theorem~\ref{theorem:master}.
It follows that condition \eqref{equation:computable-convergence} is satisfied,
in the fourth case of Theorem~\ref{theorem:master} by assumption,
in the first and second case of Theorem~\ref{theorem:master} by Lemma~\ref{lemma:cancel-comp-conv}
because the conditions \eqref{equation:range-powers-of-2}, \eqref{equation:cancellation},
and \eqref{equation:Jm-finite} are satisfied by assumption.
In each of the infinitely many stages $t$ in which $\sigma$ reaches step (IV) and satisfies $d(\sigma)[t-1]=d_0$,
the strategy $\sigma$ must first reach step (II$_\mathrm{even}$), hence, there must be a $j$ with $t=T_j$,
and we must have $A(d_0,T_j,h_{d_0}(j),x_{T_j-1})\neq 0$.
The condition \eqref{equation:computable-convergence} implies 
 $\lim_{j\to\infty} A(d_0,T_j,h_{d_0}(j),x_{T_j-1}) = 0$.
In combination with the fact that there are infinitely many $j$ with $A(d_0,T_j,h_{d_0}(j),x_{T_j-1})\neq 0$,
we conclude that there are infinitely many $m$ such that there exists a $j$ with
$|A(d_0,T_j,h_{d_0}(j),x_{T_j-1})|=2^{-m}$.
For each such $m$ there is a smallest $j$ such that $|A(d_0,T_j,h_{d_0}(j),x_{T_j-1})|=2^{-m}$.
But in these stages $T_j$ the algorithm goes from stage (II$_\mathrm{even}$) to step (I) and we obtain
$\sigma0\prefixeq\mathcal{T}[T_j]$.
\end{proof}

We turn to assertions concerning the prefixes of the true path $\mathcal{T}$.

\begin{lemma}
\label{lemma:finitely-many-initializations}
Any strategy $\tau \subseteq \mathcal{T}$ is initialized only finitely often.
\end{lemma}

\begin{proof}
We prove the assertion by induction over the length of $\tau$.
Let us fix some $\tau \subseteq \mathcal{T}$ and let us assume that the assertion is true for all $\sigma \subsetneq \tau$.
Let $s>0$ be so large such that no $\sigma$ with $\sigma <_L \tau$ is applied in any stage $t\geq s$
(see Lemma~\ref{lemma:infinite-injury-applied})
and such that no $\sigma$ with $\sigma \subsetneq \tau$ is initialized in any stage $t\geq s$.

In any stage $t>0$ any strategy is initialized at most once.
We claim that there are only finitely many $t$ such that $\tau$ is initialized in stage $t$.
For the sake of a contradiction, let us assume that there are infinitely many stages $t$ in which $\tau$ is initialized.
In any stage $t>0$ initializations take place only in step (I) or in step (IV), and then stage $t+1$ is reached.
If there were a stage $t\geq s$ such that $\tau$ is initialized in step (I) of stage $t$ then this would imply $\mathcal{T}[t] <_L \tau$,
which is impossible by our choice of $s$.
Thus, there must be infinitely many stages $t$ such that $\tau$ is initialized in step (IV) of stage $t$.
For each of these $t$ the strategy $\mathcal{T}[t]$ reaches step (IV) in stage $t$,
and we must have $\mathcal{T}[t]0 <_L \tau$.
It is sufficient to consider the case $t\geq s$.
Then by our choice of $s$ we have $\mathcal{T}[t] \not<_L \tau$, hence, we obtain $\mathcal{T}[t]1 \prefixeq \tau$.
There are only finitely many binary strings $\sigma$ with $\sigma1 \prefixeq\tau$.
Hence, there must be one string $\sigma$ with $\sigma1\prefixeq\tau$ such that 
for infinitely many $t\geq s$ we have $\sigma=\mathcal{T}[t]$ and $\sigma$ reaches step (IV) in stage $t$.

Let $|\sigma|$ be odd. Then in these stages $t$ the strategy $\sigma$ must first reach step (III$_\mathrm{odd}$).
By Lemma~\ref{lemma:odd-infinitely-oftenIV}
there are infinitely many $t'$ with $\sigma0 \subseteq \mathcal{T}[t']$.
We obtain $\sigma0 \subseteq \mathcal{T}$,
and this is not compatible with $\sigma1 \subseteq \tau$ and $\tau \subseteq \mathcal{T}$.

Let $|\sigma|$ be divisible by $2$ but not by $4$. By induction hypothesis $\sigma$ is initialized only finitely often.
As by assumption $\sigma$ reaches step (IV) in infinitely many stages,
by Lemma~\ref{lemma:4mod2-0}
there are infinitely many $t'$ with $\sigma 0 \prefixeq \mathcal{T}[t']$.
We obtain $\sigma0 \subseteq \mathcal{T}$,
and this is not compatible with $\sigma1 \subseteq \tau$ and $\tau \subseteq \mathcal{T}$.

Let $|\sigma|$ be divisible by $4$. Then, as $\sigma$ is initialized only finitely often by induction hypothesis,
there must exist a number $t_0$ such that $d(\sigma)[t]=d(\sigma)[t_0]=:d$ for all $t\geq t_0$.
Then $\sigma=\sigma_d$, and there are infinitely many $j$ such that $\sigma_d$ reaches step (IV) in stage $T^{(d)}_j$.
Hence, the function $h_d$ is unbounded.
By Corollary~\ref{corollary:Jumps-interval-bound} the sequence $(T^{(d)}_j)_j$ is $(d_0,(x_t)_t)$-good in the first three cases of Theorem~\ref{theorem:master},
and by Lemma~\ref{lemma:fourthcase} it is $(d_0,(x_t)_t)$-left-good in the fourth case of Theorem~\ref{theorem:master}.
Hence, by assumption \eqref{equation:bounded} is true,
which is a contradiction to the unboundedness of the function $h_d$.
\end{proof}

\begin{corollary}
\label{corollary:length-unbounded}
For every $\ell\in\IN$ there are infinitely many $t$ with $\mathcal{T}\restrictedto \ell \prefixeq \mathcal{T}[t]$.
\end{corollary}

\begin{proof}
For the sake of a contradiction, let us assume that the assertion is not true and that $\ell$ is the smallest natural number
such that there are only finitely many $t$ with $\mathcal{T}\restrictedto \ell \prefixeq \mathcal{T}[t]$.
Then $\ell>0$.
Let $\sigma$ be the prefix of $\mathcal{T}$ of length $\ell-1$.
Then there are infinitely many $t$ with $\sigma\prefixeq \mathcal{T}[t]$, but only
finitely many $t$ with $\sigma0\prefixeq \mathcal{T}[t]$ and only
finitely many $t$ with $\sigma1\prefixeq \mathcal{T}[t]$.
We choose some $s>\ell$ so large
such that for $t\geq s$ we have $\sigma0\not\subseteq \mathcal{T}[t]$ and
$\sigma1\not\subseteq \mathcal{T}[t]$
and such that $\sigma$ is not initialized in any stage $t\geq s$;
compare Lemma~\ref{lemma:finitely-many-initializations}.

Let us first consider the case that $|\sigma|$ is odd.
Let us consider the infinitely many $t\geq s$ with $\sigma\subseteq \mathcal{T}[t]$.
In all of these stages $t$ the strategy $\sigma$ goes from step (I) to step (II$_\mathrm{odd}$)
and then to step (III$_\mathrm{odd}$) (otherwise we would get $\sigma1\subseteq \mathcal{T}[t]$),
Thus, there are infinitely many $t$ such that $\sigma$ reaches step (III$_\mathrm{odd}$).
But then, by Lemma~\ref{lemma:odd-infinitely-oftenIV},
there are infinitely many $t$ with $\sigma 0 \subseteq \mathcal{T}[t]$. Contradiction.

Let us now consider the case that $|\sigma|$ is even.
Let $d:=d(\sigma)[s]$. Then $\sigma=\sigma_d$ and $|\sigma_d|=2\cdot\iota(d)$. 
As $\sigma$ is not initialized in any stage $t\geq s$,
we obtain $d(\sigma)[t]=d$ for all $t\geq s$.
In each of the infinitely many stages $t\geq s > \ell=|\sigma|$ with $\sigma\prefixeq\mathcal{T}[t]$
the strategy $\sigma$ reaches step (II$_\mathrm{even}$).
From there is must reach step (III$_\mathrm{even}$)
(otherwise we would get $\sigma0\subseteq \mathcal{T}[t]$).
From there it must reach step (IV)
(otherwise we would get $\sigma1\subseteq \mathcal{T}[t]$).
Thus, $\sigma$ reaches step (IV) in infinitely many stages.

If $|\sigma|$ is divisible by $2$ but not by $4$ (that means that $\iota(d)$ is odd),
then by Lemma~\ref{lemma:4mod2-0} there are infinitely many $t'$ with $\sigma 0 \prefixeq \mathcal{T}[t']$.
Contradiction.

Let us assume that $|\sigma|$ is divisible by $4$ (that means that $\iota(d)$ is even).
That $\sigma$ reaches step (IV) in infinitely many stages implies that the function $h_d$ is unbounded.
Note that by Corollary~\ref{corollary:Jumps-interval-bound}
the sequence $(T^{(d)}_j)_j$ is $d$-good in the first three cases of Theorem~\ref{theorem:master}
and that by Lemma~\ref{lemma:fourthcase}
the sequence $(T^{(d)}_j)_j$ is $d$-left-good in the fourth case of Theorem~\ref{theorem:master}.
Thus, in each case of Theorem~\ref{theorem:master}, \eqref{equation:bounded} is true,
which is a contradiction to the unboundedness of the function $h_d$.
\end{proof}

We shall need the following corollary of Corollary~\ref{corollary:even-triggered-movement}.

\begin{corollary}
\label{corollary:computablyCauchy}
Let us assume that the assumptions of the first, second, or fourth case in Theorem~\ref{theorem:master} are true.
Let $d\in\mathit{D}$ be a number such that $\iota(d)$ is odd and $J^{(d)}$ is infinite.
Then the series $\sum_{s \in \mathrm{Jumps}^{(d)}} (x_s - x_{s-1})$ is computably Cauchy.
\end{corollary}

Note that in the third case of Theorem~\ref{theorem:master} for any $d\in\mathit{D}$ with odd $\iota(d)$ the set $J^{(d)}$ is empty.

\begin{proof}
In this proof we write $T_j$ instead of $T^{(d)}_j$.
The sequence $(T_j)_j$ is infinite, increasing, computable and satisfies $T_0>0$.
The condition \eqref{equation:computable-convergence}
is true, in the first two cases of Theorem~\ref{theorem:master} by 
Corollary~\ref{corollary:Jumps-interval-bound}
and by 
Lemma~\ref{lemma:cancel-comp-conv},
in the fourth case by Lemma~\ref{lemma:fourthcase}
and by assumption.
Thus, there exists a computable function $f\colon\IN\to\IN$
such that for all $n,j_1,j_2\in\IN$ with $f(n) \leq j_1 \leq j_2$
\[ \left| \sum_{j \colon j_1 \leq j < j_2} 
    A(d,T_j,h_d(j),x_{T_j-1}) \right| \leq 2^{-n} . \]
In particular, for all $n,j$ with $j\geq f(n)$ we obtain $|A(d,T_j,h_d(j),x_{T_j-1})|\leq 2^{-n}$.
Thus, for all $n,j_1,j_2\in\IN$ with $f(n+2) \leq j_1 \leq j_2$
we obtain
\[ \left| \sum_{j \colon j_1 \leq j < j_2 \land j\not\in \mathrm{First}(d,(T_j)_j)} 
    A(d,T_j,h_d(j),x_{T_j-1}) \right| \leq 2^{-n-2} + 2^{-n-1} < 2^{-n} . \]
As we assume that the set $J^{(d)}$ is infinite, we can 
furthermore assume without loss of generality that 
$f$ is increasing and that $\range(f)\subseteq J^{(d)}$.
The function $g:\IN\to\IN$ defined by
\[ g(n) := T_{f(n+2)} \]
is computable and increasing.
We claim that it witnesses that the series $\sum_{s \in \mathrm{Jumps}^{(d)}} (x_s - x_{s-1})$ is computably Cauchy,
that is, we claim that for all $n,t_1,t_2\in\IN$ with $g(n) \leq t_1 \leq t_2$
\[ \left| \sum_{s \in \mathrm{Jumps}^{(d)} \land t_1 \leq s < t_2} (x_s - x_{s-1}) \right| \leq 2^{-n} . \]
Let us fix an $n\in\IN$ and $t_1,t_2\in\IN$ with $g(n) \leq t_1 \leq t_2$.
Let $k_1$ be the largest number in $J^{(d)}$ with $T_{k_1} \leq t_1$,
and let $k_2$ be the largest number in $J^{(d)}$ with $T_{k_2} \leq t_2$.
Then $f(n+2) \leq k_1 \leq k_2$.
By Corollary~\ref{corollary:even-triggered-movement} the sum
$\sum_{s \in \mathrm{Jumps}^{(d)} \land t_1 \leq s < t_2} (x_s - x_{s-1})$
lies between the smallest and the largest of the following four numbers:
\begin{itemize}
\item 
the sum $\sum_{j \colon k_1 <  j < k_2 \land j\not\in \mathrm{First}(d,(T_j)_j)} A(d,T_j,h_{d}(j),x_{T_j-1})$,
\item 
the sum $\sum_{j \colon k_1 \leq  j < k_2 \land j\not\in \mathrm{First}(d,(T_j)_j)} A(d,T_j,h_{d}(j),x_{T_j-1})$,
\item 
the sum $\sum_{j \colon k_1 <  j \leq k_2 \land j\not\in \mathrm{First}(d,(T_j)_j)} A(d,T_j,h_{d}(j),x_{T_j-1})$,
\item 
the sum $\sum_{j \colon k_1 \leq  j \leq k_2 \land j\not\in \mathrm{First}(d,(T_j)_j)} A(d,T_j,h_{d}(j),x_{T_j-1})$.
\end{itemize}
All four numbers lie in the interval $\left[-2^{-n},2^{-n}\right]$.
Then so does the sum $\sum_{s \in \mathrm{Jumps}^{(d)} \land t_1 \leq s < t_2} (x_s - x_{s-1})$.
That was to be shown.
\end{proof}

\begin{lemma}
\label{lemma:converging}
The sequence $(x_t)_t$ converges. 
\end{lemma}

\begin{proof}
Let us fix an arbitrary $\varepsilon >0$.
We wish to show that there is a number $t_0$ such that for all $t > t_0$
\[ \left| x_{t-1} - x_{t_0-1} \right| < \varepsilon . \]
First, we choose a natural number $d_0$ large enough such that $(d_0+2) \cdot 2^{-d_0} < \varepsilon$.
We observe that for any $s\in\IN$ with $x_s \neq x_{s-1}$
there are unique $d\in\IN$ and a unique $j\in\dom(T^{(d)})$ with $s \in \mathrm{Jumps}(T^{(d)}_j)$.
Next we consider a number $d\leq d_0$ and consider the following three cases.
\begin{itemize}
\item
$\iota(d)$ is even. 
In the proofs of Lemma~\ref{lemma:finitely-many-initializations} and of Corollary~\ref{corollary:length-unbounded} we have already seen that in this case the condition~\eqref{equation:bounded} is true.
This implies that the set $J^{(d)}$ is finite.
As for every $j\in\dom(T^{(d)})$ the set $\mathrm{Jumps}(T^{(d)}_j)$ is finite by Lemma~\ref {lemma:triggered-StepV-finite},
we can choose a number $t^{(d)}$ such that for every $s\in \mathrm{Jumps}^{(d)}$ we have $s< t^{(d)}$.
\item
$\iota(d)$ is odd and the set $J^{(d)}$ is finite.
Then as in the first case we can choose a number $t^{(d)}$ such that for every $s\in \mathrm{Jumps}^{(d)}$ we have $s< t^{(d)}$.
\item
$\iota(d)$ is odd and the set $J^{(d)}$ is infinite.
Then we are in the first, the second, or the fourth case of Theorem~\ref{theorem:master}
because in the third case of Theorem~\ref{theorem:master} the set $J^{(d)}$ is empty
for $d$ with odd $\iota(d)$.
By Corollary~\ref{corollary:computablyCauchy} there is a number $t^{(d)}$ such that for all $t_1,t_2\in\IN$ with $t^{(d)} \leq t_1 \leq t_2$
\[ \left| \sum_{s \in \mathrm{Jumps}^{(d)} \land t_1 \leq s < t_2} (x_s - x_{s-1}) \right| \leq 2^{-d_0} . \]
\end{itemize}
Finally, we define $t_0:=\max\{t^{(d)} \colon d\leq d_0\}$.
Above we observed that for any $s\in\IN$ with $x_s \neq x_{s-1}$ there are unique $d\in\IN$ and a unique $j\in\dom(T^{(d)})$ with
$s \in \mathrm{Jumps}(T^{(d)}_j)$.
We have chosen $t_0$ in such a way that for $s\geq t_0$ with $x_s \neq x_{s-1}$ the unique numbers $d\in\IN$ and $j\in\dom(T^{(d)}$ with
$s \in \mathrm{Jumps}(T^{(d)}_j)$ satisfy the following: either $d> d_0$ or $d\leq d_0$ and $\iota(d)$ is odd and the set $J^{(d)}$ is infinite.
Let us consider a number $t>t_0$. We obtain
\begin{align*}
  \lefteqn{\left| x_{t-1} - x_{t_0-1} \right|} && \\
	&= \left| \sum_{s=t_0}^{t-1} (x_s - x_{s-1}) \right| \\
	&= \left| \sum_{d>d_0} \sum_{s\in\mathrm{Jumps}^{(d)} \land t_0\leq s < t} (x_s - x_{s-1})
	   + \sum_{\substack{d\colon d< d_0 \land \iota(d)\text{ is odd } \\ \land J^{(d)} \text{ is infinite}}} \sum_{s\in\mathrm{Jumps}^{(d)} \land t_0\leq s < t} (x_s - x_{s-1}) \right| \\
	&\leq \sum_{d>d_0} 2^{-d+1} + \sum_{\substack{d\colon d< d_0 \land \iota(d)\text{ is odd } \\ \land J^{(d)} \text{ is infinite}}} 2^{-d_0} \\
	&\leq 2 \cdot 2^{-d_0} + d_0 \cdot 2^{-d_0} \\
	&< \varepsilon.
\end{align*}
In the fourth line, for the estimate of the first sum (the sum $\sum_{d>d_0} \ldots$) we used Corollary~\ref{corollary:Jumps-interval-bound}.2
and Lemma~\ref{lemma:fourthcase}.2.
\end{proof}

We define $\xi:=\lim_{t\to\infty} x_t$.

\subsection{Satisfaction of the even requirements}
\label{subsection:8}

\begin{corollary}
\label{corollary:evenrequirements}
For all $i$ the even requirement $\mathcal{R}_{2i}(\xi)$ is satisfied.
\end{corollary}

\begin{proof}
Let us consider some $i\in\IN$.
The sequence $(x_t)_t$ is computable and, according to Lemma~\ref{lemma:converging}, converging.
Let $\tau$ be the prefix of length $2i$ of $\mathcal{T}$.
By Lemma~\ref{lemma:finitely-many-initializations} $\tau$ is initialized only finitely often.
We define a number $s\in\IN$ as follows. If $\tau$ is never initialized then we set $s:=0$. Otherwise we let $s$ be the last stage in which $\tau$ is initialized.
Then $d(\tau)[t]=d(\tau)[s]=:d$ for all $t\geq s$. 
By Lemma~\ref{corollary:length-unbounded}
there are infinitely many $t$ with $\tau \prefixeq \mathcal{T}[t]$.
The sequence $(T_j^{(d)})_j$ contains all such stages $t$ with $t>s$.
Hence, the sequence $(T_j^{(d)})_j$ is infinite.
It is clearly a computable sequence.
In the first three cases of Theorem~\ref{theorem:master} it is a $(d,(x_t)_t)$-good sequence by Corollary~\ref{corollary:Jumps-interval-bound},
in the fourth case of Theorem~\ref{theorem:master} it is a $(d,(x_t)_t)$-left-good sequence by Lemma~\ref{lemma:fourthcase}.
In each case, by assumption the even requirement $\mathcal{R}_{2i}(\xi)$ is satisfied.
\end{proof}

\subsection{The sign switching and weakly computable case}
\label{subsection:9}

Next we wish to show that the number $\xi$ is weakly computable if the assumptions in the sign switching and weakly computable case in Theorem~\ref{theorem:master} are satisfied.

\begin{lemma}
\label{lemma:absolute-convergence-bound}
Let us assume that the assumptions in the sign switching and weakly computable case in Theorem~\ref{theorem:master} are satisfied.
Then for every $d\in \mathit{D}$ 
\[ \sum_{s \in\mathrm{Jumps}^{(d)}} |x_s - x_{s-1}|  \leq 2^{-d}. \]
\end{lemma}

\begin{proof}
Let us consider a $d\in \mathit{D}$ and a number $t\in\IN$.
It is sufficient to show $\sum_{s \in\mathrm{Jumps}^{(d)} \restrictedto t} |x_s - x_{s-1}|  \leq 2^{-d}$.
This is clear if $\mathrm{Jumps}^{(d)} \restrictedto t = \emptyset$.
Let us assume that $\mathrm{Jumps}^{(d)} \restrictedto t \neq \emptyset$.
Then $\dom(T^{(d)})\neq\emptyset$.
Let $J$ be the largest number with $J\in\dom(T^{(d)})$ and $T_J^{(d)} < t$.
By Corollary~\ref{corollary:even-triggered-movement} the sum
\[ \sum_{s \in\mathrm{Jumps}^{(d)}\restrictedto t} |x_s - x_{s-1}|
  = \sum_{j=0}^{J-1} \sum_{s \in\mathrm{Jumps}(T_j^{(d)})\restrictedto T_{j+1}^{(d)}} |x_s - x_{s-1}|
	   + \sum_{s \in\mathrm{Jumps}(T_J^{(d)})\restrictedto t} |x_s - x_{s-1}|
\]
is bounded from above by the sum 
\[ \sum_{j\in\{0,\ldots,J\} \setminus \mathrm{First}(d,(T^{(d)}_j)_j)} \left|A(d,T_j^{(d)},h_{d}(j),x_{T_j^{(d,0)}-1})\right| 
 \leq \sum_{j=0}^{J} \left|A(d,T_j^{(d)},h_{d}(j),x_{T_j^{(d,0)}-1})\right| . \]
By the first assertion of Lemma~\ref{lemma:even-movement} the finite sequence $(T_j^{(d)})_j|^{<t}$ is $(d,(x_t)_t)$-good.
By the assumption \eqref{equation:bound-on-absolute-sum} the previous sum is bounded by $2^{-d}$.
Then this applies to the sum $\sum_{s \in\mathrm{Jumps}^{(d)} \restrictedto t} |x_s - x_{s-1}|$ as well.
\end{proof}

\begin{corollary}
If the assumptions in the sign switching and weakly computable case in Theorem~\ref{theorem:master} are satisfied,
then the series $\sum_t |x_{t+1} - x_t|$ converges, hence, then the number $\xi$ is weakly computable.
\end{corollary}

\begin{proof}
Let us assume that the assumptions in the sign switching and weakly computable case in
Theorem~\ref{theorem:master} are satisfied.
By a result by Ambos-Spies, Weihrauch and Zheng \cite{AWZ2000}
it is sufficient to show that the series $\sum_t |x_{t+1} - x_t|$ converges; compare Definition~\ref{definition:wc} and Proposition~\ref{proposition:wc}.

Let us consider an arbitrary $\varepsilon >0$.
Let $d'\in\IN$ be so large that $2^{-d'} < \varepsilon/3$.
By Lemma~\ref{lemma:absolute-convergence-bound} 
for every $d\in\mathit{D}$ 
\[ \sum_{s \in\mathrm{Jumps}^{(d)}} |x_s - x_{s-1}|  \leq 2^{-d}. \]
Let us choose a number $t_0>0$ so large that for every $d\in D$ with $d<d'$
\[ \sum_{s \in\mathrm{Jumps}^{(d)} \wedge s\geq t_0} |x_s - x_{s-1}|  \leq 2^{-d'}/d' . \]
We obtain
\begin{eqnarray*}
  \lefteqn{\sum_{s \geq t_0} |x_s - x_{s-1}|} && \\
	&=&  \sum_{d\in \mathit{D}} \sum_{s \in\mathrm{Jumps}^{(d)} \wedge s \geq t_0} |x_s - x_{s-1}| \\
	&=& \sum_{d\in \mathit{D} \land d\geq d'} \sum_{s \in\mathrm{Jumps}^{(d)} \wedge s \geq t_0} |x_s - x_{s-1}| 
	  + \sum_{d\in \mathit{D} \land d<d'} \sum_{s \in\mathrm{Jumps}^{(d)} \wedge s \geq t_0} |x_s - x_{s-1}| \\
	&\leq & \sum_{d\in \mathit{D} \land d\geq d'} 2^{-d} + \sum_{d\in \mathit{D} \land d<d'} 2^{-d'}/d' \\
	&\leq & 2^{-d'+1} + 2^{-d'} 
	= 3 \cdot 2^{-d'} < \varepsilon . 	
\end{eqnarray*}
Thus, the series $\sum_t |x_{t+1} - x_t|$ converges.
\end{proof}

\subsection{Satisfaction of the odd requirements}
\label{subsection:10}

Finally we turn to strategies $\sigma$ of odd length.

\begin{lemma}
\label{lemma:triggered-3}
Let a strategy $\sigma$ of odd length reach step (III$_\mathrm{odd}$) in a stage $t_1$. 
Let $t_2$ be a stage with $t_1<t_2$ and $t_2\in\mathrm{StepIV}$
and such that there is no stage $t$ with $t_1<t\leq t_2$ in which $\sigma$ is initialized or reaches step  (III$_\mathrm{odd}$).
\begin{enumerate}
\item
If there is a $j\geq 1$ with $\mathrm{trigIV}^{\circ j}(t_2)=t_1$ then 
$\sigma$ reaches step (IV) in stage $t_1$ and $\mathcal{T}[t_1]=\sigma$.
\item
If there is no $j\geq 1$ with $\mathrm{trigIV}^{\circ j}(t_2)=t_1$ then $\mathrm{trigult}(t_2)>t_1$
and $\mathcal{T}[\mathrm{trigult}(t_2)]$ is either initialized in stage $t_1$ or has the property $\mathcal{T}[\mathrm{trigult}(t_2)]0 \prefixeq \sigma$.
\end{enumerate}
\end{lemma}

\begin{proof}
Let $\gamma:=\mathcal{T}[t_2]$.

\emph{Claim 1: Either $\gamma0\prefixeq\sigma$ or $\sigma0 <_L \gamma$.}

We cannot have $\gamma0 <_L \sigma$ because then $\sigma$ would be initialized in stage $t_2$.
We cannot have $\gamma=\sigma$ because then in stage $t_2$ the strategy $\sigma$ would have to reach step (IV) and before that step (III$_\mathrm{odd}$).
We cannot have $\sigma0\prefixeq\gamma$ because then $\sigma$ would have to reach step
(III$_\mathrm{odd}$) in stage $t_2$ as well.
Therefore we must have either $\gamma0\prefixeq\sigma$ or $\sigma0 <_L \gamma$.
We have proved Claim 1.

\emph{Claim 2: If $|\gamma|$ is even then either $\gamma$ is initialized in stage $t_1$ or $\gamma0\prefixeq\sigma$.}

Let us assume that $|\gamma|$ is even.
By Claim 1 either $\gamma0\prefixeq\sigma$ or $\sigma0 <_L \gamma$.
In the case $\gamma0\prefixeq\sigma$ the strategy $\gamma$ cannot be initialized in stage $t_1$ because then also $\sigma$ would be initialized in stage $t_1$,
which is not the case because $\sigma\prefixeq \mathcal{T}[t_1]$. 
Let us now consider the case $\sigma0 <_L \gamma$. 
If in stage $t_1$ the strategy $\sigma$ reaches step (IV) then $\gamma$ is initialized in stage $t_1$.
Remember that in stage $t_1$ the strategy $\sigma$ reaches step (III$_\mathrm{odd}$).
If in stage $t_1$ the strategy $\sigma$ does not reach step (IV) then $\sigma0\prefixeq\mathcal{T}[t_1]$,
hence $\mathcal{T}[t_1] <_L \gamma$ and $\gamma$ is initialized in stage $t_1$ as well.
So, in the case $\sigma0 <_L \gamma$ the strategy $\gamma$ is initialized in stage $t_1$.
We have proved Claim 2.

\emph{Claim 3: If $|\gamma|$ is odd then either $\mathrm{trigIV}(t_2)=t_1$ and 
$\sigma$ reaches step (IV) in stage $t_1$
or $t_1<\mathrm{trigIV}(t_2)$.}

Let us assume that $|\gamma|$ is odd. 
By Claim 1 we know $\gamma0\prefixeq\sigma$ or $\sigma0 <_L \gamma$.
Let us first consider the case $\sigma0 <_L \gamma$.
By the same argument as in the proof of Claim 2, $\gamma$ is initialized in stage $t_1$.
This leads to $q(\gamma)[t_1]=0$, and by Lemma~\ref{lemma:trig-1} we obtain
$t_1 < \mathrm{trigIV}(t_2)$.
Now, let us consider the case $\gamma0\prefixeq\sigma$.
Then $\gamma0\prefixeq\sigma\prefixeq\mathcal{T}[t_1]$,
hence, $q(\gamma)[t_1-1]=0$. 
By Lemma~\ref{lemma:trig-1} we obtain $t_1 \leq \mathrm{trigIV}(t_2)$,
hence, either $t_1<\mathrm{trigIV}(t_2)$ or $t_1=\mathrm{trigIV}(t_2)$.
Let us assume $t_1=\mathrm{trigIV}(t_2)$.
Then $\gamma$ is the longest string of odd length with $\gamma0\prefixeq\mathcal{T}[t_1]$.
If $\sigma$ would not reach step (IV) in stage $t_1$ then we would obtain $\sigma0\prefixeq\mathcal{T}[t_1]$
which, in combination with $\gamma0\prefixeq\sigma$ would lead to a contradiction.
We conclude that $\sigma$ reaches step (IV) in stage $t_1$.
We have proved Claim 3.

Remember that we assume $t_1 < t_2$.
If $\gamma$ has even length then $t_1 < t_2 = \mathrm{trigult}(t_2)$
and by Claim 2 either $\gamma=\mathcal{T}[\mathrm{trigult}(t_2)]$ is initialized in stage $t_1$ or $\gamma0\prefixeq\sigma$.
Let us assume that $\gamma$ has odd length.
By Claim 3 either there is a $j\geq 1$ with $\mathrm{trigIV}^{\circ j}(t_2)=t_1$,
and in this case $\sigma$ reaches step (IV) in stage $t_1$,
or for all $j\geq 0$ such that $\mathrm{trigIV}^{\circ j}(t_2)$ is defined we have
$t_1 < \mathrm{trigIV}^{\circ j}(t_2)$.
In this case also $t_1 < \mathrm{trigult}(t_2)$.
Then $\delta:=\mathcal{T}[\mathrm{trigult}(t_2)]$ has even length and reaches step (IV) in stage $\mathrm{trigult}(t_2)$.
By Claim 2, applied to $\delta$ instead of $\gamma$,
either the strategy $\delta$ is initialized in stage $t_1$ or $\delta0\prefixeq\sigma$.
This ends the proof of both assertions.
\end{proof}

\begin{lemma}
\label{lemma:main-odd}
Let $t_1,t_2$ be two numbers with $0<t_1<t_2$ and let $\sigma$ be a strategy of odd length 
such that $\sigma$ reaches step (III$_\mathrm{odd}$) in stage $t_1$ 
and such that there is no stage $t$ with $t_1<t < t_2$ in which $\sigma$ is initialized or reaches step  (III$_\mathrm{odd}$).
Then,
\[ \left|x_{t_2-1} - x_{t_1-1}\right| \leq 2\cdot 2^{-r(\sigma)[t_1-1]}
   + \sum_{e\in D^\mathrm{cancel}(t_1,t_2)} \left|\sum_{s\in \mathrm{Jumps}^{(e)} \land t_1 \leq s < t_2} (x_s - x_{s-1}) \right| , \]
where the set $D^\mathrm{cancel}(t_1,t_2)$ is defined as in the proof of  Lemma~\ref{lemma:even-movement}.
\end{lemma}

\begin{proof}
We start with several claims.

\smallskip

\emph{Claim 1: If $x_{t_1}\neq x_{t_1-1}$ then $\sigma$ reaches step (IV) in stage $t_1$,
hence, $t_1 \in\mathrm{Jumps}(t_1)\restrictedto t_2$.} 

Let us assume that $x_{t_1}\neq x_{t_1-1}$.
Then $\mathcal{T}[t_1]$ reaches step (IV) in stage $t_1$
and there is no strategy $\gamma$ of odd length
with $\gamma0\prefixeq\mathcal{T}[t_1]$.
If $\sigma$ would not reach step (IV) in stage $t_1$ then we would get $\sigma0\prefixeq\mathcal{T}[t_1]$, a contradiction.
We have proved Claim 1.

\smallskip

\emph{Claim 2: If $\sigma$ does not reach step (IV) in stage $t_1$ then
$\mathrm{Jumps}(t_1) \restrictedto t_2 = \emptyset$, hence,
\[ \left|\sum_{s\in\mathrm{Jumps}(t_1)\restrictedto t_2} (x_s - x_ {s-1})\right| = 0 . \]
}

Let us assume that $\mathrm{Jumps}(t_1)\restrictedto t_2 \neq \emptyset$.
We wish to show that then $\sigma$ reaches step (IV) in stage $t_1$.
There exists a number $s\in\mathrm{Jumps}(t_1)\restrictedto t_2$.
If $s=t_1$ then $x_{t_1}\neq x_{t_1-1}$,
and by Claim 1 this implies that $\sigma$ reaches step (IV) in stage $t_1$.
If $s \neq t_1$ then $t_1 < s < t_2$ by Lemma~\ref{lemma:Jumps}.2.
Note  that there is no stage $t$ with $t_1<t \leq s$ in which $\sigma$ is initialized or reaches step  (III$_\mathrm{odd}$) because $s<t_2$.
By the first assertion of Lemma~\ref{lemma:triggered-3} $\sigma$ reaches step (IV) in stage $t_1$.
We have proved Claim 2.

\smallskip

\emph{Claim 3: If $\sigma$ reaches step (IV) in stage $t_1$ then
\[ \left|\sum_{s\in\mathrm{Jumps}(t_1)\restrictedto t_2} (x_s - x_ {s-1})\right| \leq 2^{-r(\sigma)[t_1-1]}. \]
}

Let us assume that $\sigma$ reaches step (IV) in stage $t_1$.
By the first and the third assertion of Lemma~\ref{lemma:triggered-movement} we obtain
\[  \left|\sum_{s \in\mathrm{Jumps}(t_1)\restrictedto t_2} (x_s - x_{s-1}) \right| 
  \leq \left|v_{t_1}\right|
	\leq 2^{-r(\sigma)[t_1-1]}. \]
We have proved Claim 3.

\smallskip

We see
\[    x_{t_2-1} - x_{t_1-1} = \sum_{s=t_1}^{t_2-1} (x_s - x_{s-1}) . \]
Let us consider a number $s$ with $t_1\leq s < t_2$ and $x_s\neq x_{s-1}$.
Then $\mathcal{T}[s]$ reaches step (IV) in stage $s$.
If $s\not\in\mathrm{Jumps}(t_1)\restrictedto t_2$ then 
by Claim 1 $s\neq t_1$, hence, $t_1 < s < t_2$, and by the second assertion of 
Lemma~\ref{lemma:triggered-3}, $t_1 < \mathrm{trigult}(s)$
and either $\mathcal{T}[\mathrm{trigult}(s)]$ is initialized in stage $t_1$
or $\mathcal{T}[\mathrm{trigult}(s)]0\prefixeq\sigma$.
We obtain $s\in\mathrm{Jumps}^{(e)}\restrictedto t$ where 
\[ e:=d(\mathcal{T}[\mathrm{trigult}(s)])[\mathrm{trigult}(s)-1] . \]
Then $e\in D^\mathrm{init}(t_1,t_2) \cup D^\mathrm{cancel}(t_1,t_2)$
where the sets $D^\mathrm{init}(t_1,t_2)$ and $D^\mathrm{cancel}(t_1,t_2)$ 
are defined as in the proof of  Lemma~\ref{lemma:even-movement}.
We obtain
\begin{eqnarray*}
   \lefteqn{\left| x_{t_2-1} - x_{t_1-1} \right|} && \\
	 &=& \left| \sum_{s=t_1}^{t_2-1} (x_s - x_{s-1}) \right| \\
	 &=&\left| \sum_{s\in\mathrm{Jumps}(t_1)\restrictedto t_2} (x_s - x_{s-1})
	    + \sum_{e\in D^\mathrm{init}(t_1,t_2) \cup D^\mathrm{cancel}(t_1,t_2)} \sum_{s\in \mathrm{Jumps}^{(e)}\land t_1 \leq s < t_2} (x_s - x_{s-1}) \right| \\
	 &\leq & \left| \sum_{s\in\mathrm{Jumps}(t_1)\restrictedto t_2} (x_s - x_{s-1}) \right|
	    + \sum_{e\in D^\mathrm{init}(t_1,t_2)} \left|\sum_{s\in \mathrm{Jumps}^{(e)}\land t_1 \leq s < t_2} (x_s - x_{s-1}) \right| \\
	&&		+ \sum_{e\in D^\mathrm{cancel}(t_1,t_2)} \left|\sum_{s\in \mathrm{Jumps}^{(e)}\land t_1 \leq s < t_2} (x_s - x_{s-1}) \right| \\
	&\leq &	2^{-r(\sigma)[t_1-1]} + \sum_{e\in D^\mathrm{init}(t_1,t_2)} 2^{-e} 
	    + \sum_{e\in D^\mathrm{cancel}(t_1,t_2)} \left|\sum_{s\in \mathrm{Jumps}^{(e)}\land t_1 \leq s < t_2} (x_s - x_{s-1}) \right| \\
	&\leq & 2^{-r(\sigma)[t_1-1]} + 2^{-r(\sigma)[t_1-1]} 
	  + \sum_{e\in D^\mathrm{cancel}(t_1,t_2)} \left|\sum_{s\in \mathrm{Jumps}^{(e)}\land t_1 \leq s < t_2} (x_s - x_{s-1}) \right| ,
\end{eqnarray*}
where in the second last line we applied Claims 2 and 3 and Lemma~\ref{lemma:even-movement}.2
as well as Lemma~\ref{lemma:fourthcase}.2,
and in the last line we used that for $e\in D^\mathrm{init}(t_1,t)$ we have $e>1+r(\sigma)[t_1-1]$,
because $|\sigma|$ is odd and either $\sigma0\prefixeq\mathcal{T}[t_1]$ or $\sigma=\mathcal{T}[t_1]$.
That was to be shown.
\end{proof}

\begin{lemma}
For all $i$ the odd requirement $\mathcal{R}_{2i+1}$ is satisfied.
\end{lemma}

\begin{proof}
Let us consider some $i\in\IN$.
We assume that $\varphi_i$ is a total and increasing function because
otherwise the requirement $\mathcal{R}_{2i+1}$ is satisfied trivially.
We wish to show that the sequence
$(x_{\varphi_i(n+1)} - x_{\varphi_i(n)})_n $ converges computably to zero.

Let $\sigma$ be the prefix of length $2i+1$ of $\mathcal{T}$.
Remember that by Corollary~\ref{corollary:length-unbounded}
there are infinitely many $t>|\sigma|$ with $\sigma\prefixeq\mathcal{T}[t]$.
In all of these stages $t$ the strategy $\sigma$ reaches step (II$_\mathrm{odd}$).
We claim that in infinitely many of these stages $t$ the strategy $\sigma$ reaches even step (III$_\mathrm{odd}$).
For the sake of a contradiction, let us assume that this is not true.
But then, on the one hand, the number 
\[ \max \left( \{0\} \cup \{s<t \mid \sigma \text{ reaches step (III$_\mathrm{odd}$) in stage } s \} \right) \]
on the right hand side of \eqref{equation:III-1} will remain constant from a certain stage $t$  on.
On the other hand, for sufficiently large $t$ we have $\ell_i(t)>0$, and for increasing $t$ the number 
$ \varphi_i(\ell_i(t)-1)$ on the left hand side of \eqref{equation:III-1} tends to infinity
because $\varphi_i$ is total and increasing.
Hence, then for all sufficiently large $t$ we would have $\ell_i(t)>0$ and \eqref{equation:III-1} would be true.
But these are just the conditions stated in step (II$_\mathrm{odd}$)
that must be satisfied so that $\sigma$ would move from step (II$_\mathrm{odd}$) to step (III$_\mathrm{odd}$).
Hence, $\sigma$ will reach step (III$_\mathrm{odd}$) again, contradiction.
We have shown that $\sigma$ reaches step (III$_\mathrm{odd}$) in infinitely many stages $t$.

Let $t_0$ be so large such that
\begin{itemize}
\item
$t_0>|\sigma|$,
\item
$\ell_i(t_0)>0$,
\item
no $\rho$ with $\rho<_L\sigma$ is applied in any stage $t\geq t_0$
(remember Lemma~\ref{lemma:infinite-injury-applied}),
\item
$\sigma$ is not initialized in any stage $t\geq t_0$,
(remember Lemma~\ref{lemma:finitely-many-initializations}); note that
then no $\rho$ with $\rho\prefixeq \sigma$ is initialized in any stage $t\geq t_0$,
\end{itemize}
Then for every $\rho$ of even length with $\rho0\subseteq\sigma$ 
we have $d(\rho)[t] = \lim_{s\to\infty} d(\rho)[s]$, for all $t\geq t_0$.
Let us write $d(\rho)[\infty]$ for $\lim_{s\to\infty} d(\rho)[s]$.
Furthermore, the length of $\rho$ is not divisible by four, by Lemma~\ref{lemma:div4-no-0}.
Then $J^{(d(\rho)[\infty])}$ is either finite or infinite.
In the first case, by Lemma~\ref{lemma:triggered-StepV-finite} the set $\mathrm{Jumps}^{(d(\rho)[\infty])}$ is finite as well.
We can additionally choose $t_0$ so large such that
\begin{itemize}
\item
for every $\rho$ of even length with $\rho0\subseteq\sigma$ 
and such that the set $J^{(d(\rho)[\infty])}$ is finite, $t_0> \max \mathrm{Jumps}^{(d(\rho)[\infty])}$.
\end{itemize}
We set
\[ D^\mathrm{infinite} := \{e \in D \colon (\exists \rho) \, (\rho0\prefixeq \sigma \text{ and } e = d(\rho)[\infty] 
    \text{ and } J^{(e)} \text{ is infinite} ) \} .\]
Clearly, this set is finite: it has less than $|\sigma|$ elements.
By Lemma~\ref{lemma:main-odd}
we conclude that for any two numbers $t_1,t_2$ with $t_0\leq t_1 < t_2$
such that $\sigma$ reaches step (III$_\mathrm{odd}$) in stage $t_1$ 
and such that there is no stage $t$ with $t_1<t < t_2$ in which $\sigma$ reaches step  (III$_\mathrm{odd}$)
\begin{equation}
\label{equation:odd-infinite}
  \left|x_{t_2-1} - x_{t_1-1}\right| \leq 2\cdot 2^{-r(\sigma)[t_1-1]}
   + \sum_{e\in D^\mathrm{infinite}} \left|\sum_{s\in \mathrm{Jumps}^{(e)}\land t_1 \leq s < t_2} (x_s - x_{s-1}) \right| .
\end{equation}
By Corollary~\ref{corollary:computablyCauchy}
for any $e\in D^\mathrm{infinite}$
the series $\sum_{s \in \mathrm{Jumps}^{(e)}} (x_s - x_{s-1})$ is computably Cauchy.
That means that there exists an increasing, computable function $g^{(e)}:\IN\to\IN$ such that for all $n,t_1,t_2\in\IN$ with $g^{(e)}(n) \leq t_1 \leq t_2$
\[ \left| \sum_{s \in \mathrm{Jumps}^{(e)} \land t_1 \leq s < t_2} (x_s - x_{s-1}) \right| \leq 2^{-n} . \]

Above we have seen that there are infinitely many stages $t\geq t_0$ in which $\sigma$ reaches step (III$_\mathrm{odd}$).
By Lemma~\ref{lemma:odd-infinitely-oftenIV}
there are infinitely many $t\geq t_0$ with  $\sigma0\prefixeq \mathcal{T}[t]$,
hence, $r(\sigma)$ increases infinitely often.
Let $s_0$ be the smallest number with $s_0>t_0$ and such that
$\sigma$ reaches step (III$_\mathrm{odd}$) in stage $s_0$, and, for $m\geq 0$,
let $s_{m+1}$ be the smallest number with $s_{m+1}>s_m$ and such that
$\sigma$ reaches step (III$_\mathrm{odd}$) in stage $s_{m+1}$.
For all $m$, $\ell_i(s_m)>0$ and
\begin{equation}
\label{equation:s-m-1}
\varphi_i(\ell_i(s_{m+1})-1) \geq s_m 
\end{equation}
by \eqref{equation:III-1},
and $\varphi_i(\ell_i(s_{m+1})-1) < s_{m+1}$ by definition of $\ell_i$.
We define the function $h:\IN\to\IN$ by
\[ h(n) := \min\{m \in\IN ~: \begin{array}[t]{l}
    r(\sigma)[s_m-1] \geq n+4 \text{ and } \\
		(\forall e \in D^\mathrm{infinite}) \, s_m \geq g^{(e)}(n + 3 + \lceil \log_2(|\sigma|)\rceil) \}.
		\end{array} \]
for $n\in\IN$. It is computable, nondecreasing, and tends to infinity, and $m\geq h(n)$ implies $r(\sigma)[s_m-1] \geq n+4$
(because the function $m \mapsto r(\sigma)[s_m-1]$ is nondecreasing).
We define the function $H:\IN\to\IN$ by
\[ H(n) := \ell_i(s_{h(n)+1})-1 , \]
for $n\in\IN$. It is computable, nondecreasing, and tends to infinity.
We claim that the sequence  $(x_{\varphi_i(k+1)} - x_{\varphi_i(k)})_k$ converges to zero and 
that the function $H$ is a modulus of convergence of the
sequence $(x_{\varphi_i(k+1)} - x_{\varphi_i(k)})_k$.
Indeed, for $k,n\in\IN$ with $k\geq H(n)$ let $m$ be the largest number with $k\geq \ell_i(s_{m+1})-1$.
Then $m\geq h(n)$, hence, $r(\sigma)[s_m-1] \geq n+4$. 
Using  \eqref{equation:s-m-1} we obtain
\[ \varphi_i(k) \geq \varphi_i(\ell_i(s_{m+1})-1) \geq s_m . \]
As $k+1 \leq \ell_i(s_{m+2})-1$, we obtain
\[ \varphi_i(k+1) \leq \varphi_i(\ell_i(s_{m+2})-1) < s_{m+2} . \]
Thus, we have
\[ s_m \leq \varphi_i(k) < \varphi_i(k+1) < s_{m+2} . \]
If $\varphi_i(k) < s_{m+1}$ then by \eqref{equation:odd-infinite} we obtain
\begin{eqnarray*}
   \left| x_{\varphi_i(k)} - x_{s_m-1} \right|
   &\leq & 2 \cdot 2^{-r(\sigma)[s_m-1]} 
	 + \sum_{e\in D^\mathrm{infinite}} \left|\sum_{s\in \mathrm{Jumps}^{(e)}\land s_m \leq s < \varphi_i(k)} (x_s - x_{s-1}) \right| \\
   &\leq & 2 \cdot 2^{-n-4} +  \sum_{e\in D^\mathrm{infinite}} 2^{-n - 3 - \lceil \log_2(|\sigma|)\rceil} \\
	 &\leq &2^{-n-2} .
\end{eqnarray*}
Otherwise $s_{m+1} \leq \varphi_i(k) < s_{m+2}$, and in the same way we obtain
\[    \left| x_{\varphi_i(k)} - x_{s_{m+1}-1} \right| \leq 2^{-n-2} . \]
Furthermore, by the same argument we also obtain
\[    \left| x_{s_{m+1}-1} - x_{s_m-1} \right| \leq 2^{-n-2} . \]
This implies
\[    \left| x_{\varphi_i(k)} - x_{s_m-1} \right| \leq 2^{-n-1} , \]
which is true in the first case as well.
By applying the same argument to $k+1$ instead of to $k$ we obtain the analogous estimate
\[ \left| x_{\varphi_i(k+1)} -  x_{s_m-1} \right| \leq 2^{-n-1} . \]
Putting all of this together we obain
\begin{eqnarray*}
\lefteqn{\left| x_{\varphi_i(k+1)} - x_{\varphi_i(k)} \right|} && \\
&=& \left| x_{\varphi_i(k+1)} - x_{s_m-1} + x_{s_m-1} - x_{\varphi_i(k)} \right| \\
&\leq & \left| x_{\varphi_i(k+1)} - x_{s_m-1} \right| + \left| x_{\varphi_i(k)} - x_{s_m-1} \right| \\
& \leq & 2^{-n-1} + 2^{-n-1} \leq 2^{-n}.
\end{eqnarray*}
We have shown that the sequence
$(x_{\varphi_i(n+1)} - x_{\varphi_i(n)})_n$ converges computably to zero.
\end{proof}

\section{A Nearly Computable and Left-computable Number that is not Computable}
\label{section:not-right-computable}

In this section we illustrate how Theorem~\ref{theorem:master} can be applied by using it in order to derive the
following assertion, which follows from the main result by Downey and LaForte~\cite[Theorem 8]{DL2002}
and an observation by Stephan and Wu~\cite[Theorem 5]{SW2005}; see~\cite[Theorem 7.3]{HJ2024}.
Later we shall prove other results that imply this one.

\begin{corollary}
\label{corollary:nearly-left-not-right}
There exists a nearly computable and left-computable number that is not computable.
\end{corollary}

We wish to apply the nonnegative case, the third case, of Theorem~\ref{theorem:master}.
Instead, one could also use the nonnegative and left-good case, the fourth case, of Theorem~\ref{theorem:master}.
A real number $\xi$ is called \emph{right-computable} if there exists a computable, nonincreasing
sequence $(q_n)_n$ of dyadic numbers converging to $\xi$.
Clearly, a real number is computable iff it is left-computable and right-computable.
Let $(\varphi_i^\mathrm{rc})_i$ be an effective list of all right-computable
real numbers. By that we mean, $\varphi^\mathrm{rc}\colon\IN^2\to\ID$ is a computable function
such that
\begin{itemize}
\item
for every $i\in\IN$ the sequence $(\varphi^\mathrm{rc}_i(n))_n$ is nonincreasing and converging,
\item
for every right-computable number $\xi$ there exists an $i\in\IN$ such that
$\xi= \lim_{n\to\infty} \varphi^\mathrm{rc}_i(n)$.
\end{itemize}
We consider the following requirements:
we let $\mathcal{R}_{4i+2}(\xi)$ be true for every real number $\xi$
and define $\mathcal{R}_{4i}(\xi)$ as follows:
\begin{equation}
\label{eq:requirement-not-right-computable}
\xi \neq \lim_{n\to\infty} \varphi^\mathrm{rc}_i(n)
\tag{$\mathcal{R}_{4i}$}
\end{equation}
Obviously, if a real number $\xi$ satisfies the requirements $\mathcal{R}_{4i}(\xi)$ for all $i\in\IN$ then $\xi$ is not right-computable.
The strategy for constructing a left-computable number that is not right-computable is very simple:
we construct an nondecreasing sequence $(x_t)_t$ of dyadic numbers and, whenever it is close to but still smaller than a value
$\varphi^\mathrm{rc}_i(n)$, then we try to increase a bit so that it ``jumps over'' $\varphi^\mathrm{rc}_i(n)$. Once that is done,
the requirement $\mathcal{R}_{4i}(\lim_{t\to\infty} x_t$ is satisfied forever. This strategy is realized by the following funcion $A$.
We define a function $A:\IN^3\times\ID\to\ID$ by
\[ A(d,t,\ell,q) := \begin{cases}
    2^{-d} & \text{if } \iota(d) \text{ is even and }
		\varphi^\mathrm{rc}_{\iota(d)/2}(t) - 2^{-d-1} < q \leq \varphi^\mathrm{rc}_{\iota(d)/2}(t) , \\
		0 & \text{otherwise}.
		\end{cases}
\]
Note that in this case the parameter $\ell$ is not used.
It is clear that $A$ is computable and satisfies \eqref{equation:range-powers-of-2} and \eqref{eq:nonnegative-case}.
So, we are in the nonnegative case 
of Theorem~\ref{theorem:master}.

Let $d\in\IN$, let $(x_t)_t$ be a computable, nondecreasing sequence of dyadic numbers,
and let $(T_j)_j$ be a computable, nonempty, finite or infinite, increasing sequence with $T_0>0$.
If $\iota(d)$ is odd, then $A(d,t,\ell,q)=0$ for all $t,\ell\in\IN$ and $q\in\ID$.
Let us consider the case when $\iota(d)$ is even.
Let us assume additionally that the sequence $(T_j)_j$ is $(d,(x_t)_t)$-good.
We claim that \eqref{equation:bound-on-sum} and \eqref{equation:bounded} are true.
We set $i:=\iota(d)/2$.
We claim that there is at most one $j$ such that $A(d,T_j,h_{d}(j),x_{T_j-1})\neq 0$. 
Once this is shown, \eqref{equation:bound-on-sum} and \eqref{equation:bounded} follow.
Let us assume that there is some $j\in\dom(T)$ such that $A(d,T_j,h_{d}(j),x_{T_j-1})\neq 0$.
Let us consider the smallest such $j$. Then $A(d,T_j,h_{d}(j),x_{T_j-1})=2^{-d}$
and $\varphi^\mathrm{rc}_i(T_j) - 2^{-d-1} < x_{T_j-1} \leq \varphi^\mathrm{rc}_i(T_j)$.
If $\dom(T)=\{0,\ldots,j\}$ then the assertion is clear.
Otherwise, by using \eqref{equation:i-d-good} with $j':=j+1$ we obtain
$x_{T_{j+1}-1} \geq x_{T_j-1} + 2^{-d-1}$.
As $(x_t)_t$ is nondecreasing, we obtain for all $k\in\dom(T)$ with $k\geq j+1$
\[ x_{T_k-1} \geq x_{T_{j+1}-1} \geq x_{T_j-1} + 2^{-d-1} >
   \varphi^\mathrm{rc}_i(T_j) \geq \varphi^\mathrm{rc}_i(T_{j+1}) \geq \varphi^\mathrm{rc}_i(T_k) . \]
This shows that for all $k\in\dom(T)$ with $k\geq j+1$, $A(d,T_k,h_{d}(j),x_{T_k-1})=0$.
We have shown that there is at most one $j$ such that $A(d,T_j,h_{d}(j),x_{T_j-1})\neq 0$. 
Hence, \eqref{equation:bound-on-sum} and \eqref{equation:bounded} follow.

Remember that $\mathcal{R}_{4i+2}(\xi)$ is true for every real number $\xi$.
In order to show that the remaining assumptions of Theorem~\ref{theorem:master} are
satisfied in this case, it is sufficient to consider an even number $2i\in\IN$.
Let us assume that $(x_t)_t$ is a nondecreasing, converging, and computable
sequence, that $d\in\IN$ is a number such that $\iota(d)=2i$,
and that $(T_j)_j$ is a computable, infinite, increasing, $(d,(x_t)_t)$-good sequence with $T_0>0$.
We set $\xi := \lim_{t\to\infty} x_t$.
We claim that $\mathcal{R}_{4i}(\xi)$ is true.
For the sake of a contradiction, let us assume that
$\xi = \lim_{n\to\infty} \varphi^\mathrm{rc}_i(n)$.
As the sequence $(x_t)_t$ is nondecreasing and the sequence $(\varphi^\mathrm{rc}_i(n))_n$ is nonincreasing,
there must exist some $j\in\IN$ with
$\varphi^\mathrm{rc}_i(T_j) - 2^{-d-1} < x_{T_j-1} \leq \varphi^\mathrm{rc}_i(T_j)$.
But this implies $A(d,T_j,h_{d}(j),x_{T_j-1})=2^{-d}$.
As above, by using \eqref{equation:i-d-good} with $j':=j+1$ and as $(x_t)_t$ is nondecreasing,
we obtain 
\[ \lim_{t\to\infty} x_t \geq x_{T_{j+1}-1} \geq x_{T_j-1} + 2^{-d-1} >
   \varphi^\mathrm{rc}_i(T_j) \geq \lim_{n\to\infty} \varphi^\mathrm{rc}_i(n). \]
This is a contradiction to $\lim_{t\to\infty} x_t = \lim_{n\to\infty} \varphi^\mathrm{rc}_i(n)$.
Thus, $\mathcal{R}_{4i}(\xi)$ is true.

Thus, we can apply the third case of Theorem~\ref{theorem:master} and obtain 
the assertion that there exists a nearly computable and left-computable number that is not right-computable.
This proves Corollary~\ref{corollary:nearly-left-not-right}.

\section{A Nearly Computable Number that is not Divergence Bounded Computable}
\label{section:divergence}

\begin{definition}[{Zheng, Rettinger, Gengler~\cite[Definition 5.1]{ZRG2005}}]
\begin{enumerate}
\item
Two pairs $(i,j)$ and $(m,n)$ of natural numbers with $i< j$ and $m< n$ are called
{\em non-overlapping} if either $j\leq m$ or $n\leq i$.
\item
Let $f\colon\IN\to\IN$ be a function.
A real number $\xi$ is called {\em $f$-bounded computable} if there exists a 
computable sequence $(x_n)_n$ of dyadic numbers converging to $\xi$
such that, for every $n\in\IN$, there are at most $f(n)$
pairwise non-overlapping pairs $(i,j)$ of natural numbers with $i<j$ and
$|x_i-x_j|\geq 2^{-n}$.
\item
A real number $\xi$ is called {\em divergence bounded computable} (or {\em d.b.c.})
if there exists a total computable function $f\colon\IN\to\IN$ such that $\xi$
is $f$-bounded computable.
\end{enumerate}
\end{definition}

A real number $\xi$ is called \emph{computably approximable} if
there exists a computable sequence of dyadic numbers converging to it.
It is our goal to strengthen the following observation.

\begin{proposition}[{Zheng, Rettinger, Gengler~\cite[Proposition 5.2(3)]{ZRG2005}}]
There exists a computably approximable number that is not divergence bounded computable.
\end{proposition}

The following theorem is the main result of this section.
We shall deduce it from Theorem~\ref{theorem:master}.

\begin{theorem}
\label{theorem:divergence}
There exists a nearly computable number that is not divergence bounded computable.
\end{theorem}

\begin{proof}
We consider the following requirements:
we let $\mathcal{R}_{4i+2}(\xi)$ be true for every real number $\xi$
and define $\mathcal{R}_{4i}(\xi)$ as follows.
For $e,k\in\IN$ we define:
\begin{equation}
\label{eq:requirement-not-divergence-bounded-computable}
\begin{array}{l}
    \text{if } \varphi_e \text{ and } \varphi_k \text{ are total and } 
    \lim_{n\to\infty} \nu_\ID\varphi_k(n) = \xi \\
    \text{then there exists a number } d\in\IN \text{ such that} \\
    \text{there are at least } \varphi_e(d+2)+1 \text{ pairwise non-overlapping} \\
    \text{pairs } (i,j) \text{ with } i<j \text{ and } |\nu_\ID \varphi_k(i) - \nu_\ID \varphi_k(j)| \geq 2^{-d-2} .
    \end{array}
\tag{$\mathcal{R}_{4\langle e,k\rangle}$}
\end{equation}
If a real number $\xi$ satisfies all of these requirements $\mathcal{R}_{4\langle e,k\rangle}$
for $e,k\in\IN$ then it is not divergence bounded computable.

We wish to construct a sequence $(x_t)_t$ of dyadic numbers that converges
to a number $\xi$ that satisfies all of these requirements and is nearly computable.
The strategy for satisfying a single requirement $\mathcal{R}_{4\langle e,k\rangle}$ is simple,
once a number $d\in\IN$ is fixed:
we wait for a $t$ such that $\nu_\ID \varphi_k(t)$ is close to the current value $x_{t-1}$
and then we let $x_t$ jump far away from it. This is repeated so that it is done at least $\varphi_e(d+2)+1$ times.
Then either the sequence $(\nu_\ID \varphi_k(t))_t$ does not converge to $\xi$ or it must 
follow the jumps made by $(x_t)_t$ and must itself make
at least $\varphi_e(d+2)+1$ big jumps.
In order to make sure that the sequence $(x_t)_t$ converges, we alternate between jumps to the right and to the left.
We define a suitable function $A:\IN^3\times\ID\to\ID$.
Given $d,t,\ell\in\IN$ and $q\in\ID$ we define $A(d,t,\ell,q)$ as follows.
First, we define numbers $e,k$ by $i=\langle e,k\rangle$.
\begin{itemize}
\item
If $\iota(d)$ is even and for $i:=\iota(d)/2$ and $e,k\in\IN$ defined by $i=\langle e,k\rangle$
the following is true:
$\tilde{\ell}_k(t)>0$ and
$|q - \nu_\ID\varphi_k(\tilde{\ell}_k(t)-1)| < 2^{-d-3}$
and $\varphi_e(d+2)[t]\downarrow$ and
$\ell\leq \varphi_e(d+2)+1$,
then we set
\[ A(d,t,\ell,q):=(-1)^\ell \cdot 2^{-d} . \]
\item
Otherwise we set
\[ A(d,t,\ell,q):=0 . \]
\end{itemize}
This ends the definition of the function $A$.

It is clear that $A$ is computable and satisfies \eqref{equation:range-powers-of-2} and \eqref{eq:signswitching-case}.
So, we are in the sign switching case of Theorem~\ref{theorem:master}.

Let $d\in\IN$, let $(x_t)_t$ be a computable sequence of dyadic numbers,
and let $(T_j)_j$ be a computable, nonempty, finite or infinite, increasing sequence with $T_0>0$.
If $\iota(d)$ is odd, then $A(d,t,\ell,q)=0$ for all $t,\ell\in\IN$ and $q\in\ID$, hence, 
\eqref{equation:cancellation}, \eqref{equation:bound-on-sum}, \eqref{equation:Jm-finite} are trivially true.
Let us consider the case when $\iota(d)$ is even.
Let us assume additionally that the sequence $(T_j)_j$ is $(d,(x_t)_t)$-good.
We claim that \eqref{equation:bound-on-sum} and \eqref{equation:bounded} are true.
We set $i:=\iota(d)/2$.

First we show \eqref{equation:bounded}.
We wish to show that the function $h_{d}$ is bounded.
This is clear if the sequence $(T_j)_j$ is finite. Let us assume that the sequence $(T_j)_j$ is infinite.
We claim that $h_{d}$ is bounded by $\varphi_e(d+2)+2$.
Let us assume that this is not the case.
It is clear that for every $j$ we have $h_{d}(j+1)\in\{h_{d}(j),h_{d}(j)+1\}$.
Thus, if $h_{d}$ is not bounded by $\varphi_e(d+2)+2$, then there exists a number $j$ such that $h_{d}(j)=\varphi_e(d+2)+2$ and $h_{d}(j+1)=\varphi_e(d+2)+3$.
This implies $0 \neq A(d,T_j,h_{d}(j),x_{T_j-1}) = A(d,T_j,\varphi_e(d+2)+2,x_{T_j-1})$.
But this contradicts the definition of the function $A$:
when the third input parameter is greater than $\varphi_e(d+2)+1$ then it should produce the value $0$.
We have shown that \eqref{equation:bounded} is satisfied.

Next, we wish to show that \eqref{equation:bound-on-sum} is satisfied.
Let us additionally fix some $J\in\dom(T)$.
If $A(d,T_j,h_{d}(j),x_{T_j-1})=0$ for all $j\leq J$ or if $A(d,T_j,h_{d}(j),x_{T_j-1})\neq 0$ for just one $j \leq J$ then
\eqref{equation:bound-on-sum} is obviously true.
If there are at least two numbers $j,j'$ with $j,j' \leq J$
and with $A(d,T_j,h_{d}(j),x_{T_j-1})\neq 0$ and $A(d,T_{j'},h_{d}(j'),x_{T_{j'}-1})\neq 0$
and if $j'$ is the smallest number larger than $j$ with $A(d,T_{j'},h_{i,d}(j'),x_{T_{j'}-1})\neq 0$
then $h_{d}(j')=h_{d}(j)+1$, hence, 
\[ A(d,T_{j'},h_{d}(j'),x_{T_{j'}-1})= -A(d,T_j,h_{d}(j),x_{T_j-1}) , \]
that is, consecutive nonzero values cancel each other out.
Thus, in this case, 
\[    \sum_{j=0}^{J} A(d,T_j,h_{d}(j),x_{T_j-1}) \in \{0,2^{-d}\},  \]
and \eqref{equation:bound-on-sum} is true again.

Remember that $\mathcal{R}_{4i+2}(\xi)$ is true for every real number $\xi$.
In order to show that the remaining assumptions of Theorem~\ref{theorem:master} are
satisfied in this case, it is sufficient to consider an even number $2i\in\IN$.
Let us assume that $(x_t)_t$ is a converging, and computable
sequence, that $d\in\IN$ is a number such that $\iota(d)=2i$,
and that $(T_j)_j$ is a computable, increasing, infinite, $(d,(x_t)_t)$-good sequence with $T_0>0$.
We claim that $\mathcal{R}_{4i}(\lim_{t\to\infty} x_t)$ is true.
We set $\xi:=\lim_{t\to\infty} x_t$, and we define $e,k$ by $\langle e,k \rangle = i$.
So, we wish to prove that $\mathcal{R}_{4\langle e,k\rangle}(\xi)$ is true.
We can assume that
$\varphi_e$ and $\varphi_k$ are total and $\lim_{n\to\infty} \nu_\ID(\varphi_k(n))=\xi$
because otherwise $\mathcal{R}_{4\langle e,k\rangle}(\xi)$ is true trivially.

We claim that it is sufficient to show that there exists a $j$ with $h_{d}(j)\geq\varphi_e(d+2)+2$.
Let us assume that such a $j$ exists.
Then let $j_0,\ldots,j_{\varphi_e(d+2)+1}$ with $j_0<j_1< \ldots < j_{\varphi_e(d+2)+1}$
be the $\varphi_e(d+2)+2$ smallest numbers $j$ with $A(d,T_j,h_{d}(j),x_{T_j-1})\neq 0$.
For any $j\in\{j_0,j_1,\ldots,j_{\varphi_e(d+2)}\}$,
by using \eqref{equation:i-d-good} with $j':=j+1$ we obtain
\[ \left| x_{T_{j+1}-1} - x_{T_j-1} - A(d,T_j,h_{d}(j),x_{T_j-1}) \right|
   \leq \frac{1}{2} \cdot \left|A(d,T_j,h_{d}(j),x_{T_j-1})\right| . \]
Using additionally the fact that $A(d,T_j,h_{d}(j),x_{T_j-1}) \in \{-2^{-d},2^{-d}\}$
we obtain
\[ \left| x_{T_{j+1}-1} - x_{T_j-1} \right| \geq 2^{-d-1} . \]
We also know that 
$\left|x_{T_j-1} - \nu_\ID\varphi_k(\tilde{\ell}_k(T_j)-1)\right| < 2^{-d-3}$ for
$j\in \{j_0,j_1,\ldots,j_{\varphi_e(d+2)+1}\}$.
Hence, for $j\in\{j_0,j_1,\ldots,j_{\varphi_e(d+2)}\}$ we obtain
\[ \left| \nu_\ID\varphi_k(\tilde{\ell}_k(T_{j+1})-1) - \nu_\ID\varphi_k(\tilde{\ell}_k(T_j)-1) \right| > 2^{-d-2} . \]
Note also that the $\varphi_e(d+2)+1$ pairs 
$(\tilde{\ell}_k(T_j)-1,\tilde{\ell}_k(T_{j+1})-1)$ for $j\in\{j_0,j_1,\ldots,j_{\varphi_e(d+2)}\}$
are pairwise nonoverlapping and satisfy $\tilde{\ell}_k(T_j)-1 < \tilde{\ell}_k(T_{j+1})-1$.
Thus, we have shown that $\mathcal{R}_{4\langle e,k\rangle}(\xi)$ is true
if there exists a $j$ with $h_{i,d}(j)\geq\varphi_e(d+2)+2$.

For the sake of a contradiction, let us assume that there is no such $j$.
Then the function $h_{d}$ is bounded by $\varphi_e(d+2)+1$.
As the sequence $(h_{d}(j))_j$ is nondecreasing, it converges.
Let $H:=\lim_{j\to\infty} h_{d}(j)$. Then $H\leq \varphi_e(d+2)+1$.
By choosing a $j$ sufficiently large we obtain
\begin{itemize}
\item
$\tilde{\ell}_k(T_j)>0$,
\item
$|x_{T_j-1} - \nu_\ID\varphi_k(\tilde{\ell}_k(T_j)-1)| < 2^{-d-3}$
(remember that we assume $\lim_{n\to\infty} \nu_\ID(\varphi_k(n))=\xi=\lim_{t\to\infty} x_t$),
\item
$\varphi_e(d+2)[T_j]\downarrow$, and
\item
$h_{d}(j)=H\leq\varphi_e(d+2)+1$.
\end{itemize}
But that implies that $A(d,T_j,h_{d}(T_j),x_{T_j-1})\neq 0$, hence,
$h_{d}(j+1) = h_{d}(j)+1 = H+1$ in contradiction to the definition of $H$.
We have shown that the function $A$ leads to the fulfillment of the requirement $\mathcal{R}_{4i}$.

Thus, we can apply the first case of Theorem~\ref{theorem:master} and obtain the assertion of the theorem.
\end{proof}

\section{Solovay Reducibility}
\label{section:Solovay-reducibility}

For left-computable real numbers there is a natural reducibility relation due to Solovay~\cite{Sol1975}, which leads to a classification of left-computable real numbers according to the question
how easy it is to approximate them from the left effectively. 
In this section we present the definition and some basic, well known properties of Solovay reducibility.

\begin{definition}[Solovay~\cite{Sol1975}]
Let $\alpha$ and $\beta$ be left-computable numbers.
We say that $\alpha$ is \emph{Solovay-reducible} to $\beta$ (denoted by $\alpha \leqS \beta$) if
there exist a computable function $g:\subseteq\ID\to\ID$ 
and a number $C\in\IN$ such that
$\{q\in\ID \colon q<\beta\} \subseteq\dom(g)$ and,
for all $q\in\ID$ with $q < \beta$, we have
$g(q) < \alpha$ and $\alpha - g(q) \leq C\cdot \left(\beta - q\right)$.
\end{definition}

We shall use the following characterizations of Solovay reducibility.

\begin{lemma}[{\cite{CHKW2001,DHN2002}}]
\label{lemma:Sreducibility}
Let $\alpha$ and $\beta$ be left-computable numbers.
Let $(a_n)_n$ and $(b_n)_n$ be computable, nondecreasing, but not ultimately constant
sequences of dyadic numbers converging to $\alpha$ and $\beta$, respectively.
The following conditions are equivalent.
\begin{enumerate}
\item
$\alpha\leqS\beta$.
\item
There exist a number $C\in\IN$ and a total computable function $f\colon\IN\to\IN$ such that, for all $n\in\IN$,
\[ \alpha - a_{f(n)} < C \cdot (\beta - b_n) . \]
\item
There exist computable sequences $(c_n)_n$ and $(d_n)_n$ of dyadic numbers 
with $\sum_{n=0}^\infty c_n=\alpha$ and $\sum_{n=0}^\infty d_n=\beta$
and a constant $C\in\IN$ such that, for all $n\geq 1$, 
\[ 0 \leq c_n \leq C \cdot d_n . \]
\item
There exists a number $D\in\IN$ such that the number $D\cdot\beta-\alpha$ is left-computable.
\end{enumerate}
\end{lemma}

It is obvious that the Solovay reducibility relation is reflexive and transitive.
As usual, one defines an equivalence relation $\equivS$ on the left-computable numbers
by $\alpha \equivS \beta$ if, and only if, $\alpha \leqS \beta$ and $\beta \leqS \alpha$.
Its equivalence classes are called \emph{Solovay degrees}.
There is a least Solovay degree, the set of computable numbers~\cite{CHKW2001}.
And there is a highest Solovay degree whose elements can be described in several different ways.

\begin{proposition}[Chaitin~\cite{Cha1975}, Solovay~\cite{Sol1975},
Calude, Hertling, Khoussainov, Wang~\cite{CHKW2001},
Ku\v{c}era, Slaman~\cite{KS2001}]
\label{prop:omega}
For a left-computable number $\beta$ the following are equivalent:
\begin{enumerate}
\item
For all left-computable numbers $\alpha$ we have $\alpha \leqS \beta$.
\item
The fractional part of $\beta$ (that is, the unique number $\gamma\in[0,1)$ with $\beta-\gamma\in\IZ$) is an Omega number in the sense of Chaitin~\cite{Cha1975}.
\item
$\beta$ is Martin-L{\"o}f random~\cite{Mar1966}.
\end{enumerate}
\end{proposition}

Furthermore $\leqS$ is an upper semilattice:
If $\alpha$ and $\beta$ are left-computable numbers then
the Solovay degree of the left-computable number $\alpha+\beta$ is the $\leqS$-supremum
of the Solovay degrees of $\alpha$ and $\beta$~\cite{CHKW2001}.
For further results on Solovay degrees see~\cite{DHN2002} and the monographs~\cite{Nie2009,DH2010}.

\section{Nearly Computable Numbers and Solovay Incomparability}
\label{section:Solovay-incomparable}

The following theorem is the main result of this section.
We shall deduce it from Theorem~\ref{theorem:master}.
By $\leqS$ we denote the Solovay reducibility relation between left-computable numbers.

\begin{theorem}
\label{theorem:technical-incomparable-Sdegrees}
Let $\alpha$ be a left-computable number that is not computable,
and let $\beta$ be a left-computable number that is not Martin-L{\"o}f random.
Then there exists a nearly computable and left-computable number $\xi$
with $\alpha\not\leqS \xi$ and $\xi \not\leqS \beta$.
\end{theorem}

The motivation for this result is the following corollary.

\begin{corollary}
\label{corollary:incomparable-Sdegrees}
Let $\alpha$ be a left-computable number that 
is neither in the least nor in the highest Solovay degree of left-computable numbers.
Then there exists a left-computable and nearly computable number $\xi$ such that 
$\alpha\not\leqS\xi$ and $\xi\not\leqS\alpha$.
\end{corollary}

So, the number $\xi$ is Solovay-incomparable with $\alpha$.
To the best of our knowledge, even the assertion of this corollary
without the statement that $\xi$ is nearly computable, is new.

\begin{proof}
Apply Theorem~\ref{theorem:technical-incomparable-Sdegrees}
with $\beta:= \alpha$.
\end{proof}

Let us discuss one interesting consequence of Corollary~\ref{corollary:incomparable-Sdegrees}.
It is well known that there is a highest Solovay degree of left-computable numbers
and that its elements are exactly the left-computable, Martin-L{\"o}f random numbers; compare Proposition~\ref{prop:omega}.
It is also known that the left-computable nearly computable numbers are closed downwards
under Solovay reducibility~\cite[Prop. 7.5]{HJ2024}.
In particular, any Solovay degree of left-computable numbers that contains a left-computable and nearly computable number
contains only left-computable, nearly computable numbers.
This raises the question whether similarly to the maximal Solovay degree of all left-computable numbers
there is perhaps a maximal Solovay degree of nearly computable, left-computable numbers.
The answer to this question is no.

\begin{corollary}
There is no left-computable, nearly computable number $\beta$ such that all left-computable, nearly computable numbers $\alpha$ satisfy $\alpha \leq_{S} \beta$.
\end{corollary}

\begin{proof}
For the sake of a contradiction, let us assume that there is a left-computable, nearly computable number $\beta$ such that all left-computable, nearly computable numbers $\alpha$ satisfy $\alpha \leq_{S} \beta$.
Such a left-computable number $\beta$ cannot be an element of the maximal Solovay degree of left-computable numbers
because the maximal Solovay degree does not contain a nearly computable number because the left-computable, nearly computable numbers are not Martin-Löf random; see~\cite[Theorem 6]{SW2005}.
And such a left-computable number $\beta$ cannot be an element of the smallest Solovay degree of left-computable numbers
because it is known that there exist left-computable, nearly computable numbers that are not in the smallest Solovay degree; 
see Corollary~\ref{corollary:nearly-left-not-right} or \cite[Theorem 7.3]{HJ2024}.
Finally, due to Corollary~\ref{corollary:incomparable-Sdegrees} such a left-computable number $\beta$ cannot be an element of any other Solovay degree.
\end{proof}

Theorem~\ref{theorem:technical-incomparable-Sdegrees}
will be deduced from Theorem~\ref{theorem:master}.

\begin{proof}[Proof of Theorem~\ref{theorem:technical-incomparable-Sdegrees}]
Let $(\varphi_i^\mathrm{lc})_i$ be an effective list of all left-computable
real numbers. By that we mean, $\varphi^\mathrm{lc}\colon\IN^2\to\ID$ is a computable function
such that
\begin{itemize}
\item
for every $i\in\IN$ the sequence $(\varphi^\mathrm{lc}_i(n))_n$ is increasing and converging,
\item
for every left-computable number $\xi$ there exists an $i\in\IN$ such that
$\xi= \lim_{n\to\infty} \varphi^\mathrm{lc}_i(n)$.
\end{itemize}
We set $\varphi^\mathrm{lc}_i(\infty) :=  \lim_{n\to\infty} \varphi^\mathrm{lc}_i(n)$.

We consider the following requirements:
for $i\in\IN$ 
we let $\mathcal{R}_{4i+2}(\xi)$ be true for every real number $\xi$
and define $\mathcal{R}_{4i}(\xi)$ as follows.
For $\xi\in\IR$ and $e,k\in\IN$ we define:
\begin{align}
& 2^{-e} \cdot \xi + \varphi^\mathrm{lc}_k(\infty) \neq \beta
\label{eq:requirement_xi-notleqS-beta}
\tag{$\mathcal{R}_{8\langle e,k\rangle}$} \\
& 2^{-e} \cdot \alpha + \varphi^\mathrm{lc}_k(\infty) \neq \xi
\label{eq:requirement_alpha-notleqS-xi}
\tag{$\mathcal{R}_{8\langle e,k\rangle+4}$} 
\end{align}
If a real number $\xi$ satisfies all of these requirements then by Lemma~\ref{lemma:Sreducibility} $\xi \not\leqS\beta$ and $\alpha\not\leqS\xi$.

Let $(a_t)_t$ be a computable, increasing sequence of dyadic numbers converging to $\alpha$.
Let $(b_t)_t$ be a computable, increasing sequence of dyadic numbers converging to $\beta$.
Let $\Omega$ be a left-computable and Martin-L{\"o}f random number with $0<\Omega<1$,
and let $(\omega_n)_n$ be a nondecreasing, computable sequence of natural numbers such that $\sum_{n=0}^\infty 2^{-\omega_n}=\Omega$.

The strategy for satisfying the requirements $\mathcal{R}_{8\langle e,k\rangle}$ is as follows:
whenever the value $2^{-e} \cdot x_{t-1} + \varphi^\mathrm{lc}_k(t)$ is close to $b_t$ (where $x_{t-1}$ is our current approximation of the 
left-computable number $\xi$ that we wish to construct) then we try to add a certain fraction of $\Omega$ (with a scaling factor $2^{-d}$)
to $x_{t-1}$. If the requirement $\mathcal{R}_{8\langle e,k\rangle}(\xi)$ were not true then,
on the one hand, we would have $\xi \leqS \beta$,
on the other hand, this would happen infinitely often,
and this would have the effect that the number $\xi$ would be in the same Solovay degree as $\Omega$,
which is in contradiction to $\xi \leqS \beta$ and to the fact that $\beta$ is not in the same Solovay degree as $\Omega$.
The strategy for satisfying the requirements $\mathcal{R}_{8\langle e,k\rangle+4}$ is as follows:
whenever the value $2^{-e} \cdot a_t + \varphi^\mathrm{lc}_k(t)$ is close to $x_{t-1}$ (where $x_{t-1}$ is our current approximation of the 
left-computable number $\xi$ that we wish to construct) then we add only a very small number to $x_{t-1}$.
If the requirement $\mathcal{R}_{8\langle e,k\rangle}(\xi)$ were not true then,
on the one hand, we would have $\alpha \leqS \xi$,
on the other hand, this would happen infinitely often,
and this would have the effect that there is a subsequence of the sequence $(x_t)_t$ that converges computably,
hence, the limit $\xi=\lim_{t\to\infty} x_t$ would be a computable number,
which is in contradiction to $\alpha \leqS \xi$ and to the fact that $\alpha$ is not computable.

We define a suitable function $A:\IN^3\times\ID\to\ID$.
Given $d,t,\ell\in\IN$ and $q\in\ID$ we define $A(d,t,\ell,q)$ as follows.
\begin{itemize}
\item
If $\iota(d)$ is odd then we define $A(d,t,\ell,q):=0$.
\end{itemize}
In the following we assume that $\iota(d)$ is even.
We define $i:=\iota(d)/2$, and we define numbers $e,k$ by 
\[ \langle e,k\rangle = \begin{cases}
   i/2 & \text{if } i \text{ is even}, \\
	(i-1)/2 & \text{if } i \text{ is odd}.
\end{cases} \]
\begin{itemize}
\item
If $i$ is even then we define 
\[ A(d,t,\ell,q):= \begin{cases}
   2^{-d-\omega_\ell} & \text{if } \left|2^{-e} \cdot q + \varphi^\mathrm{lc}_k(t) - b_t\right| \leq 2^{-\ell} , \\
	 0 & \text{otherwise}.
	\end{cases}
\]
\item
If $i$ is odd then we define
\[ A(d,t,\ell,q):= \begin{cases}
   2^{-d-\ell-1} & \text{if } \left|2^{-e} \cdot a_t + \varphi^\mathrm{lc}_k(t) - q\right| \leq 2^{-\ell} , \\
	 0 & \text{otherwise}.
	\end{cases}
\]
\end{itemize}
This ends the definition of the function $A$.

It is clear that $A$ is computable and satisfies \eqref{equation:range-powers-of-2} and \eqref{eq:nonnegative-case}.
So, we are in the nonnegative case, the third case, of Theorem~\ref{theorem:master}.

Let $d\in\IN$, let $(x_t)_t$ be a computable sequence of dyadic numbers,
and let $(T_j)_j$ be a computable, nonempty, finite or infinite, increasing sequence with $T_0>0$.
If $\iota(d)$ is odd, then $A(d,t,\ell,q)=0$ for all $t,\ell\in\IN$ and $q\in\ID$.
Let us consider the case when $\iota(d)$ is even.
Let us assume additionally that the sequence $(T_j)_j$ is $(d,(x_t)_t)$-good.
We claim that \eqref{equation:bound-on-sum} and \eqref{equation:bounded} are true.
We set $i:=\iota(d)/2$.
As above, we define $e,k\in\IN$ by
\[ \langle e,k\rangle = \begin{cases}
   i/2 & \text{if } i \text{ is even}, \\
	(i-1)/2 & \text{if } i \text{ is odd}.
\end{cases} \]

First we show \eqref{equation:bound-on-sum}.
Let us additionally consider some $J\in\dom(T)$.
In the case of odd $i$ we see
\[ \sum_{j=0}^{J} A(d,T_j,h_{d}(j),x_{T_j-1}) \leq \sum_{\ell=0}^\infty 2^{-d-\ell-1} = 2^{-d} . \]
In the case of even $i$ we see
\[ \sum_{j=0}^{J} A(d,T_j,h_{d}(j),x_{T_j-1}) \leq \sum_{\ell=0}^\infty 2^{-d-\omega_\ell} = 2^{-d} \cdot \Omega < 2^{-d}. \]
Thus, \eqref{equation:bound-on-sum} is true.

Next, we show that \eqref{equation:bounded} is satisfied as well.
We wish to show that the function $h_{d}$ is bounded.
This is clear if the sequence $(T_j)_j$ is finite. Let us assume that the sequence $(T_j)_j$ is infinite.
We distinguish the cases of odd and of even $i$.

Let us first assume that $i$ is odd. As above, we define $e,k\in\IN$ by $\langle e,k\rangle = (i-1)/2$.
For the sake of a contradiction, let us assume that $h_{d}$ is not bounded.
Then there are infinitely many numbers $j$ with $A(d,T_j,h_{d}(j),x_{T_j-1})\neq 0$, hence, with
\[ \left| 2^{-e} \cdot a_{T_j} + \varphi_k^\mathrm{lc}(T_j) - x_{T_j-1} \right| \leq 2^{-h_{d}(j)} . \]
As the sequences $(a_t)_t$ and $(\varphi_k^\mathrm{lc}(t))_t$ converge and the sequence $(x_t)_t$ is nondecreasing we conclude that the
sequence $(x_t)_t$ must converge as well and 
\[ \xi := \lim_{t\to\infty} x_t = 2^{-e}\cdot\alpha + \varphi_k^\mathrm{lc}(\infty) . \]
This implies $\alpha \leqS \xi$.
Let $(j(m))_m$ be the increasing sequence of the numbers $j$ with $A(d,T_j,h_{d}(j),x_{T_j-1})\neq 0$.
This is a computable sequence.
Then $h_{d}(j(m))=m$, for all $m$, and $A(d,T_j,h_{d}(j(m)),x_{T_{j(m)}-1}) = 2^{-d} \cdot 2^{-h_{d}(j(m))} = 2^{-d-m}$.
For all $m$, by using \eqref{equation:i-d-good} with $j:=j(m)$ and $j':=j(m+1)$ we obtain 
\[ \left| x_{T_{j(m+1)}-1} - x_{T_{j(m)}-1} - 2^{-d-m} \right| \leq \frac{1}{2}\cdot 2^{-d-m} , \]
hence, 
\[ \left| x_{T_{j(m+1)}-1} - x_{T_{j(m)}-1} \right| \leq \frac{3}{2}\cdot 2^{-d-m} . \]
This implies that the limit $\xi=\lim_{m\to\infty} x_{T_{j(m)}-1}$ is a computable number.
As $\alpha$ is a left-computable number that is not computable, we obtain $\alpha \not\leqS \xi$. Contradiction.

Now let us assume that $i$ is even. As above, we define $e,k\in\IN$ by $\langle e,k\rangle = i/2$.
For the sake of a contradiction, let us assume that $h_{d}$ is not bounded.
Then there are infinitely many numbers $j$ with $A(d,T_j,h_{d}(j),x_{T_j-1})\neq 0$, hence, with
\[ \left| 2^{-e} \cdot x_{T_j-1} + \varphi_k^\mathrm{lc}(T_j) - b_{T_j} \right| \leq 2^{-h_{d}(j)} . \]
As the sequences $(b_t)_t$ and $(\varphi_k^\mathrm{lc}(t))_t$ converge and the sequence $(x_t)_t$ is nondecreasing we conclude that the
sequence $(x_t)_t$ must converge as well and 
\[ \xi := \lim_{t\to\infty} x_t = 2^e\cdot\beta - 2^e \cdot \varphi_k^\mathrm{lc}(\infty) . \]
This implies $\xi \leqS \beta$.
Let $(j(m))_m$ be the increasing sequence of the numbers $j$ with $A(d,T_j,h_{d}(j),x_{T_j-1})\neq 0$.
This is a computable sequence.
Then $h_{d}(j(m))=m$, for all $m$, and $A(d,T_{j(m)},h_{d}(j(m)),x_{T_{j(m)}-1}) = 2^{-d-\omega_{h_{d}(j(m))}} = 2^{-d-\omega_m}$.
For all $m$, by using \eqref{equation:i-d-good} with $j:=j(m)$ and $j':=j(m+1)$ we obtain 
\[ \left| x_{T_{j(m+1)}-1} - x_{T_{j(m)}-1} - 2^{-d-\omega_m} \right| \leq \frac{1}{2}\cdot 2^{-d-\omega_m} , \]
hence, 
\[ \left| x_{T_{j(m+1)}-1} - x_{T_{j(m)}-1} \right| \geq \frac{1}{2}\cdot 2^{-d-\omega_m} . \]
By Lemma~\ref{lemma:Sreducibility} this implies $\Omega \leqS \xi$.
But from $\Omega \leqS \xi$ and $\xi \leqS \beta$ we obtain $\Omega\leqS\beta$ in contradiction to the assumption $\Omega\not\leqS\beta$.
We have shown that \eqref{equation:bounded} is satisfied.

Remember that $\mathcal{R}_{4i+2}(\xi)$ is true for every real number $\xi$ and every $i\in\IN$.
In order to show that the remaining assumptions of the nonnegative case, the third case, of Theorem~\ref{theorem:master} are
satisfied, it is sufficient to consider an even number $2i\in\IN$.
Let us assume that $(x_t)_t$ is a nondecreasing, converging, and computable
sequence, that $d\in\IN$ is a number such that $\iota(d)=2i$,
and that $(T_j)_j$ is a computable, infinite, increasing, $(d,(x_t)_t)$-good sequence with $T_0>0$.
We set $\xi:=\lim_{t\to\infty} x_t$.
We claim that $\mathcal{R}_{4i}(\xi)$ is true.
We define $e,k$ by 
\[ \langle e,k\rangle = \begin{cases}
   i/2 & \text{if } i \text{ is even}, \\
	(i-1)/2 & \text{if } i \text{ is odd}.
\end{cases} \]
Above we have shown that the nondecreasing function $h_{d}$ is bounded. Hence, it converges to some number $H:=\lim_{j\to\infty} h_{d}(j)$.
There exists a number $j_0$ such that, for all $j\geq j_0$ we have $h_{d}(j)=H$ and $A(d,T_j,h_{d}(j),x_{T_j-1})=0$.
Now we distinguish the cases of odd and of even $i$.
First we consider the case when $i$ is odd. 
For all $j\geq j_0$ we have $A(d,T_j,h_{d}(j),x_{T_j-1})=0$, hence, 
\[ \left| 2^{-e} \cdot a_{T_j} + \varphi_k^\mathrm{lc}(T_j) - x_{T_j-1} \right| > 2^{-h_{d}(j)} = 2^{-H} . \]
As the sequence $\left(2^{-e} \cdot a_{T_j} + \varphi_k^\mathrm{lc}(T_j)\right)_j$ converges to $2^{-e} \cdot \alpha + \varphi_k^\mathrm{lc}(\infty)$,
we conclude that the sequence $(x_{T_j-1})_j$ does not converge to this number, hence,
$\xi \neq 2^{-e} \cdot \alpha + \varphi_k^\mathrm{lc}(\infty)$.
We have shown that $\mathcal{R}_{8\langle e,k\rangle+4}(\xi)$ is true.
Finally we consider the case when $i$ is even. 
For all $j\geq j_0$ we have $A(d,T_j,h_{d}(j),x_{T_j-1})=0$, hence, 
\[ \left| 2^{-e} \cdot x_{T_j-1} + \varphi_k^\mathrm{lc}(T_j) - b_{T_j} \right| > 2^{-h_{d}(j)} = 2^{-H}. \]
As the sequence $\left(- \varphi_k^\mathrm{lc}(T_j) + b_{T_j}\right)_j$ converges to $- \varphi_k^\mathrm{lc}(\infty) + \beta$,
we conclude that the sequence $(2^{-e}\cdot x_{T_j-1})_j$ does not converge to this number, hence,
$2^{-e} \cdot \xi \neq - \varphi_k^\mathrm{lc}(\infty) + \beta$.
We have shown that $\mathcal{R}_{8\langle e,k\rangle}(\xi)$ is true.

Thus, we can apply the nonnegative case of Theorem~\ref{theorem:master}
and obtain the assertion of the theorem.
\end{proof}

\section{A left-computable and nearly computable number that is not regainingly approximable}
\label{section:not-regainingly-approximable}

\begin{definition}[\cite{HHJ2025}]
A real number $\alpha$ is called \emph{regainingly approximable} if there exists a computable, nondecreasing sequence of dyadic numbers
$(a_n)_n$ converging to $\alpha$ with $\alpha - a_n < 2^{-n}$ for infinitely many $n\in\IN$.
\end{definition}

H{\"o}lzl and Janicki~\cite{HJ2026} proved the following theorem.

\begin{theorem}[{H{\"o}lzl and Janicki~\cite[Theorem 37]{HJ2026}}]
\label{theorem:notregaininglyappr}
There exists a left-computable and nearly computable number that is not regainingly approximable.
\end{theorem}

They also proved that a left-computable and nearly computable number is regainingly approximable if, and only if, it is speedable in the sense of Merkle and Titov~\cite{MT2020}.
Merkle and Titov~\cite{MT2020} had shown that speedable numbers are not Martin-L{\"o}f random and had asked whether the converse is true as well.
As a left-computable and nearly computable number cannot be Martin-L{\"o}f random (see~\cite[Theorem 6]{SW2005}),
Theorem~\ref{theorem:notregaininglyappr} gave a negative answer to this question.

It is the goal of this section to demonstrate that Theorem~\ref{theorem:notregaininglyappr}
can be deduced from the fourth assertion of Theorem~\ref{theorem:master}.

\begin{proof}[{Proof of Theorem~\ref{theorem:notregaininglyappr}}]
Let $(\varphi_i^\mathrm{lc})_i$ be an effective list of all left-computable real numbers
as in the proof of Theorem~\ref{theorem:technical-incomparable-Sdegrees}.
We shall define a computable function $A\colon\IN^3\times\ID\to\ID$
that satisfies the assumptions of the fourth case, the nonnegative and left-good case,
of Theorem~\ref{theorem:master} for the following requirements on a real number $\xi$ for $i\in\IN$.
We let $\mathcal{R}_{4i}(\xi)$ be true for every real number $\xi$ and define $\mathcal{R}_{4i+2}(\xi)$ by:
\begin{equation}
\label{eq:req-theorem-not-regainingly}
\tag{$\mathcal{R}_{4i+2}$}
   \text{if } \varphi^\mathrm{lc}_i(\infty) = \xi \text{ then } (\exists m\in\IN)\,(\forall n\geq m) \, \xi - \varphi^\mathrm{lc}_i(n) \geq 2^{-n} . 
\end{equation}
Clearly, if for a left-computable number $\xi$ the requirements $\mathcal{R}_{4i+2}(\xi)$ are satisfied for all $i\in\IN$,
then $\xi$ is not regainingly approximable.
The strategy for satisfying the requirement $\mathcal{R}_{4i+2}$ is as follows:
for some sufficiently large $d$ we do the following:
whenever $x_{t-1} - 2^{-s} < \varphi^\mathrm{lc}_{i}(s) \leq x_{t-1}$ is true for some $s$ with $d+1<s\leq t$
(where $x_{t-1}$ is our current approximation of the left-computable number $\xi$ that we wish to construct) then 
we try to add $2^{-s}$ to $x_{t-1}$. When this succeeds always, then $\mathcal{R}_{4i+2}$ is satisfied.

We define a suitable function $A\colon\IN^3\times\ID\to\ID$.
Given $d,t,\ell\in\IN$ and $q\in\ID$ we define $A(d,t,\ell,q)$ as follows:
\begin{itemize}
\item
If $\iota(d)$ is even then we set $A(d,t,\ell,q) := 0$.
\item
If $\iota(d)$ is odd then we set $i:=(\iota(d)-1)/2$ and we check whether the set
\[ M := \{s \in \IN \colon d+1 < s \leq t \text{ and } q - 2^{-s} < \varphi^\mathrm{lc}_{i}(s) \leq q \} \]
is empty or not. Then we define
\[ A(d,t,\ell,q) := \begin{cases}
    2^{- \min M} & \text{if } M\neq\emptyset, \\
		0 & \text{if } M = \emptyset.
		\end{cases}
\]
\end{itemize}
This ends the definition of the function $A$.

It is clear that $A$ is computable and satisfies \eqref{equation:range-powers-of-2} and \eqref{eq:nonnegative-case}.
We will show that the assumptions of the fourth case, the nonnegative and left-good case,
of Theorem~\ref{theorem:master} are satisfied.

Let us consider some $d\in\IN$.
Let us consider some computable, nondecreasing sequence $(x_t)_t$ of dyadic numbers,
and some computable, nonempty, finite or infinite, increasing, $(d,(x_t)_t)$-left-good sequence $(T_j)_j$ with $T_0>0$.
If $\iota(d)$ is even, then $A(d,\ldots)=0$ for any choice of the other input parameters.
Therefore, the assertions~\eqref{equation:bound-on-sum} and \eqref{equation:bounded} are trivially true.

Let us assume that $\iota(d)$ is odd and let us set $i:=(\iota(d)-1)/2$.
We claim that \eqref{equation:bound-on-sum} is true.
By definition of $A$, for $s\leq d+1$ there is no $j$ with $A(d,T_j,h_d(j),x_{T_j-1})=2^{-s}$.
We claim that for any $s \geq d+2$ there are at most two numbers $j$ with $A(d,T_j,h_d(j),x_{T_j-1})=2^{-s}$.
Let us assume that for some $s \geq d+2$ there are at least two numbers $j$ with $A(d,T_j,h_d(j),x_{T_j-1})=2^{-s}$,
and let $j$ be the second one. 
This implies 
\[ x_{T_j-1} - 2^{-s} < \varphi^\mathrm{lc}_{i}(s) \leq x_{T_j-1} . \]
If $\dom(T)=\{0,\ldots,j\}$ then the assertion is clear.
Otherwise, as the sequence $(T_k)_k$ is $(d,(x_t)_t)$-left-good,
we obtain 
\[  x_{T_{j+1}-1} - x_{T_j-1} \geq A(d,T_j,h_d(j),x_{T_j-1}) = 2^{-s} , \]
hence, for all $j'\geq j+1$ we obtain
\[ x_{T_{j'-1}} \geq x_{T_{j+1}-1} \geq x_{T_j-1} + 2^{-s} \geq \varphi^\mathrm{lc}_{i}(s) + 2^{-s} . \]
Hence, there can be no $j'>j$ with $A(d,T_{j'},h_d(j'),x_{T_{j'}-1}) = 2^{-s}$.
We have shown that for any $s \geq d+2$ there are at most two numbers $j$ with $A(d,T_j,h_d(j),x_{T_j-1})=2^{-s}$.
For any $J\in\dom(T)$ we obtain
\[ 0 \leq \sum_{j=0}^J A(d,T_j,h_d(j),x_{T_j-1}) \leq 2 \cdot \sum_{s = d+2}^\infty 2^{-s} = 2^{-d} . \]
Thus, \eqref{equation:bound-on-sum} is true.

Now let us assume that the sequence $(T_j)_j$ is infinite.
We wish to show that \eqref{equation:computable-convergence} is true, that is, we wish to show that
the series $\sum_j A(d,T_j,h_d(j),x_{T_j-1})$ is computably Cauchy.
This is clear if there are only finitely many $j$ with $A(d,T_j,h_d(j),x_{T_j-1})\neq 0$.
Let us assume that there are infinitely many $j$ with $A(d,T_j,h_d(j),x_{T_j-1})\neq 0$.
Above we have proved that for any $s$ there can be at most two numbers $j$ with $A(d,T_j,h_d(j),x_{T_j-1})=2^{-s}$.
Therefore, given a number $n$, one can compute a number $J$ such that $A(d,T_J,h_d(J),x_{T_J-1}) \leq 2^{-n-2}$.
Let $s_0$ be the natural number with $A(d,T_J,h_d(J),x_{T_J-1})=2^{-s_0}$.
Then $s_0$ is the smallest number $s'$ with $d+1<s'\leq T_J$ and
\[ x_{T_J-1} - 2^{-s'} < \varphi^\mathrm{lc}_{i}(s') \leq x_{T_J-1} . \]
We claim that for all $j\geq J$ we have $A(d,T_j,h_d(j),x_{T_j-1})\leq 2^{-s_0}$.
For the sake of a contradiction, let us assume that there is a $j>J$ such that $A(d,T_j,h_d(j),x_{T_j-1})> 2^{-s_0}$.
Let $s$ be the natural number with $A(d,T_J,h_d(J),x_{T_J-1})=2^{-s}$.
Then $d+1<s<s_0$ and 
\[ x_{T_j-1} - 2^{-s} < \varphi^\mathrm{lc}_{i}(s) \leq x_{T_j-1} . \]
But $j>J$ implies $x_{T_j-1} \geq x_{T_J-1}$ and $d+1<s<s_0$ implies $\varphi^\mathrm{lc}_{i}(s) \leq \varphi^\mathrm{lc}_{i}(s_0)$,
so we obtain
\[ x_{T_J-1} - 2^{-s}
  \leq x_{T_j-1} - 2^{-s} 
	< \varphi^\mathrm{lc}_{i}(s)
  \leq \varphi^\mathrm{lc}_{i}(s_0) 
	\leq x_{T_J-1} \]
in contradiction to the fact that $s_0$ was the smallest number $s'$ with $d+1<s'\leq T_J$
and 
\[ x_{T_J-1} - 2^{-s'} < \varphi^\mathrm{lc}_{i}(s') \leq x_{T_J-1} . \]
We have shown that for all $j\geq J$ we have $A(d,T_J,h_d(J),x_{T_J-1})\leq 2^{-s_0}$.
Using again that for any $s$ there can be at most two numbers $j$ with $A(d,T_j,h_d(j),x_{T_j-1})=2^{-s}$,
we conclude that for any $j_1,j_2$ with $J\leq j_1 \leq j_2$
\[ 0 \leq \sum_{j\colon j_1 \leq j < j_2} A(d,T_j,h_d(j),x_{T_j-1})
  \leq 2 \cdot \sum_{s\geq s_0} 2^{-s} \leq 4 \cdot 2^{-s_0} \leq 2^{-n} . \] 
Thus, the series $\sum_j A(d,T_j,h_ d(j),x_{T_j-1})$ is computably Cauchy
and \eqref{equation:computable-convergence} is true.

Remember that $\mathcal{R}_{4i}(\xi)$ is true for every real number $\xi$ and every $i\in\IN$.
In order to show that the remaining assumptions of the nonnegative and left-good case,
the fourth case, of Theorem~\ref{theorem:master} are satisfied,
it is sufficient to consider an odd number $2i+1\in\IN$.
Let us assume that $(x_t)_t$ is a nondecreasing, converging, and computable
sequence, that $d\in\IN$ is a number such that $\iota(d)=2i+1$,
and that $(T_j)_j$ is a computable, infinite, increasing, $(d,(x_t)_t)$-left-good sequence with $T_0>0$.
We set $\xi:=\lim_{t\to\infty} x_t$.
We claim that $\mathcal{R}_{4i+2}(\xi)$ is true.
We can assume that $\varphi^\mathrm{lc}_i(\infty) = \xi$
because otherwise the requirement $\mathcal{R}_{4i+2}(\xi)$ is satisfied trivially.
We claim that then $\xi - \varphi^\mathrm{lc}_i(s) \geq 2^{-s}$ is even true for all $s>d+1$.
For the sake of a contradiction let us assume that this is not the case.
Let $s_0$ be the smallest number with $d+1 < s_0$ and $\xi - \varphi^\mathrm{lc}_i(s_0) < 2^{-s_0}$.
As $\varphi^\mathrm{lc}_i(s_0) < \varphi^\mathrm{lc}_i(\infty) = \xi$
and the sequence $(x_{T_j})_j$ converges to $\xi$, for all sufficiently large $j$ we have
$x_{T_j} - 2^{-s_0} < \varphi^\mathrm{lc}_i(s_0) < x_{T_j}$, and also $d+1 < s_0 \leq T_j$, of course.
But it follows from \eqref{equation:computable-convergence} that
$\lim_{j\to\infty} A(d,T_j,h_ d(j),x_{T_j-1}) = 0$, hence, for all sufficiently large $j$ we have
$A(d,T_j,h_ d(j),x_{T_j-1}) < 2^{-s_0}$. 
This is a contradiction.

We have shown that all assumptions of the fourth case, the nonnegative and left-good case, of Theorem~\ref{theorem:master}
are satisfied. Thus, we obtain Theorem~\ref{theorem:notregaininglyappr}.
\end{proof}

\section{A nearly and weakly computable number with only random decompositions}
\label{section:decomp-random}

We are interested in the interplay between weak computability (compare Definition~\ref{definition:wc} and Proposition~\ref{proposition:wc}) and Martin-L{\"o}f randomness.
First, let us state some well known results that we are going to use.

\begin{proposition}
\label{proposition:randomness-sum}
For two left-computable numbers $\alpha$ and $\beta$ the following conditions are equivalent:
\begin{enumerate}
\item
At least one of $\alpha$ and $\beta$ is Martin-L{\"o}f random.
\item
The sum $\alpha+\beta$ is Martin-L{\"o}f random.
\end{enumerate}
\end{proposition}

\begin{proof}
The direction 1 $\Rightarrow$ 2 follows from the well known fact $\alpha\leqS\alpha+\beta$
and $\beta\leqS\alpha+\beta$ and from Proposition~\ref{prop:omega}.
The direction 1 $\Leftarrow$ 2 was proved by Downey, Hirschfeldt, and Nies~\cite[Corollary 3.6]{DHN2002}.
\end{proof}

\begin{lemma}[{Barmpalias and Lewis-Pye~\cite{BL2017}}]
\label{lemma:leftcomp-diff-random}
If $\alpha$ and $\beta$ are left-computable reals and $\alpha$ is Martin-L{\"o}f random but $\beta$ is not
then $\alpha-\beta$ is left-computable and Martin-L{\"o}f random.
\end{lemma}

\begin{proposition}[{Rettinger and Zheng~\cite[Theorem 2.5]{RZ2005}}]
\label{proposition:weaklycomputable-random}
If a weakly computable number is Martin-L{\"o}f random then it is left-computable or right-computable.
\end{proposition}

Miller looked at the interplay between weak computability and Martin-L{\"o}f randomness.
The equivalence of the first and the third condition in the following corollary was observed by him~\cite[Prop.~3.4]{Mil2017}.

\begin{corollary}
\label{corollary:variation-random}
Let $\gamma$ be a weakly computable number. Then the following conditions are equivalent:
\begin{enumerate}
\item
For every computable sequence $(c_n)_n$ of dyadic numbers converging to $\gamma$
and such that the series $\sum_n |c_{n+1} - c_n|$ converges,
the number $\sum_{n=0}^\infty |c_{n+1} - c_n|$ is Martin-L{\"o}f random.
\item
For every pair $\alpha,\beta$ of left-computable numbers with $\gamma=\alpha-\beta$
the number $\alpha+\beta$ is Martin-L{\"o}f random.
\item
For every pair $\alpha,\beta$ of left-computable numbers with $\gamma=\alpha-\beta$
at least one of $\alpha$ and $\beta$ is Martin-L{\"o}f random.
\end{enumerate}
\end{corollary}

For completeness sake we give the proof.

\begin{proof}
2 $\iff$ 3: This follows directly from Proposition~\ref{proposition:randomness-sum}.

2 $\Rightarrow$ 1:
Let $(c_n)_n$ be a computable sequence $(c_n)_n$ of dyadic numbers converging to $\gamma$
and such that the series $\sum_n |c_{n+1} - c_n|$ converges.
We define two sequences $(a_n)_n$ and $(b_n)_n$ of dyadic numbers by
\[ 
a_n := c_0 + \sum_{i=0}^{n-1} (c_{n+1} \dotminus c_{n}) , \quad
b_n := \sum_{i=0}^{n-1} (c_{n} \dotminus c_{n+1})
\]
(where $\dotminus : \IR^2 \to \IR$ is the modified subtraction defined by $x \dotminus y := \max \{0, x - y\}$).
It is clear that both sequences are computable, nondecreasing and bounded, hence, convergent,
and that their limits $\alpha:=\lim_{n\to\infty} a_n$ and $\beta:=\lim_{n\to\infty} b_n$
satisfy $\alpha - \beta=\gamma$ and $\alpha + \beta - c_0 = \sum_{n=0}^\infty |c_{n+1} - c_n|$.
Note that $\alpha+\beta-c_0$ is Martin-L{\"o}f random if, and only if, $\alpha+\beta$ 
is Martin-L{\"o}f random.

1 $\Rightarrow$ 2:
Let $\alpha,\beta$ be left-computable numbers with $\gamma=\alpha-\beta$.
By adding a sufficiently large positive integer to both $\alpha$ and $\beta$ we can assume
that both $\alpha$ and $\beta$ are positive (note that for any real $\delta$ either all numbers $\delta+n$ for $n\in\IZ$
are Martin-L{\"o}f random or none of them is).
Let $(a_n)_n$ and $(b_n)_n$ be computable sequences of positive dyadic numbers with
$\alpha=\sum_{n=0}^\infty a_n$ and $\beta=\sum_{n=0}^\infty b_n$.
We define a sequence $(c_n)_n$ of dyadic numbers by
\[ c_0:=0,\ c_{2n+1} := c_{2n} + a_n, \ c_{2n+2} := c_{2n+1} - b_n . \]
It is clear that this sequence is computable and converges to $\gamma$ and that
the series $\sum_n  |c_{n+1} - c_n|$ converges and that the limit of this series
satisfies $\sum_{n=0}^\infty |c_{n+1} - c_n| = \alpha + \beta$.
\end{proof}

Miller~\cite{Mil2017} called a real number $\gamma$ \emph{variation random} if it satisfies the first (and then all three)
of the conditions in Corollary~\ref{corollary:variation-random}.
One may now ask whether there is a weakly computable real that is variation random but not Martin-L{\"o}f random.

\begin{proposition}[{Miller~\cite[Theorem 3.5]{Mil2017}}]
\label{proposition:Miller}
There exists a weakly computable number that is variation random but not Martin-L{\"o}f random.
\end{proposition}

In fact, such a number can be characterized in an interesting, slightly different way.

\begin{lemma}
\label{lemma:variationrandom-notrandom}
Let $\gamma$ be a weakly computable number. Then the following conditions are equivalent:
\begin{enumerate}
\item
The number $\gamma$ is variation random but not Martin-L{\"o}f random.
\item
For every pair $\alpha,\beta$ of left-computable numbers with $\gamma=\alpha-\beta$
both $\alpha$ and $\beta$ are Martin-L{\"o}f random.
\end{enumerate}
\end{lemma}

\begin{proof}
1 $\Leftarrow$ 2:
Let $\gamma$ satisfy the second condition.
It is clear that $\gamma$ is variation random.
For the sake of a contradiction, let us assume that $\gamma$ is Martin-L{\"o}f random.
By Proposition~\ref{proposition:weaklycomputable-random}
$\gamma$ is left-computable or right-computable.
In the first case we can write it as $\gamma=\gamma - 0$, in the second case we can write it as
$\gamma = 0 - (-\gamma)$. Thus, in both cases we can write $\gamma$ as the difference of two 
left-computable numbers of which only one is Martin-L{\"o}f random, a contradiction.

1 $\Rightarrow$ 2:
Let $\gamma$ satisfy the first condition, that is, $\gamma$ is variation random but not Martin-L{\"o}f random.
Let $\alpha,\beta$ be a pair of left-computable numbers with $\gamma=\alpha-\beta$.
Then $\alpha$ or $\beta$ is Martin-L{\"o}f random.
If $\alpha$ were Martin-L{\"o}f random but $\beta$ not then by Lemma~\ref{lemma:leftcomp-diff-random}
$\gamma=\alpha-\beta$ would be Martin-L{\"o}f random, a contradiction.
If $\beta$ were Martin-L{\"o}f random but $\alpha$ not then by Lemma~\ref{lemma:leftcomp-diff-random}
$-\gamma=\beta-\alpha$ would be Martin-L{\"o}f random, and so would be $\gamma$, a contradiction.
Thus, both $\alpha$ and $\beta$ must be Martin-L{\"o}f random.
\end{proof}

In this section, we are going to deduce the following theorem from 
the second assertion of Theorem~\ref{theorem:master}.

\begin{theorem}
\label{theorem:main3}
There exists a weakly computable number that is nearly computable and variation random.
\end{theorem}

Note that by \cite[Corollary 8.4(2)]{HJ2024} a nearly computable number cannot be Martin-L{\"o}f random.
So, Theorem~\ref{theorem:main3} strengthens \cite[Theorem 3.5]{Mil2017} by Miller
cited above as Proposition~\ref{proposition:Miller}. 
Due to Lemma~\ref{lemma:variationrandom-notrandom},
a number as in Theorem~\ref{theorem:main3} has the property that
for every pair $\alpha,\beta$ of left-computable numbers with $\gamma=\alpha-\beta$
both $\alpha$ and $\beta$ are Martin-L{\"o}f random.

\begin{proof}[Proof of Theorem~\ref{theorem:main3}]
Let us consider a left-computable, Martin-L{\"o}f random real $\Omega$ with $0 < \Omega < 1$.
We shall define a computable function $A\colon\IN^3\times\ID\to\ID$
that satisfies the assumptions of the second case, the sign switching and weakly computable case,
of Theorem~\ref{theorem:master} for the following requirements on a real number $\xi$ for $i\in\IN$.
We let $\mathcal{R}_{4i}(\xi)$ be true for every real number $\xi$ and define $\mathcal{R}_{4i+2}(\xi)$ by:
\begin{equation}
\label{eq:req-theorem3}
\tag{$\mathcal{R}_{4i+2}$}
    \begin{array}{l}
    \text{if } \varphi_i \text{ is total and } \lim_{n\to\infty} \nu_\ID\varphi_i(n) = \xi \\
		\text{and the series } \sum_n |\nu_\ID\varphi_i(n+1) - \nu_\ID\varphi_i(n)| \text{ converges} \\
    \text{then } \Omega\leqS \sum_{n=0}^\infty |\nu_\ID\varphi_i(n+1) - \nu_\ID\varphi_i(n)|.
    \end{array} 
\end{equation}
By Proposition~\ref{prop:omega} this is sufficient.

The strategy for satisfying the requirement $\mathcal{R}_{4i+2}$ is as follows:
whenever $x_{t-1}$ is close to $\nu_\ID\varphi_i(\tilde{\ell}_i(t)-1)$ 
(where $x_{t-1}$ is our current approximation of the left-computable number $\xi$ that we wish to construct) 
then we will in turns add or subtract a certain fraction of $\Omega$ (multiplied with a scaling factor $2^{-d}$) to $x_{t-1}$.
We add and subtract in turns so that the sum converges.
If $\lim_{n\to\infty} \nu_\ID\varphi_i(n) = \xi = \lim_{t\to\infty} x_t$
then the sequence
$(\nu_\ID\varphi_i(n))_n$ must follow the movements of a subsequence of $(x_t)_t$
and we will add or subtract almost all fractions of $\Omega$. 
This will have the effect that $\Omega \leqS \sum_{n=0}^\infty \left| \nu_\ID\varphi_i(n+1) - \nu_\ID\varphi_i(n) \right|$,
hence, $\mathcal{R}_{4i+2}(\xi)$ is satisfied.

Let $(\omega_n)_n$ be a nondecreasing, computable sequence of natural numbers with
$\sum_{n=0}^\infty 2^{-\omega_n}=\Omega$.
As $\Omega<1$ we see that $\omega_0>0$.
Let $(p_n)_n$ be the nondecreasing sequence of natural numbers such that for all $m\in\IN$
\[ |p^{-1}(\{m\})| = \begin{cases}
     0 & \text{if } \omega^{-1}(\{m\}) = \emptyset, \\
		\left|\omega^{-1}(\{m\})\right| + 3 & \text{if } \omega^{-1}(\{m\}) \neq \emptyset \text{ and } 
		\left|\omega^{-1}(\{m\})\right| \text{ is even}, \\
		\left|\omega^{-1}(\{m\})\right| + 2 & \text{if } \omega^{-1}(\{m\}) \neq \emptyset \text{ and } 
		\left|\omega^{-1}(\{m\})\right| \text{ is odd}.
		\end{cases} \]
The sequence $(p_n)_n$ is uniquely determined and computable.
Furthermore, $p_0>0$.
We define a suitable function $A\colon\IN^3\times\ID\to\ID$.
Given $d,t,\ell\in\IN$ and $q\in\ID$ we define $A(d,t,\ell,q)$ as follows.
\begin{itemize}
\item
If $\iota(d)$ is odd and if for $i:=(\iota(d)-1)/2$ we have $\tilde{\ell}_i(t)>0$ and 
\[ \left| q - \nu_\ID\varphi_i(\tilde{\ell}_i(t)-1) \right| < 2^{-5-d-p_\ell} \]
then we set $A(d,t,\ell,q) := (-1)^\ell\cdot 2^{-2-d-p_\ell}$.
\item
Otherwise we set $A(d,t,\ell,q):=0$.
\end{itemize}
This ends the definition of the function $A$.

It is clear that $A$ is computable and satisfies \eqref{equation:range-powers-of-2} and \eqref{eq:signswitching-case}.
We will show that the assumptions of the second case, the sign switching and weakly computable case,
of Theorem~\ref{theorem:master} are satisfied.

Let us consider some $d\in\IN$.
Let us consider some computable sequence $(x_t)_t$ of dyadic numbers,
and some computable, nonempty, finite or infinite, increasing sequence $(T_j)_j$ with $T_0>0$.
If $\iota(d)$ is even, then $A(d,\ldots)=0$ for any choice of the other input parameters.
Therefore, the assertions~\eqref{equation:bound-on-sum}, \eqref{equation:bound-on-absolute-sum}, and \eqref{equation:bounded} are true trivially.
Let us assume that $\iota(d)$ is odd and let us set $i:=(\iota(d)-1)/2$.
We wish to show that the assertion~\eqref{equation:cancellation} is true.
Let us assume that for some $m$ the set $J_m$ is not empty and that
$j$ and $j'$ are two consecutive numbers in $J_m$.
Then $h_d(j')=h_d(j)+1$. By definition of $A$ we obtain
\begin{eqnarray*}
 A(d,T_{j'},h_d(j'),x_{T_{j'}-1})
  &=& (-1)^{h_d(j')}\cdot 2^{-2-d-p_{h_d(j')}} \\
  &=& (-1)^{h_d(j')}\cdot 2^{-m} \\
  &=& (-1)^{h_d(j)+1} \cdot 2^{-2-d-p_{h_d(j)}} \\
	&=& - A(d,T_{j},h_d(j),x_{T_{j}-1}) . 
\end{eqnarray*}
Now let us assume that $m,m'\in\IN$ are numbers with $m<m'$ and $J_{m'}\neq\emptyset$ and $J_m\neq\emptyset$.
Then the set $J_m$ has exactly as many elements as the set $p^{-1}(\{m-2-d\})$,
which has an odd number of elements.
Furthermore, as the sequence $(p_n)_n$ is nondecreasing and by the definition of $A$ we obtain
$\max J_m < \min J_{m'}$.
We have shown that the assertion~\eqref{equation:cancellation} is true.
Actually, by the definition of $(p_n)_n$ and of $A$ it is clear that, for any $m\in\IN$,
the set $J_m$ can contain at most as many elements as the set $p^{-1}(\{m-2-d\})$, which is always finite.
Therefore, \eqref{equation:Jm-finite} is true
(we do not need to assume that the sequence $(T_j)_j$ is $(d,(x_t)_t)$-good).
Next, we wish to show that the assertions~\eqref{equation:bound-on-sum}
and \eqref{equation:bound-on-absolute-sum} are true.
As \eqref{equation:bound-on-absolute-sum} implies \eqref{equation:bound-on-sum}
it is sufficient to prove \eqref{equation:bound-on-absolute-sum}.
By definition of $(\omega_n)_n$, of $(p_n)_n$ and of $A$ we see that
\[ \sum_{j\in\dom(T)} \left|A(d,T_j,h_d(j),x_{T_j-1}) \right| 
    \leq 2^{-d-2} \cdot \left (3 + \Omega\right)
    < 2^{-d}. \]
Thus, \eqref{equation:bound-on-absolute-sum} is true
(again, we do not need to assume that the sequence $(T_j)_j$ is $(d,(x_t)_t)$-good).

Remember that $\mathcal{R}_{4i}(\xi)$ is true for every real number $\xi$ and every $i\in\IN$.
In order to show that the remaining assumptions of the sign switching and weakly computable case,
the second case, of Theorem~\ref{theorem:master} are satisfied,
it is sufficient to consider an odd number $2i+1\in\IN$.
Let us assume that $(x_t)_t$ is a converging and computable
sequence, that $d\in\IN$ is a number such that $\iota(d)=2i+1$,
and that $(T_j)_j$ is a computable, infinite, increasing, $(d,(x_t)_t)$-good sequence with $T_0>0$.
We set $\xi:=\lim_{t\to\infty} x_t$.
We claim that $\mathcal{R}_{4i+2}(\xi)$ is true.
We can assume that $\varphi_i$ is total and $\lim_{n\to\infty} \nu_\ID\varphi_i(n) = \xi$
and that the series $\sum_n |\nu_\ID\varphi_i(n+1) - \nu_\ID\varphi_i(n)|$ converges
because otherwise the requirement $\mathcal{R}_{4i+2}(\xi)$ is satisfied trivially.
We wish to show that 
\[ \Omega\leqS \sum_{n=0}^\infty |\nu_\ID\varphi_i(n+1) - \nu_\ID\varphi_i(n)|. \]

Let us consider an arbitrary number $q\in\range(\omega)$.
Let 
\[ n_0 := \min\{n\in\IN \colon \omega_n=q\} \]
and
\[ k := |\omega^{-1}(\{q\})| . \]
Then the numbers $n_0+i$ for $i=0,\ldots, k-1$
are exactly the numbers $n$ with $\omega_n=q$.
We define $m:=2+d+q$.
Let $(j_e)_e$ be the increasing and infinite sequence of all numbers $j$ with
\[ A(d,T_j,h_{d}(j),x_{T_j-1}) \neq 0 . \]
There is a smallest number $E$ with 
\[ \left| A(d,T_{j_E},h_{d}(j_E),x_{T_{j_E}-1}) \right| = 2^{-2-d-q} = 2^{-m}. \]
Actually, there are exactly 
\[ p^{-1}(\{q\}) = \omega^{-1}(\{q\}) + 3 - \left( \omega^{-1}(\{q\}) \modd 2\right)
   = k + 3 - \left( k \modd 2\right) \]
numbers $j$ with $\left|A(d,T_j,h_{d}(j),x_{T_j-1})\right| = 2^{-m}$, and these are exactly
the numbers $j_{E+h}$ for $h=0,\ldots, k + 2 - \left( k \modd 2\right)$.
For $h=1,\ldots,k+1$ (actually, even for $h=0,\ldots,k + 2 -\left( k \modd 2\right)$, but that is not important here)
we have
\[ \left| x_{T_{j_{E+h}}-1} - \nu_\ID\varphi_i(\tilde{\ell}_i(T_{j_{E+h}})-1) \right| < 2^{-5-d-q} = 2^{-3-m}. \]
And for $h=1,\ldots,k$ we have 
\begin{eqnarray*}
   \lefteqn{ \left| x_{T_{j_{E+h+1}}-1} - x_{T_{j_{E+h}}-1} - (-1)^{h_d(j_{E+h})} \cdot 2^{-2-d-q} \right| } && \\
	&=&  \left| x_{T_{j_{E+h+1}}-1} - x_{T_{j_{E+h}}-1} - A(d,T_{j_{E+h}},h_{d}(j_{E+h}),x_{T_{j_{E+h}}-1}) \right| \\
	&\geq & \frac{1}{2} \cdot \left| A(d,T_{j_{E+h}},h_{d}(j_{E+h}),x_{T_{j_{E+h}}-1}) \right| \\
	&=& 2^{-3-d-q}
\end{eqnarray*}
because the sequence $(T_j)_j$ is $(d,(x_t)_t))$-good.
For $h=1,\ldots,k$ we obtain
\[ \left| x_{T_{j_{E+h+1}}-1} - x_{T_{j_{E+h}}-1} \right| \geq 2^{-3-d-q} . \]
We conclude that for $h=1,\ldots,k$
\begin{eqnarray*}
2^{-4-d-\omega_{n_0+h-1}}
  &=& 2^{-4-d-q} \\
	&\leq & \left| \nu_\ID\varphi_i(\tilde{\ell}_i(T_{j_{E+h+1}})-1)
        - \nu_\ID\varphi_i(\tilde{\ell}_i(T_{j_{E+h}}-1)) \right| \\
  &\leq & \sum_{n=\tilde{\ell}_i(T_{j_{E+h}})-1}^{\tilde{\ell}_i(T_{j_{E+h+1}})-2}
	   \left| \nu_\ID\varphi_i(n+1) - \nu_\ID\varphi_i(n) \right| . 
\end{eqnarray*}
Finally, we obtain
\[ \Omega \equivS 2^{-4-d} \cdot \Omega
  \leqS \sum_{n=0}^\infty \left| \nu_\ID\varphi_i(n+1) - \nu_\ID\varphi_i(n) \right| . \]
We have shown that under the given assumptions,
the requirement $\mathcal{R}_{4i+2}(\xi)$ is satisfied.
Thus, Theorem~\ref{theorem:main3} follows from the second assertion of Theorem~\ref{theorem:master}.
\end{proof}

\section{Nearly and Weakly Computable Numbers versus Differences of Left-computable Nearly Computable Numbers}
\label{section:wnc-dlnc}

First, let us have a look at some well-known computability-theoretic sets of real numbers.
The perhaps most fundamental one of them is the set $\EC$ of computable real numbers.
It is well known that $\EC$ is a real closed field.
Let us recall the definition of a real closed field.

\begin{definition}
A subset $K \subseteq \IR$ is called a \emph{real closed field} if 
$0,1\in K$, if $K$ is closed under addition, subtraction, multiplication, and division,
and if for every $n\in\IN$, every vector $(a_0, \dots, a_n) \in K^{n+1} \setminus\{(0,\ldots,0)\}$, and every $x \in \IR$
with $\sum_{k=0}^{n} a_k x^k = 0$ we have $x \in K$.
\end{definition}

For example, the set $\NC$ of nearly computable numbers is a real closed field as well~\cite{HJ2024}.
We already mentioned that a real number $\alpha$ is computable if, and only if, both $\alpha$ and $-\alpha$ are left-computable.
The set $\LC$ of the left-computable numbers is not a subfield of the field $\IR$ of all real numbers
because there are left-computable numbers that are not computable,
and this implies that $\LC$ is not closed under the operation $\alpha\mapsto -\alpha$.
Let us denote the set of weakly computable numbers by $\WC$.

\begin{proposition}
\begin{enumerate}
\item
The set $\WC$ is the smallest field containing $\LC$
(Ambos-Spies, Weihrauch and Zheng~\cite{AWZ2000})
\item
The set $\WC$ is a real closed field (Raichev~\cite{Rai2005b}).
\end{enumerate}
\end{proposition}

The set of the left-computable numbers that are nearly computable will be denoted by $\LNC := \LC \cap \NC$.
There are numbers in $\LNC$ that are not computable. We already mentioned that this
follows from the main result by Downey and LaForte~\cite{DL2002}
and an observation by Stephan and Wu~\cite[Theorem 5]{SW2005},
and we deduced this from our main result, Theorem~\ref{theorem:master}; see the proof of Corollary~\ref{corollary:nearly-left-not-right}.
For further results on left-computable, nearly computable numbers see \cite{DL2002} and \cite{HJ2024}.
Since the set $\LNC$ is a subset of $\LC$ and has also incomputable members, it does not form a field either.
But there exists a smallest field containing all left-computable nearly computable numbers. 
This is the set
\[ \DLNC := \{ \beta - \gamma \colon \beta,\gamma \in \LNC  \} \]
of all differences of two real numbers that are left-computable and nearly computable.

\begin{theorem}
\label{theorem:DLNC}
The set $\DLNC$ is a real closed field and the smallest subfield of $\IR$ containing $\LNC$.
\end{theorem}

Before we prove this, we show that $\DLNC$ is closed under computable functions that are Lipschitz continuous or differentiable.
Remember that a function $f:D\to\IR$ defined on a subset $D\subseteq\IR^m$ (for some $m\geq 1$) is 
\emph{Lipschitz continuous} if there exists a constant $L \geq 0$ such that
\[ \left| f(x) - f(y) \right| \leq L \cdot \left\| x - y \right\|, \]
for all $x, y \in D$
(here $\|\cdot\|$ is the usual Euclidean norm on $\IR^m$).
For functions $f:D\to\IR$ defined on a subset $D\subseteq\IR^m$ we consider the usual computability notion;
compare~\cite{Wei2000}.

\begin{theorem}
\label{theorem:dlnc-closure}
Let $m\geq 1$, let $D\subseteq \IR^m$ be an open subset, and let $f:D\to\IR$ be a computable function.
If $(\alpha_1,\ldots,\alpha_m) \in D \cap \DLNC^m$ and 
\begin{enumerate}
\item
if $f$ is Lipschitz continuous in a neighbourhood of $(\alpha_1,\ldots,\alpha_m)$, or
\item
if $f$ is differentiable in $(\alpha_1,\ldots,\alpha_m)$,
\end{enumerate} 
then $f(\alpha_1,\ldots,\alpha_m)\in \DLNC$.
\end{theorem}

\begin{proof}
We shall give the proof only in the case $m=1$. The modifications needed for the general case are obvious.
In the proof we use an idea contained in the proof of Theorem 2.5 in \cite{Mil2017}.
Let $\alpha\in D \cap\DLNC$.
Let $f$ be Lipschitz continuous in a neighbourhood of $\alpha$
or differentiable in $\alpha$. In both cases there exist a dyadic number $\varepsilon>0$ and 
a natural number $L$ such that
$\{x \in\IR \colon |\alpha-x|< \varepsilon\}\subseteq D$ and
\begin{equation}
\label{equation:diffLipschitz}
 (\forall x \in\IR) (\text{if } |\alpha - x| < \varepsilon \text{ then } |f(\alpha) - f(x)| < L \cdot |\alpha - x|) . 
\end{equation}
Note also that in both cases $f$ is continuous in $\alpha$.
Let $\beta$ and $\gamma$ be numbers in $\LNC$ with $\alpha=\beta-\gamma$.
Let $(b_n)_n$ and $(c_n)_n$ be computable, increasing sequences converging to $\beta$ and to $\gamma$, respectively.
The sequence $(a_n)_n$ defined by $a_n:=b_n-c_n$, for all $n\in\IN$, converges to $\beta-\gamma=\alpha$.
Let $N$ be so large that $|\alpha - a_m| < \varepsilon$, for all $m\geq N$.
We define a computable and increasing function $g:\IN\to\IN$ as follows.
First, we set $g(0):=N$.
Once $g(n)$ is known, for some $n\in\IN$, we search for a number $m>g(n)$ such that
$|f(a_m) - f(a_{g(n)})|< L \cdot |a_m - a_{g(n)}|$.
Such a number $m$ exists because of \eqref{equation:diffLipschitz} and because $f$ is continuous in $\alpha$,
and one can find such a number $m$ because $f$ is computable.
Let $g(n+1)$ be such a number $m$.
Then, for all $n\in\IN$, we have
\[ |\alpha - a_{g(n)}|<\varepsilon \text{ and } |f(a_{g(n+1}) - f(a_{g(n)})| < L \cdot |a_{g(n+1)} - a_{g(n)}|  . \]
In the following, we will work with the sequences $(a_{g(n)})_n$, $(b_{g(n)})_n$, and $(c_{g(n)})_n$.
For notational simplicity we will rename them and call them
$(a_n)_n$, $(b_n)_n$, and $(c_n)_n$. So,
$(b_n)_n$ and $(c_n)_n$ are computable increasing sequences converging to $\beta$ and to $\gamma$, respectively,
and the sequence $(a_n)_n$ satisfies $a_n=b_n-c_n$ and
\[ |\alpha - a_n|<\varepsilon \text{ and } |f(a_{n+1}) - f(a_n)| < L \cdot |a_{n+1} - a_n|,   \]
for all $n\in\IN$.

We have to show that there exist $\sigma, \tau \in \LNC$ with $f(\alpha) = \sigma - \tau$.
Since both $f$ and $(a_n)_n$ are computable, 
there exists a computable sequence $(r_n)_n$ of dyadic numbers
with $\left|f(a_n) - r_n\right| < 2^{-(n+2)}$, for all $n \in \IN$.
Since $f$ is continuous, we have also $\lim_{n \to \infty} r_n = f(\alpha)$.
We define the following two sequences $(s_n)_n$ and $(t_n)_n$:
\[ 
s_n := r_0 + \sum_{i=0}^{n-1} (r_{n+1} \dotminus r_{n}) , \quad
t_n := \sum_{i=0}^{n-1} (r_{n} \dotminus r_{n+1})
\]
(where $\dotminus : \IR^2 \to \IR$ is the modified subtraction defined by $x \dotminus y := \max \{0, x - y\}$).
Both sequences $(s_n)_n$ and $(t_n)_n$ are computable and non-decreasing.
We wish to show that they are nearly computably Cauchy.
We will only give the proof for the sequence $(s_n)_n$, the other proof is done analogously.
Let $h:\IN \to \IN$ be a computable, increasing function. Then we have:
\begin{align*}
		s_{h(n+1)} - s_{h(n)} &= \sum_{i=h(n)}^{h(n+1)-1} (r_{i+1} \dotminus r_{i}) \\
		&\leq \sum_{i=h(n)}^{h(n+1)-1} \left|r_{i+1} - r_{i}\right| \\
		&= \sum_{i=h(n)}^{h(n+1)-1} |r_{i+1} - f(a_{i+1}) + f(a_{i+1}) - f(a_i) + f(a_i) - r_i| \\
		&< \sum_{i=h(n)}^{h(n+1)-1} \left(2^{-(i+3)} + L \cdot \left|a_{i+1} - a_{i}\right| + 2^{-(i+2)} \right) \\
		&< \sum_{i=h(n)}^{h(n+1)-1} \left(2^{-(i+1)} + L \cdot (b_{i+1} - b_{i} + c_{i+1} - c_i) \right) \\
		&< 2^{-h(n)} + L \cdot (b_{h(n+1)} - b_{h(n)}) + L \cdot (c_{h(n+1)} - c_{h(n)})
\end{align*}
As $\beta$ is nearly computable and $(b_n)_n$ is a computable increasing sequence converging to $\beta$,
by~\cite[Proposition 7.2]{HJ2024} the sequence $(b_n)_n$ is nearly computably Cauchy.
Hence, the sequence $(b_{h(n+1)} - b_{h(n)})_n$ converges computably to zero.
By the same reasoning, the sequence $(c_{h(n+1)} - c_{h(n)})_n$ converges computably to zero as well.
Hence, the sequence $(s_{h(n+1)}-s_{h(n)})_n$ converges computably to zero, too.
We have shown that the sequence $(s_n)_n$ is nearly computably Cauchy.
Due to Lemma \ref{lem:computable-convergent} it converges.
So its limit $\sigma:=\lim_{n\to\infty} s_n$ is a left-computable, nearly computable number.
Similarly, the sequence $(t_n)_n$ is nearly computably Cauchy and convergent,
and its limit $\tau:=\lim_{n\to\infty} t_n$ is a left-computable, nearly computable number.
It is clear that $f(\alpha)=\sigma-\tau$.
Thus, $f(\alpha) \in\DLNC$. This completes the proof.
\end{proof}

We still have to prove Theorem~\ref{theorem:DLNC}.
The following lemma is a slight strengthening of an observation implicitly contained
in the article~\cite{Rai2005b} by Raichev and was explicitly stated in \cite{HJ2024}.

\begin{lemma}[{\cite[Lemma 6.6]{HJ2024}}]
\label{lemma:real-closed-field}
If a subset $K\subseteq\IR$ contains a number $x_0\neq 0$ and is closed under Lipschitz continuous computable functions
$f:\dom(f)\to\IR$ with open domain $\dom(f)\subseteq\IR^m$, $m\geq 1$ arbitrary, then
$K$ is a real closed subfield of $\IR$.
\end{lemma}

\begin{proof}[{Proof of Theorem~\ref{theorem:DLNC}}]
The assertion that $\DLNC$ is a real closed field
follows directly from $1 \in \DLNC$, Lemma~\ref{lemma:real-closed-field},
and Theorem~\ref{theorem:dlnc-closure}.

We claim that $\DLNC$ is also the smallest field containing $\LNC$.
Let $K \subseteq \IR$ be a field with $\LNC \subseteq K$.
Let us consider a number $\alpha \in \DLNC$.
By definition, there are $\beta, \gamma \in \LNC$ with $\alpha = \beta - \gamma$.
Due to our assumption, $\beta$ and $\gamma$ are elements of $K$.
Since $K$ is a field, we obtain $\alpha \in K$. Thus, $\DLNC \subseteq K$.
\end{proof}

Finally we turn to nearly computable numbers that are weakly computable.
Let us denote by $\WNC := \WC \cap \NC$
the set of real numbers that are weakly computable and nearly computable.
Obviously $\DLNC$ is a subset of $\WNC$.
In view of Proposition~\ref{proposition:wc}
one might suspect that the sets $\DLNC$ and $\WNC$ are identical.
But we will show that this is not the case. 

\begin{theorem}
\label{theorem:wnc}
\begin{enumerate}
\item
The set $\WNC$ is a real closed field.
\item
$\DLNC\subsetneq \WNC$.
\item
$\WNC\subsetneq\WC$ and $\WNC\subsetneq \NC$. 
\end{enumerate}
\end{theorem}

\begin{proof}
Since $\WC$~\cite{Rai2005b} and $\NC$~\cite{HJ2024} are both real closed fields,
the set $\WNC$ is a real closed field as well.

It is clear that $\DLNC$ is a subset of $\WNC$.
According to Theorem~\ref{theorem:main3} and Lemma~\ref{lemma:variationrandom-notrandom}
there exists a weakly computable and nearly computable number $\alpha$
such that for all left-computable numbers $\beta$ and $\gamma$ with $\alpha = \beta - \gamma$
both $\beta$ and $\gamma$ are Martin-L{\"o}f random.
On the one hand, we have $\alpha\in\WNC$,
on the other hand, $\alpha\not\in\DLNC$ because a left-computable and Martin-L{\"o}f random number cannot be
nearly computable~\cite[Theorem 6]{SW2005}.

The inclusions $\WNC\subseteq\WC$ and $\WNC\subseteq\NC$ are clear by definition.
Let $K\subseteq\IN$ be a computably enumerable set that is not decidable.
Then the number $0.K:=\sum_{n\in K} 2^{-(n+1)}$ is left-computable, hence also weakly computable.
But it is not computable, hence, by \cite[Pages 105, 106]{DHL2004} (see also \cite[Prop.~7.9]{HJ2024}) it is not nearly computable.
This shows $\WNC\subsetneq\WC$.
According to Theorem~\ref{theorem:divergence} there exists a number that is nearly computable but not divergence bounded
computable. As every weakly computable number is divergence bounded computable~\cite[Proposition 5.2(2)]{ZRG2005},
this shows $\WNC\subsetneq \NC$. 
\end{proof}

\bibliography{nearly-computable}
\bibliographystyle{abbrv}

\end{document}